\documentclass{article}

\usepackage[english]{babel}

\usepackage[scale=0.8]{geometry}

\usepackage{mathtools}
\usepackage{amssymb}
\usepackage{bm}
\usepackage{graphicx}
\usepackage{amsthm}
\usepackage{amsfonts}
\usepackage{algorithm}
\usepackage{algpseudocode}
\usepackage{cases}
\usepackage{booktabs}
\usepackage{siunitx}
\usepackage{subcaption}
\usepackage{authblk}

\usepackage{xcolor}

\definecolor{review1}{RGB}{228,26,28}
\definecolor{review2}{RGB}{0,0,195}
\definecolor{review12}{RGB}{200,30,200}

\newcommand{\ra}{}
\newcommand{\rb}{}
\newcommand{\rab}{}

\usepackage[colorlinks=true, allcolors=blue]{hyperref}

\numberwithin{equation}{section}

\algrenewcommand\algorithmicrequire{\textbf{Input:}}
\algrenewcommand\algorithmicensure{\textbf{Output:}}

\newcommand\restr[2]{{
\left.\kern-\nulldelimiterspace 
#1 
\vphantom{\big|} 
\right|_{#2} 
}}

\makeatletter
\renewcommand\paragraph{%
  \@startsection{paragraph}
    {4}
    {\z@}
    {3.25ex \@plus1ex \@minus.2ex}
    {-1em}
    {\normalfont\normalsize\bfseries\maybe@addperiod}%
}
\newcommand{\maybe@addperiod}[1]{%
  #1\@addpunct{.}%
}
\makeatother

\title{Approximately well-balanced Discontinuous Galerkin methods using bases enriched with Physics-Informed Neural Networks}
\author[1]{Emmanuel Franck}
\author[1]{Victor Michel-Dansac}
\author[1,2]{Laurent Navoret}
\affil[1]{\small Université de Strasbourg, CNRS, Inria, IRMA, F-67000 Strasbourg, France}
\affil[2]{\small IRMA, Université de Strasbourg, CNRS UMR 7501, 7 rue René Descartes, 67084 Strasbourg, France}

\usepackage[capitalize, nameinlink, noabbrev]{cleveref}

\newtheorem{theorem}{Theorem}[section]

\newtheorem{lemma}[theorem]{Lemma}
\newtheorem{proposition}[theorem]{Proposition}

\begin{document}

\maketitle

\begin{abstract}
	This work concerns the enrichment of Discontinuous Galerkin (DG) bases,
so that the resulting scheme provides a much better approximation
of steady solutions to hyperbolic systems of balance laws.
The basis enrichment leverages a prior
-- an approximation of the steady solution --
which we propose to compute using a Physics-Informed Neural Network (PINN).
To that end, after presenting the classical DG scheme,
we show how to enrich its basis with a prior.
Convergence results and error estimates follow,
in which we prove that the basis with prior
does not change the order of convergence,
and that the error constant is improved.
To construct the prior, we elect to use parametric PINNs,
which we introduce, as well as the algorithms to construct a prior from PINNs.
We finally perform several validation experiments
on four different hyperbolic balance laws
to highlight the properties of the scheme.
Namely, we show that the DG scheme with prior
is much more accurate on steady solutions
than the DG scheme without prior,
while retaining the same approximation quality on unsteady solutions.

\end{abstract}

\section{Objectives and model}

In the last decades, much work has been devoted to
proposing numerical methods for hyperbolic systems with source terms,
which correctly capture stationary solutions of the system,
as well as perturbations of flows around these steady states.
If the perturbation is smaller than the scheme error,
traditional numerical schemes are not able
to provide a good approximation of the perturbed steady solution.
To address such an issue,
a first possibility is to refine the mesh in space.
However, for small perturbations,
this would greatly increase the computational overhead.
To avoid this,
schemes specifically dedicated to capturing stationary solutions
have been introduced.
They are called \emph{well-balanced schemes}.

There are two families of well-balanced (WB) schemes:
exactly and approximately WB schemes.
Exactly~WB schemes give an exact representation of the equilibria.
Such schemes are usually developed for subclasses of steady solutions,
especially for complex balance laws, or multidimensional problems.
For instance, first- and second-order accurate exactly WB schemes have been developed for
the shallow water equations
\cite{AudBouBriKlePer2004,KurPet2007,MicBerClaFou2016,MicBerClaFou2017}
or the Euler equations with gravity
\cite{KaeMis2014,ThoPupKli2020}.
High-order exactly well-balanced schemes include
\cite{GalParCas2007,NoeXinShu2007,GabCasDum2018a,MicBerClaFou2020,berberich2020high,BerBulFouMbaMic2021,GomCasParRus2021}
with finite volume methods or related approaches,
{\ra or \cite{XinShu2006,XinZhaShu2010,BriXin2020} with discontinuous Galerkin methods,
		including \cite{ManNoe2021} which relies on rewriting the scheme
		in equilibrium variables
		and \cite{ManOefRic2024} where a global flux strategy is used.}
	{\rb We specifically mention \cite{CasPar2020,GomCasPar2021,GomCasParRus2021,GueEscCas2021},
		where the discrete steady solution is
		found at each time step and in each cell,
		before computing the high-order fluctuations
		around this local steady solution.}
The second family, approximately WB schemes,
consist in ensuring a better approximation of the equilibria
compared to traditional numerical schemes.
This better approximation can be under the form
of a better order of
accuracy~\cite{DesZenBerKli2016,franck2016finite,GomCasPar2021}
or a better error constant \cite{AbbLorBerIapSchWou2020,CiaTorRic2023}.
Both families of WB schemes may incur significant additional computational cost
compared to traditional schemes,
due to the extensive modifications necessary to ensure the WB property,
especially for complex systems and equilibria.

	{\rab
		In this work, we focus on providing a well-balanced scheme for
		the following parametric partial differential equation (PDE):
		\begin{equation}\label{eq:main_model}
			\begin{dcases}
				\partial_t u + \partial_x F_{{\mu}_1}(u) = S_{{\mu}_2}({\ra x,}u), \\
				u(t=0,x)=u_0(x),
			\end{dcases}
		\end{equation}
		with ${\mu}_1$ and ${\mu}_2$
		the parameters of the PDE.
		We set ${\mu} = \{{\mu}_1,{\mu}_2\}$,
		and we assume that
		${\mu} \in \mathbb{P} \subset \mathbb{R}^m$.
		In~\eqref{eq:main_model}, the unknown function is $u$;
		$F_{{\mu}_1}$ is {\ra called} the physical flux function,
		while $S_{{\mu}_2}$ is the source term.
			{\ra We emphasize that~$S_{{\mu}_2}$ may intrinsically depend on space,
				and not only through $u$.}
		We assume that the equation is hyperbolic,
		that is to say that the Jacobian matrix of $F_{{\mu}_1}$
		is diagonalizable with real eigenvalues.
		Our goal will be to construct a
		well-balanced approach for the general steady state
		$\partial_x F_{{\mu}_1}(u)= S_{{\mu}_2}({\ra x,}u)$.

		To that end,
		we endeavor to improve the classical Discontinuous Galerkin (DG) method,
		which usually relies on a discontinuous approximation of the solution
		in a suitable polynomial basis.
		More information on the~DG method can be found
		in~\cite{hesthaven2007nodal,PieErn2012} for instance.
		A natural way of improving the traditional DG method
		to improve the accuracy on some family of solutions
		is to enrich the basis with a prior.
		This is for example the case of the Trefftz
		method~\cite{kretzschmar2014discontinuous,barucq2020,buet2020trefftz,imbert2023space},
		or the non-polynomial bases studied in \cite{YUAN2006295}.

		In the present work, we will consider a modal basis and enrich it with a prior on the steady solution.
		This prior will either be the exact steady solution
		(thus making the scheme exactly WB)
		or an approximation of the steady solution
		(thus making the scheme approximately WB).
		With such an enriched basis,
		the preservation of the exact steady state
		is significantly improved,
		and the method is able to accurately capture
		the dynamics of perturbations around this state.
		However, note that the modal basis under consideration is not orthogonal.
		Hence, the mass matrix of the DG method is block diagonal instead of diagonal.
		In addition, the enrichment do not include any orthogonalization process.
		For a given number of discretization points,
		the scheme is thus slightly less efficient algorithmically than
		DG methods based on orthogonal polynomial bases,
		but we will show that the gains obtained through the enrichment process
		clearly outweigh this effect.

		To perform this enrichment process,
		we require either an exact or an approximate steady state,
		to be evaluated in a preprocessing step, at each Gauss point,
		to compute the DG mass matrix.
		When the DG basis is enriched with the exact steady solution, it turns out that the discretized exact solution is nothing but a discrete steady solution of the scheme,
		which then is exactly WB.
		However, in practice, the exact steady solution is rarely known and we use instead an approximate steady state, which is no longer exactly preserved by the scheme.
		It then becomes only approximately WB.
		However, as evidenced in the numerical experiments,
		even with an approximate steady state,
		the enriched basis is able to provide a very convincing approximation
		of perturbations around the steady solution.
		Of course, this requires the approximation to have a smaller amplitude than that of the perturbation.
		Aside from the basis enrichment,
		the second main contribution of this work is
		to provide a method to build this approximation.

		To get approximate steady states and then perform the basis enrichment,
		we use a learning-based offline computation
		with a neural network to build a prior
		which approximates a parametrized family of equilibria.
		To that end, we use Physics-Informed Neural Networks (PINNs), see e.g.
		\cite{raissi2019physics,cai2021physics},
		and parametric neural networks, see e.g.~\cite{sun2020surrogate}.
		This prior is then introduced into the
		Discontinuous Galerkin basis,
		to increase the accuracy of the scheme
		around this family of equilibria.
		Note that the prior construction could be handled without
		the use of neural networks,
		but we will show that the neural network approach is more efficient.
		Namely, PINNs are particularly well-suited for parametric problems,
		and their mesh-less nature allows for an easy and efficient computation
		of the prior on all the Gauss points of the DG mesh.
		This framework could require significant offline calculation cost
		(depending on the problem),
		but will generate a very small additional cost in the online phase, i.e.,
		when actually using the modified scheme.

		Therefore, the approach proposed in this paper is based
		on the hybridization of classical approaches and neural networks.
		The combination of learning and numerical methods
		(known as Scientific Machine Learning)
		has produced good results for hyperbolic PDEs.
		Examples include work on
		the design of limiters or shock
		detection~\cite{ray2018artificial,beck2020neural,yu2022multi},
		artificial viscosity
		\cite{DisHesRay2020,schwander2021controlling,yu2022data,LeoBois2023},
		or numerical schemes \cite{bar2019learning}.
		More specifically, other learning-based basis enrichment techniques
		have also been successfully implemented
		for other applications.
		In~\cite{ainsworth2021galerkin}, for elliptic problems,
		the authors use a network to provide a finite element basis
		that is dynamically enriched.
		In~\cite{sun2022local,sun2023local}, the authors show that
		random neural networks can be used as DG bases,
		and can be more accurate than classical ones
		for a sufficiently large number of basis function in each cell.
	}

The paper is constructed as follows.
First, we assume that we know a prior (an approximation)
of a family of equilibria,
and we introduce the modification of the DG basis.
Theoretical results show that this modification
does not change the order of accuracy of the method,
but decreases the error constant close to steady solutions.
Then, we introduce the learning methods
that will enable us to build our prior for a family of equilibria,
and finally we perform numerical experiments,
in one and two space dimensions,
on several linear and nonlinear systems of balance laws.
A conclusion ends this paper.

\section{Modified Discontinuous Galerkin scheme}
\label{sec:modified_dg_scheme}

This section is devoted to the presentation of the modified DG scheme.
We start by quickly introducing the classical DG scheme
in \cref{sec:dg_scheme},
and then move on to proposing the modification
in \cref{sec:modal_enriched_dg_scheme}.
Theoretical convergence results related to this modification
will be presented in \cref{sec:theoretical_results}.
In this section and the following one,
we write the scheme in the case of a scalar and one-dimensional PDE,
but the method is easily extendable to systems and to higher dimensions.

\subsection{Classical Discontinuous Galerkin scheme}
\label{sec:dg_scheme}

The goal of this section is to present the classical DG scheme
in order to discretize the PDE~\eqref{eq:main_model}.
To that end, we discretize the space domain
	{\ra $\Omega \subset \mathbb{R}$} in cells
$\Omega_k = (x_{k - 1/2}, x_{k + 1/2})$ of size $\Delta x_k$,
and of centers $x_k$.

The idea behind the classical DG scheme is to
first compute the weak form of the considered PDE,
and then to locally approximate the solution in each cell,
by projecting it onto a finite-dimensional vector space~$V_h$.
We consider a space $V_h$ of dimension $q+1$:
\begin{equation*}
	V_h=\operatorname{Span} \left( \phi_{k,0}, \dots, \phi_{k,q} \right).
\end{equation*}
Note that the space $V_h$ can be different for each cell $k$.

The first assumption of DG scheme is to approximate
the solution $u$ to the PDE, in each cell, with a value in~$V_h$:
\begin{equation*}
	\forall k, \quad \restr{u}{\Omega_k}(t,x) \simeq u_k(t,x).
\end{equation*}
Since $u_k \in V_h$, we can write
\begin{equation}
	\label{eq:dg_representation}
	u_k(t,x) \coloneqq
	\sum_{j=0}^q u_{k,j}(t) \phi_{k,j}(x).
\end{equation}
To obtain the DG scheme,
we first write the weak form of the equation in each cell:
\begin{equation}\label{eq:weak1}
	\int_{\Omega_k}\partial_t u(t,x) \phi(x) \, dx +
	\int_{\Omega_k}\partial_x F_{{\mu}_1}(u(t,x)) \phi(x) \, dx =
	\int_{\Omega_k}S_{{\mu}_2}({\ra x,}u(t,x)) \phi(x) \, dx,
\end{equation}
with $\phi(x)$ a smooth test function.
Performing an integration by parts, the above equation is equivalent to
\begin{equation}\label{eq:weak2}
	\begin{aligned}
		\partial_t \left( \int_{\Omega_k} u(t,x) \, \phi(x) \, dx \right)
		 & - \int_{\Omega_k} F_{{\mu}_1}(u(t,x)) \, \partial_x \phi(x) \, dx
		+ \big[F_{{\mu}_1}(u(t,x)) \, \phi(x)\big]_{x_{k - 1/2}}^{x_{k + 1/2}} \\
		 & = \int_{\Omega_k} S_{{\mu}_2}({\ra x,}u(t,x)) \, \phi(x) \, dx.
	\end{aligned}
\end{equation}

We now plug the DG representation \eqref{eq:dg_representation}
in the weak form \eqref{eq:weak2},
using $\phi_{k,i}$ as test function,
for any $i \in \{0, \dots, q\}$.
\begin{enumerate}
	\item We begin with the first term:
	      \[
		      \int_{\Omega_k} u \, \phi_{k,i} =
		      \sum_{j=0}^q \left( \int_{\Omega_k} u_{k,j}(t) \phi_{k,j}(x) \phi_{k,i}(x) \, dx \right)  =
		      \sum_{j=0}^q u_{k,j}(t) \left( \int_{\Omega_k} \phi_{k,j}(x) \phi_{k,i}(x) \, dx \right).
	      \]
	      To handle the integral in the expression above,
	      we introduce the following quadrature formula,
	      with weights~$w_{k,p}$ and points $x_{k,p}$,
	      valid for any smooth function $\phi$:
	      \begin{equation*}
		      \int_{\Omega_k} \phi(x) \, dx \simeq
		      \sum_{p=1}^{n_q} w_{k,p} \, \phi(x_{k,p}).
	      \end{equation*}
	      We assume that the first and last quadrature points
	      coincide with the cell boundaries,
	      i.e. $x_{k,1} = x_{k - 1/2}$ and $x_{k,n_q} = x_{k + 1/2}$.
	      In practice, we use the well-known
	      Gauss-Lobatto quadrature rule,
	      see e.g. \cite{AbrSte1992} for more information.
	      Equipped with this quadrature formula, we introduce
	      \[
		      M_{k,i,j} =
		      \sum_{p=1}^{n_q}w_{k,p} \, \phi_{k,j}(x_{k,p}) \, \phi_{k,i}(x_{k,p}) \simeq
		      \int_{\Omega_k}\phi_{k,j} \phi_{k,i}
		      ,
	      \]
	      so that the first term of \eqref{eq:weak2} becomes
	      \begin{equation*}
		      \int_{\Omega_k} u(t, x) \phi_{k,i}(x) \, dx \simeq
		      \sum_{j=0}^q M_{k,i,j} u_{k,j}(t) {\ra{}= \mathcal{M}_k u_k(t),}
	      \end{equation*}
	      {\ra where $\mathcal{M}_k = (M_{k,i,j})$  denotes the local mass matrix in the $k$-th element, which is of size $(q+1)\times (q+1)$.}
	\item Using the same techniques,
	      the second term is approximated in the following way:
	      \[
		      \int_{\Omega_k} F_{{\mu}_1}(u(t,x)) \, \partial_x \phi_{k,i}(x) dx \simeq
		      \sum_{p=1}^{n_q} \left[w_{k,p} F_{{\mu}_1}\left(\sum_{j=0}^q u_{k,j}(t) \phi_{k,j}(x_{k,p})\right) \partial_x \phi_{k,i}(x_{k,p})\right] {\ra = \left(\mathcal{V}_{\mu_1, k}(u_k(t))\right)_i,}
	      \]
	      {\ra where, for any vector $u \in \mathbb{R}^q$, $\mathcal{V}_{\mu_1,k}(u) \in \mathbb{R}^q$ denotes the local volume operator in the $k$-th element.}
	\item We note that the third term reduces to
	      \[
		      \big[F_{{\mu}_1}(u) \, \phi_{k,i}\big]_{x_{k - 1/2}}^{x_{k + 1/2}} =
		      F_{{\mu}_1}\left(u_k(t,x_{k+\frac12})\right) \phi_{k,i}(x_{k+\frac12}) -
		      F_{{\mu}_1}\left(u_k(t,x_{k-\frac12})\right) \phi_{k,i}(x_{k-\frac12}),
	      \]
	      where the physical flux $F_{{\mu}_1}$
	      has to be approximated at the cell boundaries.
	      To that end,
	      like the well-known finite volume method,
	      the DG method requires the introduction of
	      a consistent numerical flux
	      \[
		      G_{{\mu}_1}(u_L,u_R) \quad \text{ such that }
		      G_{{\mu}_1}(u,u)=F_{{\mu}_1}(u).
	      \]
	      This numerical flux is then used to approximate
	      the interface flux, as follows
	      \begin{equation*}
		      F_{{\mu}_1}\left(u_k(t,x_{k+\frac12})\right) \simeq
		      G_{{\mu}_1}\left(
		      u_k(t,x_{k+\frac12}),
		      u_{k+1}(t,x_{k+\frac12})
		      \right).
	      \end{equation*}
	      {\ra Therefore, given the unknowns $u_{k-1}(t),u_k(t),u_{k+1}(t) \in \mathbb{R}^q$ in the $(k-1)$-th, $k$-th and $(k+1)$-th cells, we introduce the interface operator between the $k$-th cell and its neighbors:
	      \begin{align*}
		       & \left(\mathcal{I}_{\mu_1,k}(u_{k-1}(t),u_k(t),u_{k+1}(t))\right)_i = G_{{\mu}_1}\left(
		      u_{k,q}(t),
		      u_{k+1,0}(t)
		      \right) \phi_{k,i}(x_{k+\frac12}) -
		      G_{{\mu}_1}\left(
		      u_{k-1,q}(t),
		      u_{k,0}(t)
		      \right) \phi_{k,i}(x_{k-\frac12}).
	      \end{align*}
	      }
	\item Finally, for the last term,
	      we use a straightforward application of the quadrature rule:
	      \[
		      \int_{\Omega_k} S_{{\mu}_2}({\ra x,}u(t,x)) \, \phi(x) \, dx \simeq
		      \sum_{p=1}^{n_q} \left[
		      w_{k,p} S_{{\mu}_2}\left({\ra x_{k,p},}
		      \sum_{j=0}^q u_{k,j}(t) \phi_{k,j}(x_{k,p})
		      \right)\phi_{k,i}(x_{k,p})
		      \right] {\ra = \left(\mathcal{S}_{\mu_1, k}(u_k(t))\right)_i,}
	      \]
	      {\ra where, for any vector $u \in \mathbb{R}^q$, $\mathcal{S}_{\mu_1,k}(u) \in \mathbb{R}^q$ denotes the local source operator in the $k$-th element.}
\end{enumerate}

Gathering all these terms, we show that, in each cell,
the DG scheme can be written as an ordinary differential equation,
where the interface flux term
couples the cell $\Omega_k$ with its neighbors:
\[
	\ra
	\mathcal{M}_k \, \partial_t u_k(t) -
	\mathcal{V}_{{\mu}_1,k}(u_k) +
	\mathcal{I}_{{\mu}_1,k}(u_{k-1},u_k,u_{k+1}) =
	\mathcal{S}_{{\mu}_2,k}( u_k).
\]
Now that we have recalled the classical DG space discretization,
we have all the tools we need to introduce a modification
to this discretization that will enable us
to provide an approximately WB scheme.

\subsection{Enrichment of the modal DG basis}
\label{sec:modal_enriched_dg_scheme}

There are many vector spaces
able to represent the solution in each cell.
For instance, nodal DG schemes~\cite{hesthaven2007nodal}
use Lagrange polynomials or
other polynomials based on nodes chosen within each cell.
Legendre polynomials or Taylor expansions around the cell centers
lead to \emph{modal DG schemes}.
In this work, we focus on the Taylor basis,
given on each cell $\Omega_k$ by
\begin{equation}
	\label{eq:modal_basis}
	V_h =
	\operatorname{Span}
	\left(\phi_{k,0},\phi_{k,1},\phi_{k,2},\dots,\phi_{k,q}\right) =
	\operatorname{Span}
	\left(1, (x-x_k), \frac12(x-x_k)^2,\dots,\frac{1}{q!}(x-x_k)^q\right).
\end{equation}

In the remainder of this section,
we assume that we have access to a prior on the equilibrium,
denoted by~$u_{\theta}(x,{\mu})$.
Obtaining such a prior is discussed in \cref{sec:prior_construction}.
For the moment, suffice it to say that $u_{\theta}$ provides
an approximation of the steady solution for $x \in \Omega$
and for ${\mu}$ in some parameter space $\mathbb{P}$ to be defined.

Given the prior $u_{\theta}$, we modify the local basis $V_h$
to incorporate the prior:
for that, we propose two possibilities.
\begin{itemize}
	\item The \emph{additive correction} $V_h^{+}$
	      consists in replacing the first element of $V_h$ by the prior:
	      \begin{equation}
		      \label{eq:modal_basis_additive}
		      V_h^{+} =
		      \operatorname{Span}
		      \big({\phi}_{k,0}^+,{\phi}_{k,1}^+,{\phi}_{k,2}^+,\dots,{\phi}_{k,q}^+\big) =
		      \operatorname{Span}
		      \left(
		      u_{\theta}(x,{\mu}),
		      (x-x_k),
		      \dots,
		      \frac{1}{q!}(x-x_k)^q
		      \right).
	      \end{equation}
	\item The \emph{multiplicative correction} $V_h^{*}$
	      consists in multiplying each element of $V_h$ by the prior:
	      \begin{equation}
		      \label{eq:modal_basis_multiplicative}
		      V_h^{*} =
		      \operatorname{Span}
		      \big({\phi}_{k,0}^*,{\phi}_{k,1}^*,{\phi}_{k,2}^*,\dots,{\phi}_{k,q}^*\big) =
		      \operatorname{Span}
		      \left(u_{\theta}(x,{\mu}),
		      (x-x_k) \, u_{\theta}(x,{\mu}),
		      \dots,
		      \frac{1}{q!}(x-x_k)^q u_{\theta}(x,{\mu})
		      \right).
	      \end{equation}
\end{itemize}

\noindent A first remark is that,
if the prior is exactly equal to the steady solution,
then it can be exactly represented by an element of
$V_h^{+}$ or $V_h^{*}$ (namely, the first one)
in each cell, which is not the case for the classical space~$V_h$.
However, whether the prior is exact or not,
the method will only be of interest
if the projector onto the modified vector space is accurate
(or even exact in the case of an exact prior).
The second point to note is that, unlike conventional DG approaches,
the bases are not polynomial.
We must therefore ensure that
this does not hinder the convergence of the DG method.
In the next section, we follow Yuan and Shu's work
\cite{YUAN2006295} to study the convergence of the modified DG method,
and provide error estimates.

\section{Error estimates}
\label{sec:theoretical_results}

In this section, we prove some convergence results
on the modified DG scheme.
We assume that our prior $u_{\theta}$ is~$p$ times continuously
differentiable, i.e.,
that it has differentiability class $\mathcal{C}^p$,
with $p \geqslant q+1$.
This hypothesis is compatible with the construction of the prior
from \cref{sec:prior_construction}.

In \cite{YUAN2006295}, the authors study the convergence of
the DG scheme for non-polynomial bases.
They show that,
if the non-polynomial basis can be represented in a specific way
by a polynomial basis, then the convergence of
the local and global projection operators is not hampered.
Using some stability results
(given in \cite{YUAN2006295} for the transport equation)
together with these estimations, convergence can be recovered.

These theoretical results will be split in two parts.
To begin with, in \cref{sec:convergence_non_polynomial},
{\ra we prove} that the bases proposed in \cref{sec:modal_enriched_dg_scheme}
fit into the hypotheses of \cite{YUAN2006295},
which ensures convergence.
However, this study is insufficient to show that the better the prior,
the more accurate the modified DG scheme.
To that end, in \cref{sec:estimate_prior},
we derive the projector estimates
in the case of $V_h^{*}$,
in order to show the potential gains of the method.

\subsection{Convergence in non-polynomial DG bases}
\label{sec:convergence_non_polynomial}

In \cite{YUAN2006295}, the authors prove the following lemma.
\begin{lemma}
	\label{lem:yuan_shu_Taylor_of_basis}
	Consider an approximation vector space $V_h$
	with local basis $(v_{k,0},\dots,v_{k,q})$,
	which may depend on the cell $\Omega_k$.
	If there exists constant real numbers $a_{j\ell}$ and $b_j$
	independent of the size of the cell $\Delta x_k$
	such that, in each cell $\Omega_k$,
	\begin{equation}\label{eq:assumption_shu}
		\forall j \in \{ 0, \dots, q \}, \quad
		\left| v_{k,j}(x) - \sum_{\ell=0}^q a_{j\ell}(x-x_k)^\ell \right| \leq
		b_j (\Delta x_k)^{q+1},
	\end{equation}
	then for any function $u\in H^{q+1}(\Omega_k)$,
	there exists $v_h \in V_h$ and a constant real number $C$
	independent of $\Delta x_k$, such that
	\[
		\| v_h - u \|_{L^\infty(\Omega_k)}
		\leq
		C \|u\|_{H^{q+1}(\Omega_k)}(\Delta x_k)^{q+{\frac12}}.
	\]
\end{lemma}

Using this result, the authors show that
the global projection error in the DG basis converges
with an error in $(\Delta x_k)^{q+1}$ in the Sobolev norm $H^{q+1}$,
and later prove the convergence of the whole scheme
using a monotone flux for a scalar equation.
In the remainder of this section, we prove that the two new bases
proposed in \cref{sec:modal_enriched_dg_scheme}
satisfy the assumptions of \cref{lem:yuan_shu_Taylor_of_basis}.
Using these results together with the proofs of \cite{YUAN2006295},
we will obtain that both bases lead to a convergent scheme.

\begin{proposition}
	\label{prop:Vh_mod_add}
	If the prior $u_{\theta}(x;{\mu})$
	has differentiability class $\mathcal{C}^{q+1}(\mathbb{R})$
	with respect to $x$,
	then the approximation space $V_h^{+}$
	satisfies the assumption \eqref{eq:assumption_shu}.
\end{proposition}

\begin{proof}
	Since the prior is $\mathcal{C}^{q+1}(\mathbb{R})$,
	we can write its Taylor series expansion
	around the cell center $x_k$.
	Namely, there exists a constant
	$c\in [x_{k-1/2},x_{k+1/2}]$ such that
	\begin{equation}
		\label{eq:Taylor_series_prior}
		u_{\theta}(x) =
		u_{\theta}(x_k) +
		(x-x_k) u_{\theta}'(x_k) +
		\dots +
		\frac{1}{q!}(x-x_k)^{q} u^{(q)}(x_k) +
		\frac{(x-x_k)^{q+1}}{(q+1)!} u^{(q+1)}(c).
	\end{equation}
	With that expansion, we can write our basis $V_h^{+}$
	with respect to the classical modal basis $V_h$ as follows:
	\[
		\begin{pmatrix}
			u_{\theta}(x) \\
			(x-x_k)       \\
			\vdots        \\
			(x-x_k)^q
		\end{pmatrix} =
		\underbrace{
			\begin{pmatrix}
				u_{\theta}(x_k)  &
				u_{\theta}'(x_k) &
				\ldots           &
				\frac{1}{q!} u_{\theta}^{(q)}(x_k)          \\
				0                & 1      & \ldots & 0      \\
				\vdots           & \vdots & \ddots & \vdots \\
				0                & 0      & \hdots & 1      \\
			\end{pmatrix}
		}_{A_+}
		\begin{pmatrix}
			1       \\
			(x-x_k) \\
			\vdots  \\
			(x-x_k)^q
		\end{pmatrix} +
		(x-x_k)^{q+1}
		\underbrace{
			\begin{pmatrix}
				\frac{u^{(q+1)}(c)}{(q+1)!} \\
				0                           \\
				\vdots                      \\
				0
			\end{pmatrix}.
		}_{b_+}.
	\]
	{\rb We remark that the matrix $A_+$ is independent of $\Delta x_k$.
	Moreover, the first component of the vector~$b_+$ is bounded by
	$\smash{\frac{1}{(q+1)!}\| u^{(q+1)} \|_{L^\infty(\Omega_k)}}$,
	and its other components are zero.
	Therefore, each component of~$b_+$ is bounded by a constant
	value, independent of $\Delta x_k$.}
	Hence, assumption~\eqref{eq:assumption_shu} is verified,
	and \cref{lem:yuan_shu_Taylor_of_basis} can be applied.
\end{proof}

\begin{proposition}
	\label{prop:Vh_mod_mul}
	If the prior $u_{\theta}(x;{\mu})$
	has differentiability class $\mathcal{C}^{q+1}(\mathbb{R})$
	with respect to $x$,
	then the approximation space $V_h^{*}$
	satisfies the assumption \eqref{eq:assumption_shu}.
\end{proposition}

\begin{proof}
	The proof follows the same lines as the proof of
	the previous proposition.
	Namely, \eqref{eq:Taylor_series_prior} is still satisfied
	since the prior is $\mathcal{C}^{q+1}(\mathbb{R})$.
	Then, the basis $V_h^{*}$ is written
	with respect to the classical modal basis $V_h$ as follows:
	\[
		\renewcommand{\arraystretch}{1.5}
		\begin{pmatrix}
			u_{\theta}(x)              \\
			(x-x_k) \, u_{\theta}(x)   \\
			\vdots                     \\
			(x-x_k)^q \, u_{\theta}(x) \\
		\end{pmatrix} =
		\underbrace{
			\begin{pmatrix}
				u_{\theta}(x_k)  &
				u_{\theta}'(x_k) &
				\hdots           &
				\frac{u^{(q)}(x_k)}{q!}                              \\
				0                &
				u_{\theta}(x_k)  &
				\hdots           &
				\frac{u_{\theta}^{(q-1)}(x_k)}{(q-1)!}               \\
				\vdots           & \vdots & \ddots & \vdots          \\
				0                & 0      & \hdots & u_{\theta}(x_k) \\
			\end{pmatrix}
		}_{A_*}
		\begin{pmatrix}
			1         \\
			(x-x_k)   \\
			\vdots    \\
			(x-x_k)^q \\
		\end{pmatrix} +
		(x-x_k)^{q+1}
		\underbrace{
			\begin{pmatrix}
				\frac{u_{\theta}^{q+1}(c)}{(q+1)!} \\
				\frac{u_{\theta}^{q}(c)}{q!}       \\
				\vdots                             \\
				1
			\end{pmatrix}
		}_{b_*}
	\]
	{\rb Just like before, the matrix $A_*$ is independent of $\Delta x_k$,
	and the vector $b_*$ is bounded by values independent of~$\Delta x_k$.}
	Hence, assumption~\eqref{eq:assumption_shu} is verified,
	and \cref{lem:yuan_shu_Taylor_of_basis} can be applied.
\end{proof}

These two propositions show that,
if the prior is sufficiently smooth,
we can apply the results of~\cite{YUAN2006295},
which shows the convergence of the method.
However, this approach does not give an estimation of
the error with respect to the quality of the prior.
Indeed, we expect the modified DG scheme to be more accurate
when the prior is closer to the solution.
Obtaining such an estimate is the objective of the following section.

\subsection{Estimate with prior dependency}
\label{sec:estimate_prior}

The goal of this section is to refine the error estimates
from \cref{sec:convergence_non_polynomial}
for a specific modified basis.
We consider the case of $V_h^{*}$,
since it is easier to write the projector onto the classical basis.
This will enable us to quantify the gains
that can be expected when using this new basis.
The case of $V_h^{+}$ is more complicated,
since the projector is harder to write.
Nevertheless, we will show in the numerical experiments
from \cref{sec:application_numerical_results}
that both modified bases exhibit similar behavior.

Recall that the basis $V_h^{*}$
is obtained by multiplying each element of $V_h$ by the prior.
	{\rb
		Therefore, its basis functions are given
		for each cell $\Omega_k$ and for $j \in \{0, \dots, q\}$ by
		\begin{equation}
			\label{eq:def_phi_star_wrt_phi}
			{\phi}^*_{k,j} = \phi_{k,j} u_{\theta}.
		\end{equation}}

\begin{lemma}
	\label{lem:local_error_estimate}
	Assume that the prior $u_{\theta}$ satisfies
	\[
		u_{\theta}(x; {\mu})^2 > m^2 > 0,
		\quad \forall x \in \Omega,
		\quad \forall {\mu} \in \mathbb{P}.
	\]
	For a given cell $\Omega_k$,
	for any function $u\in H^{q+1}(\Omega_k)$,
	the $L^2$ projector onto $V_h^{*}$,
	denoted by $P_h$ and such that $P_h(u) \in V_h^{*}$,
	satisfies the inequality
	\[
		\big\| u - P_h(u) \big\|_{L^{\infty}(\Omega_k)}
		\lesssim
		\left| \frac{u (\cdot)}{u_{\theta}(\cdot \, ;{\mu})} \right|_{H^{q+1}(\Omega_k)}
		(\Delta x_k)^{q+\frac12} \,
		\left( 1 + \frac{\left\|u_{\theta}(\cdot \, ;{\mu})^2 \right\|_{L^{\infty}(\Omega_k)}}{m^2} \right)
		\big\| u_{\theta}(\cdot \, ;{\mu}) \big\|_{L^\infty}.
	\]
\end{lemma}

\begin{proof}
	The proof uses a strategy similar to \cite{YUAN2006295}.
	We consider the cell $\Omega_k$.
	For any smooth function $f$ defined on~$\Omega_k$,
	we define, for all $x \in \Omega_k$, the operator $T$ by
	\[
		T (f) (x) =
		\sum_{j=0}^q f^{(j)}(x_k) \frac{1}{j!}(x-x_k)^j
	\]
	and the operator $T_{\theta}$ by
	\begin{equation}
		\label{eq:def_T_theta}
		T_{\theta} (f) (x) =
		\left(
		\sum_{j=0}^q \left(
			\frac{f}{u_{\theta}}
			\right)^{(j)}(x_k;{\mu})
		\frac{1}{j!} (x-x_k)^j
		\right) u_{\theta}(x;{\mu}).
	\end{equation}
	For simplicity, we no longer explicitly write the dependence
	in ${\mu}$ in this proof.
	Let $u\in H^{q+1}(\Omega_k)$.
	Using $T_{\theta}$, we write the following estimation:
	\begin{equation}\label{eq:estimate_priorproof1}
		\| u-P_h(u)\|_{L^{\infty}(\Omega_k)} \leq
		\| u- T_{\theta} (u) \|_{L^{\infty}(\Omega_k)} +
		\| T_{\theta} (u) - P_h(u) \|_{L^{\infty}(\Omega_k)}
		\eqqcolon N_1 + N_2.
	\end{equation}
	To complete the proof,
	we need to estimate both terms $N_1$ and $N_2$.

	We start with the estimation of $N_1$.
	We obtain, according to the relationship
	between $T$ and $T_{\theta}$,
	\begin{equation}
		\label{eq:def_N1}
		N_1 = \| u- T_{\theta} (u) \|_{L^{\infty}(\Omega_k)} =
		\left\|
		\frac u {u_\theta} u_\theta -
		T \left( \frac u {u_\theta} \right) u_\theta
		\right\|_{L^{\infty}(\Omega_k)} \leq
		\left\|
		\frac u {u_\theta} -
		T \left( \frac u {u_\theta} \right)
		\right\|_{L^{\infty}(\Omega_k)}
		\left\| u_\theta \right\|_{L^{\infty}(\Omega_k)}.
	\end{equation}
	We can now use an intermediate result from \cite{YUAN2006295}:
	for all $f$ smooth enough, the Taylor formula and
	the Cauchy-Schwartz inequality,
	followed by a direct computation, gives
	\begin{align}
		\notag
		\| f - T(f) \|_{L^{\infty}(\Omega_k)} & =
		\sup_{x \in \Omega_k} \left|
		\int_{x_k}^x f^{(q+1)}(\xi)\frac{(x-\xi)^q}{q!} d\xi
		\right|                                           \\
		\notag
		                                      & \leq
		\sup_{x\in \Omega_k} \left[\left(
			\int_{x_k}^x \left| f^{(q+1)}(\xi)  \right|^2 d\xi
			\right)^{\frac 1 2} \left(
			\int_{x_k}^x \left| \frac{(x-\xi)^q}{q!} \right|^2 d\xi
		\right)^{\frac 1 2} \right],                      \\
		\label{majo_chiwangshu}
		                                      & \lesssim
		| f |_{H^{q+1}(\Omega_k)} \,
		(\Delta x_k)^{q+\frac12}.
	\end{align}
	Going back to $N_1$ and plugging \eqref{majo_chiwangshu}
	into the estimate \eqref{eq:def_N1}, we obtain
	\begin{equation}
		\label{eq:estimate_priorproof3}
		N_1 \lesssim
		\left| \frac u {u_\theta} \right|_{H^{q+1}(\Omega_k)} \,
		(\Delta x_k)^{q+\frac12}
		\left\| u_\theta \right\|_{L^{\infty}(\Omega_k)}.
	\end{equation}

	Now, we proceed with estimating $N_2$,
	the second term of \eqref{eq:estimate_priorproof1}.
	{\rb
	The $L^2$ projector $P_h$ onto $V_h^{*}$
	is defined by
	\[
		P_h(u)=\sum_{j=0}^q \tilde \alpha_j {\phi}^*_{k,j},
	\]
	with $\tilde \alpha = (\tilde \alpha_j)_{j \in \{0, \dots, q\}}$
	is defined for all $\ell \in \{0, \dots, q\}$ by
	\begin{equation*}
		\sum_{j=0}^q
		\tilde \alpha_j \left\langle
		{\phi}^*_{k, \ell}, {\phi}^*_{k, j}
		\right\rangle_{L^2(\Omega_k)}
		=
		\left\langle
		{\phi}^*_{k, \ell}, u
		\right\rangle_{L^2(\Omega_k)}.
	\end{equation*}
	Normalizing by $\Delta x_k$ and defining
	$\alpha_j = \Delta x_k^j \tilde \alpha_j$, we obtain
	\begin{equation*}
		\sum_{j=0}^q
		\alpha_j \left\langle
		\frac{{\phi}^*_{k,\ell}}{(\Delta x_k)^{\ell}},
		\frac{{\phi}^*_{k,j}}{(\Delta x_k)^{j}}
		\right\rangle_{L^2(\Omega_k)}
		=
		\left\langle
		\frac{{\phi}^*_{k,\ell}}{(\Delta x_k)^{\ell}}, u
		\right\rangle_{L^2(\Omega_k)}.
	\end{equation*}
	In the end, we obtain the following expression for $P_h(u)$:
	\[
		P_h(u)=\sum_{j=0}^q \frac{\alpha_j}{(\Delta x_k)^{j}} {\phi}^*_{k,j},
	\]
	with $\alpha =
		(\alpha_j)_{j \in \{1, \dots, q\}} =
		(M^{\ast})^{-1} b^\ast$, where
	\begin{equation}
		\label{eq:def_M_star_and_b}
		M^\ast_{\ell j} = \int_{\Omega_k}
		\frac{{\phi}^*_{k,\ell}(x)}{(\Delta x_k)^{\ell}}
		\frac{{\phi}^*_{k,j}(x)}{(\Delta x_k)^{j}}
		dx,
		\text{\quad and \quad}
		b^\ast_\ell = \int_{\Omega_k} u(x) \frac{{\phi}^*_{k,\ell}(x)}{(\Delta x_k)^{\ell}} dx.
	\end{equation}
	Note that, for all $x \in \Omega_k$,
	arguing the definitions \eqref{eq:def_phi_star_wrt_phi}
	and \eqref{eq:modal_basis}
	of ${\phi}^*_{k,j}$ and $\phi_{k,j}$, we have
	\[
		P_h(u)(x)
		=
		\sum_{j=0}^q \frac{\alpha_j}{(\Delta x_k)^{j}} {\phi}^*_{k,j}(x)
		=
		\sum_{j=0}^q \frac{\alpha_j}{(\Delta x_k)^{j}} {\phi}_{k,j}(x) u_\theta(x)
		=
		\sum_{j=0}^q \frac{\alpha_j}{(\Delta x_k)^{j}} \frac{1}{j!}(x-x_k)^j u_\theta(x).
	\]
	We are now ready to start estimating $N_2$.
	The definition \eqref{eq:def_T_theta} of $T_{\theta}$ yields
	\begin{equation*}
		T_{\theta} (u) (x) =
		\left(
		\sum_{j=0}^q \left(
		\frac{u}{u_{\theta}}
		\right)^{(j)}\!\!(x_k)\,
		\frac{1}{j!} (x-x_k)^j
		\right) u_{\theta}(x).
	\end{equation*}
	Therefore, $N_2$ satisfies}
	\begin{equation*}
		\begin{aligned}
			N_2 =
			\|
			T_{\theta} (u) - P_h(u)
			\|_{L^\infty(\Omega_k)} & =
			\sup_{x \in \Omega_k} \left| \sum_{j=0}^q \left(
			\left(
			\frac{u}{u_{\theta}}
			\right)^{(j)}\!\!(x_k) - \frac{\alpha_j}{(\Delta x_k)^{j}} \right)
			\frac{1}{j!}(x-x_k)^j
			u_{\theta}(x) \right|                             \\
			                        & \leq
			\sup_{x \in \Omega_k} \left| \sum_{j=0}^q \left(
			(\Delta x_k)^{j}\left(
			\frac{u}{u_{\theta}}
			\right)^{(j)}\!\!(x_k) - \alpha_j \right)
			\frac{1}{j!}\frac{(x-x_k)^j}{(\Delta x_k)^{j}} \right|
			\left\| u_{\theta} \right\|_{L^\infty(\Omega_k)}. \\
		\end{aligned}
	\end{equation*}
	Using the Cauchy-Schwartz inequality on the sum,
	and bounding the resulting polynomial on the cell,
	we obtain the estimate
	\begin{equation}\label{eq:second_estimateproof_1}
		N_2 \lesssim
		\left[
			\sum_{j=0}^q
			\left(
			(\Delta x_k)^{j}
			\left(\frac{u}{u_{\theta}}\right)^{(j)}\!\!(x_k) - \alpha_j
			\right)^2
			\right]^{\frac 1 2}
		\left\| u_{\theta} \right\|_{L^\infty(\Omega_k)} = \left\| \delta - \alpha \right\|_2 	\left\| u_{\theta} \right\|_{L^\infty(\Omega_k)},
	\end{equation}
	where
	the vector $\delta =  (\delta_j)_{j \in \{0, \dots, q\}}$ is defined by
	\begin{equation*}
		\delta_j =
		(\Delta x_k)^{j}
		\left( \frac{u}{u_{\theta}}\right)^{(j)}\!\!(x_k)\,.
	\end{equation*}
	Recalling the definition $\alpha = (M^{\ast})^{-1} b^\ast$, we obtain
	\begin{equation}
		\label{eq:def_d_minus_alpha}
		\big\| \delta - \alpha \big\|_2 =
		\big\| (M^{\ast})^{-1} (M^\ast \delta - b^\ast) \big\|_2 \leq
		\big\| (M^{\ast})^{-1} \big\|_2
		\big\| M^\ast \delta - b^\ast \big\|_2.
	\end{equation}
	We first take care of the term in $M^\ast \delta - b^\ast$.
	We have
	\begin{equation*}
		\begin{aligned}
			\left\| M^\ast \delta - b^\ast \right\|_2^2
			 & =
			\sum_{\ell=0}^q \left[
				\sum_{j=0}^q M^\ast_{\ell j} \delta_{j} - b^\ast_\ell
			\right]^2 \\
			 & =
			\sum_{\ell=0}^q \left[
				\sum_{j=0}^q
				\left(
				\int_{\Omega_k}
				\frac{{\phi}^*_{k,\ell}(x)}{(\Delta x_k)^{\ell}}
				\frac{{\phi}^*_{k,j}(x)}{(\Delta x_k)^{j}}
				dx
				\right)
				(\Delta x_k)^{j}
				\left( \frac{u}{u_{\theta}}\right)^{(j)}\!\!(x_k)\,
				-
				\int_{\Omega_k}
				u(x)
				\frac{{\phi}^*_{k,\ell}(x)}{(\Delta x_k)^{\ell}}
				dx
				\right]^2 \eqqcolon
			\sum_{\ell=0}^q \Xi_\ell^2
		\end{aligned}
	\end{equation*}
	We denote the summand by $\Xi_\ell$,
	and we use the definition of the basis to obtain
	\begin{align}
		\notag
		\forall \ell \in \{0, \dots, q\}, \quad
		\Xi_\ell & \coloneqq
		\sum_{j=0}^q
		(\Delta x_k)^j
		\left(\frac{u}{u_{\theta}}\right)^{(j)}\!\!(x_k)\,
		\int_{\Omega_k}
		\frac{{\phi}^*_{k,\ell}(x)}{(\Delta x_k)^{\ell}}
		\frac{{\phi}^*_{k,j}(x)}{(\Delta x_k)^{j}}
		dx -
		\int_{\Omega_k} u(x)
		\frac{{\phi}^*_{k,\ell}(x)}{(\Delta x_k)^{\ell}} dx \\
		\notag
		         & \, =
		\sum_{j=0}^q
		(\Delta x_k)^j
		\left(\frac{u}{u_{\theta}}\right)^{(j)}\!\!(x_k)\,
		\int_{\Omega_k}
		\frac{{\phi}_{k,\ell}(x)}{(\Delta x_k)^{\ell}}
		\frac{{\phi}_{k,j}(x)}{(\Delta x_k)^{j}}
		u_{\theta}^2(x) dx -
		\int_{\Omega_k} \frac{u(x)}{u_{\theta}(x)}
		\frac{{\phi}_{k,\ell}(x)}{(\Delta x_k)^{\ell}}
		u_{\theta}^2(x) dx                                  \\
		\label{eq:def_xi_j}
		         & \, =
		\int_{\Omega_k}  \left(
		\sum_{j=0}^q
		\left( \frac{u}{u_{\theta}}\right)^{(j)}\!\!(x_k)\, \phi_{k,j}(x)
		-
		\frac{u(x)}{u_{\theta}(x)}
		\right) \frac{\phi_{k,\ell}(x)}{(\Delta x_k)^{\ell}} u_{\theta}^2(x) dx.
	\end{align}
	Using a Taylor expansion, we obtain, for all $\ell \in \{0, \dots, q\}$,
	\begin{equation*}
		\begin{aligned}
			\Xi_\ell = 
			 & -
			\int_{\Omega_k} \left(
			\int_{x_k}^x
			\left(\frac{u}{u_{\theta}}\right)^{(q+1)}\!\!(\xi) \,
			\frac{(x-\xi)^{q}}{q!} d\xi \right)
			\frac{\phi_{k,\ell}(x)}{(\Delta x_k)^{\ell}} u_{\theta}^2(x) dx, \\
		\end{aligned}
	\end{equation*}
	from which we get the following upper bound
	\begin{equation*}
		\forall \ell \in \{0, \dots, q\}, \quad
		| \Xi_\ell | \leq
		\sup_{x \in \Omega_k} \Biggl|
		\int_{x_k}^x
		\left(\frac{u}{u_{\theta}}\right)^{(q+1)}\!\!(\xi) \,
		\frac{(x-\xi)^{q}}{q!} d\xi
		\Biggr| \,\,
		\Biggl|\int_{\Omega_k} \frac{\phi_{k,j}(x)}{(\Delta x_k)^{j}} u_{\theta}^2(x) dx \Biggr|.
	\end{equation*}
	Using the same ingredients as
	in the computation of \eqref{majo_chiwangshu}
	for the leftmost term
	and bounding the rightmost term
	by the $L^{\infty}$ norm of the prior and by
	noting that the classical basis functions are bounded,
	we obtain the estimate
	\begin{equation*}
		\forall \ell \in \{0, \dots, q\}, \quad
		| \Xi_\ell | \lesssim
		\left| \frac{u}{u_{\theta}} \right|_{H^{q+1}(\Omega_k)} \,
		(\Delta x_k)^{q+\frac12} \,
		(\Delta x_k)\left\|u_{\theta}^2 \right\|_{L^{\infty}(\Omega_k)}.
	\end{equation*}
	Going back to what we had set out to prove, we get
	\begin{equation}
		\label{eq:second_estimateproof_3}
		\left\| M^\ast \delta - b^\ast \right\|_2 =
		\left(\sum_{\ell=0}^q | \Xi_\ell |^2 \right)^{\frac 1 2} \lesssim
		\left| \frac{u}{u_{\theta}} \right|_{H^{q+1}(\Omega_k)} \,
		(\Delta x_k)^{q+\frac12}\, (\Delta x_k)
		\left\|u_{\theta}^2 \right\|_{L^{\infty}(\Omega_k)}.
	\end{equation}
	Plugging \eqref{eq:second_estimateproof_3}
	into \eqref{eq:def_d_minus_alpha} and then into
	\eqref{eq:second_estimateproof_1}, we get
	\begin{equation}
		\label{eq:second_estimateproof_4}
		N_2 \lesssim
		\left\| (M^\ast)^{-1} \right\|_2
		\left| \frac{u}{u_{\theta}} \right|_{H^{q+1}(\Omega_k)} \,
		(\Delta x_k)^{q+\frac12} \, (\Delta x_k)
		\left\|u_{\theta}^2 \right\|_{L^{\infty}(\Omega_k)}
		\left\|u_{\theta}\right\|_{L^{\infty}(\Omega_k)}.
	\end{equation}
	Finally, we note that, for any $y \in \mathbb{R}^{q+1}$,
	given the expression \eqref{eq:def_M_star_and_b} of $M^\ast$,
	\begin{equation*}
		\begin{aligned}
			\left\langle
			M^\ast y, M^\ast y
			\right\rangle =
			\int_{\Omega_k} \left(\sum_{j=0}^q \frac{{\phi}^*_{k,j}(x)}{(\Delta x_k)^j} y_j\right)^2 dx & =
			\int_{\Omega_k} \left(\sum_{j=0}^q \frac{{\phi}_j(x)}{(\Delta x_k)^j} y_j\right)^2 u_\theta(x)^2 dx      \\
			                                                                                            & \geqslant
			m^2\int_{\Omega_k} \left(\sum_{j=0}^q \frac{{\phi}_j(x)}{(\Delta x_k)^j} y_j\right)^2  dx = m^2\, \left\langle M y, M y \right\rangle,
		\end{aligned}
	\end{equation*}
	where $M$ is the mass matrix associated with the classical basis functions
	\[
		M_{j\ell}
		=
		\int_{\Omega_k} \frac{{\phi}_j(x)}{(\Delta x_k)^j} \frac{{\phi}_\ell(x)}{(\Delta x_k)^\ell} dx
		=
		\frac{\Delta x_k}{1 + j +\ell}
		=
		\Delta x_k H_{j\ell},
	\]
	where $H = (H_{j\ell})_{j \ell}$ is the Hilbert matrix.
	Then we deduce the following inequality
	\begin{equation}
		\label{eq:bound_on_M_star}
		\| (M^\ast)^{-1} \|_2
		\leqslant
		\frac{1}{m^2} \| M^{-1}\|_2
		=
		\frac{1}{m^2} \frac 1 {\Delta x_k} \| H^{-1}\|_2.
	\end{equation}
	Combining \eqref{eq:second_estimateproof_4} and \eqref{eq:bound_on_M_star},
	we obtain
	\begin{equation}
		\label{eq:second_estimateproof_5}
		N_2 \lesssim
		\left| \frac{u}{u_{\theta}} \right|_{H^{q+1}(\Omega_k)} \,
		(\Delta x_k)^{q+\frac12}
		\frac{\left\|u_{\theta}^2 \right\|_{L^{\infty}(\Omega_k)}}{m^2}
		\left\|u_{\theta}\right\|_{L^{\infty}(\Omega_k)}.
	\end{equation}

	We get, from \eqref{eq:estimate_priorproof3} and
	\eqref{eq:second_estimateproof_5},
	the expected result.
\end{proof}

The above proof relies on the smoothness of the prior.
This may seem counter-intuitive in a hyperbolic context.
However, since the prior will be obtained from a neural network
in \cref{sec:prior_construction},
this smoothness assumption becomes reasonable.

\begin{lemma}
	\label{lem:global_error_estimate}
	We make the same assumptions as in the previous lemma,
	and still consider the vector space $V_h^{*}$.
	For any function $u \in H^{q+1}(\Omega)$,
	\[
		\left\| u - P_h(u) \right\|_{L^2(\Omega)} \lesssim
		\left|
		\frac{u}{u_{\theta}}
		\right|_{H^{q+1}(\Omega)} \,
		(\Delta x_k)^{q+1} \,
		\| u_{\theta} \|_{L^\infty(\Omega)}.
	\]
\end{lemma}

\begin{proof}
	We begin by stating the definition of the discrete $L^2$ norm:
	by assuming that
	\begin{equation*}
		\Omega = \bigcup_{k=1}^N \Omega_k,
	\end{equation*}
	we obtain
	\[
		\left\| u - P_h(u) \right\|_{L^2(\Omega)}^2 \leqslant
		\sum_{k=1}^N \Delta x_k \, \| u - P_h(u) \|^2_{L^\infty(\Omega_k)}.
	\]
	Using the result from \cref{lem:local_error_estimate}, we get
	\begin{equation*}
		\begin{aligned}
			\| u - P_h(u) \|_{L^2(\Omega)}^2 & \lesssim
			\sum_{k=1}^N \Delta x_k \left(
			\left| \frac{u}{u_{\theta}} \right|_{H^{q+1}(\Omega_k)}
			(\Delta x_k)^{q+\frac12} \,
			\left( 1 + \frac{\left\|u_{\theta}^2 \right\|_{L^{\infty}(\Omega_k)}}{m^2} \right)
			\left\| u_{\theta} \right\|_{L^\infty(\Omega_k)}
			\right)^2                                           \\
			                                 & \lesssim
			\sum_{k=1}^N
			\left| \frac{u}{u_{\theta}} \right|^2_{H^{q+1}(\Omega_k)}
			(\Delta x_k)^{2q+2} \,
			\left( 1 + \frac{\left\|u_{\theta}^2 \right\|_{L^{\infty}(\Omega_k)}}{m^2} \right)^2
			\left\| u_{\theta} \right\|_{L^\infty(\Omega_k)}^2. \\
		\end{aligned}
	\end{equation*}
	We assume that there exists
	$\delta_-$, $\delta_+$ and $\Delta x$ such that,
	for all $k \in \{1, \dots, N\}$,
	$\delta_- \Delta x \leqslant \Delta x_k \leqslant \delta_+ \Delta x$.
	Then, since
	$\|u_{\theta}\|_{L^\infty(\Omega_k)} \leqslant \|u_{\theta}\|_{L^\infty(\Omega)}$,
	we obtain
	\begin{equation*}
		\begin{aligned}
			\| u - P_h(u) \|_{L^2(\Omega)}^2 & \lesssim
			(\Delta x)^{2q+2}
			\sum_{k=1}^N
			\left| \frac{u}{u_{\theta}} \right|^2_{H^{q+1}(\Omega_k)}
			\left( 1 + \frac{\left\|u_{\theta}^2 \right\|_{L^{\infty}(\Omega)}}{m^2} \right)^2
			\left\| u_{\theta} \right\|_{L^\infty(\Omega_k)}^2         \\ & \lesssim
			(\Delta x)^{2q+2}
			\left( 1 + \frac{\left\|u_{\theta}^2 \right\|_{L^{\infty}(\Omega_k)}}{m^2} \right)^2
			\left\| u_{\theta} \right\|_{L^\infty(\Omega)}^2
			\sum_{k=1}^N
			\left| \frac{u}{u_{\theta}} \right|^2_{H^{q+1}(\Omega_k)}. \\
		\end{aligned}
	\end{equation*}
	The proof is concluded by recognizing the $H^{q+1}(\Omega)$ seminorm.
\end{proof}

The global error estimate provided in \cref{lem:global_error_estimate}
shows that the projection error onto the basis $V_h^{*}$
is bounded by a term depending on the prior:
\[
	\left|
	\frac{u}{u_{\theta}}
	\right|_{H^{q+1}(\Omega)}.
\]
This bound is equal to zero if the prior is exact,
since it is nothing but the ($q+1$)\textsuperscript{th}
derivative of the constant function equal to one.
This estimate also proves that the closer the prior is to the solution,
the smaller the bound of the projection error.
However, to obtain an even smaller bound, we need the prior
and the solution to be close in the sense of the
$H^{q+1}(\Omega)$ seminorm.
This means that the prior must be
constructed in such a way that it also gives a good approximation
of the derivatives of the solution.

As a summary, we have shown that the $L^2$ projection error
tends to zero when the prior tends to the solution.
This result gives an idea of the expected gains in error ensured
by using the modified basis $V_h^{*}$.
The final convergence error depends on this projection error,
as has been shown in \cite{YUAN2006295}.
The proof to obtain the final convergence result
is the same as in \cite{YUAN2006295}.

For the additive basis $V_h^{+}$,
such error estimates are harder to obtain,
since the projection in the new basis is harder to write
with respect to the traditional one.
We expect the error to be bounded by a term
in $| u-u_{\theta}|_{H^{q+1}(\Omega)}$,
which would enable us to draw similar conclusions
as for the multiplicative basis $V_h^{*}$.
Namely, the error would also tend to zero when
the prior tends to the solution,
and the derivatives of the prior would need to
be close to the derivatives of the solution.
Proving this result is out of the scope of this paper,
even though it should be ensured by the results of \cite{YUAN2006295}.
However, we will extensively study the behavior of the additive basis
in \cref{sec:application_numerical_results}.

\section{Prior construction and algorithm}
\label{sec:prior_construction}

Equipped with the modified bases from
\cref{sec:modal_enriched_dg_scheme}
and with the theoretical results from \cref{sec:theoretical_results},
what is left to do is to propose a way to obtain a suitable
prior $u_\theta$.

Note that the approach described in \cref{sec:modified_dg_scheme}
will be interesting if the prior $u_\theta$ is a good
approximation of the steady solution to \eqref{eq:main_model}
for a wide range of parameters.
In addition, according to \cref{sec:estimate_prior},
the derivatives of the prior must also
be good approximations of the derivatives of the steady solution.

This means that we wish to capture large families of solutions,
i.e., we want to be able to calculate an approximation
for several parameters.
For example, assuming that \eqref{eq:main_model} depends on $4$
physical parameters leads to $\mu \in \mathbb{R}^4$,
and considering a problem in two space dimensions,
leads to $x \in \mathbb{R}^2$.
Therefore, we are looking for a prior $u_\theta(x; \mu)$,
where $u_\theta \in \mathcal{C}^q(\mathbb{R}^2 \times \mathbb{R}^4, \mathbb{R})$.
Approaching such a function using polynomials
defined on a mesh would be a very difficult task,
especially if the space or parameter domains have a complex geometry.
Neural networks have demonstrated their ability to approximate
functions in fairly high dimensions,
notably thanks to their intrinsic regularity.
PINNs are a mesh-free approach to solving PDEs using neural networks.
Their properties make them good candidates
for approaching solutions to high-dimensional problems.

To build our prior, we propose to solve the parametric
steady problem with PINNs.
To that end, we now briefly introduce this method
in \cref{sec:parametric_PINNs},
and we show how to compute and store the prior.
Then, our algorithm is summarized in \cref{sec:algorithm}.

\subsection{Parametric PINNs}
\label{sec:parametric_PINNs}

Note that the steady solutions to \eqref{eq:main_model} are given by
\begin{equation*}\label{eq:main_model_steady_solutions}
	\partial_x F_{{\mu}_1}(u)= S_{{\mu}_2}({\ra x,}u).
\end{equation*}
{\rb This is nothing but a parametric space-dependent PDE}.
Therefore, we introduce PINNs for
the following generic boundary value problems (BVPs):
\begin{equation}
	\label{eq:generic_BVP}
	\begin{dcases}
		\mathcal{D}({u}, x; {\mu}) = f(x , {\mu})
		 & \text{ for } x \in \Omega,          \\
		{u}(t,x) = g(x, {\mu}),
		 & \text{ for } x \in \partial \Omega, \\
	\end{dcases}
\end{equation}
where $\mathcal{D}$ is a differential operator
involving the solution ${u}$ and its space derivatives,
and with ${\mu}$
some physical parameters.
We recall that ${\mu} \in \mathbb{P} \subset \mathbb{R}^m$.
PINNs use the fact that classical fully-connected neural networks
are smooth functions of their inputs,
as long as their activation functions are also smooth,
to approximate the solution to \eqref{eq:generic_BVP}.
Contrary to traditional numerical schemes such as the DG method,
where the degrees of freedom encode some explicit
modal or nodal values of the solutions,
the degrees of freedom of PINNs representation are the weights $\theta$
of the neural network, and so they o not explicitly represent the solution.
Equipped with both of these remarks,
the idea behind PINNs is to plug the network,
which represents the solution to \eqref{eq:generic_BVP},
into the equation.
Then, the degrees of freedom (i.e., the weights $\theta$ of the network)
are found by minimizing a loss function.
Since the neural network is differentiable,
its derivatives can be exactly computed.
In our case, the PINN is thus a smooth neural network
that takes as input the space variable $x$
and the parameter vector ${\mu}$,
which we denote by~${u}_{\theta}(x; {\mu})$.

Thanks to these definitions, solving the PDE can be rewritten
as the following minimization problem:
\begin{equation}\label{eq:chap67_MC}
	\min_{\theta} \mathcal{J}(\theta),
	\text{\quad where }
	\mathcal{J}(\theta) =
	\mathcal{J}_r (\theta) +
	\mathcal{J}_b(\theta) +
	\mathcal{J}_\text{data}(\theta).
\end{equation}
In \eqref{eq:chap67_MC}, we have introduced three different terms:
the residual loss function $\mathcal{J}_r$,
the boundary loss function $\mathcal{J}_b$,
and the data loss function $\mathcal{J}_\text{data}$.
For parameters ${\mu} \in \mathbb{P}$,
the residual loss function is defined by
\begin{equation}\label{eq:chap67_MC1}
	\mathcal{J}_r(\theta) =
	\int_{\mathbb{P}}
	\int_{\Omega}
	\big\| \mathcal{D}({u}_{\theta}, x; {\mu}) - f(x; {\mu}) \big\|_2^2
	\, dx d{\mu},
\end{equation}
while the boundary loss function is given by
\begin{equation}\label{eq:chap67_MC2}
	\mathcal{J}_b(\theta) =
	\int_{\mathbb{P}}
	\int_{\partial \Omega}
	\big \| {u}_{\theta}(x, {\mu})-g(x, {\mu})\big\|_2^2
	\, dx d{\mu}.
\end{equation}
Finally, to define the data loss function,
we assume that we know the exact solution
to \eqref{eq:generic_BVP} at some points $x_i$
and for some parameters ${\mu}_i$, and we set
\begin{equation*}\label{eq:data_loss}
	\mathcal{J}_\text{data}(\theta) =
	\sum_i \big\|
	{u}_{\theta}(x_i, {\mu}_i) -
	{u}(x_i, {\mu}_i)
	\big\|_2^2.
\end{equation*}
In practice, the integrals in
\eqref{eq:chap67_MC1} and \eqref{eq:chap67_MC2}
are approximated
using a Monte-Carlo method.
This method relies on sampling a certain number of so-called
``collocation points'' in order to approximate the integrals.
Then, the minimization problem on $\theta$ is solved
using a gradient-type method,
which corresponds to the learning phase.

The main advantage of PINNs is that they are mesh-free and
less sensitive to dimension than classical methods.
Indeed, neural networks easily deal with large input dimensions,
and the Monte-Carlo method converges independently of the dimension.
Consequently, PINNs are particularly well-suited to
solving parametric PDEs such as \eqref{eq:generic_BVP}.
Thanks to that, we do not solve for a single equilibrium
but rather for families of equilibria
indexed by the parameters ${\mu}$.

Traditional PINNs use this method to approximate both
\eqref{eq:chap67_MC1} and \eqref{eq:chap67_MC2}.
However, for the boundary conditions,
we elected to use another approach,
which makes it possible to completely eliminate $\mathcal{J}_b$
from the minimization algorithm.
The idea is to define the approximate solution through a
boundary operator~$\mathcal{B}$,
which can for instance be a multiplication
by a function which satisfies the boundary condition.
We obtain
\[
	\widetilde{{u}}_{\theta}(x; {\mu}) =
	\mathcal{B} \big({u}_{\theta}, x ; {\mu} \big),
\]
with ${u}_{\theta}$ the neural network
and $\mathcal{B}$ a simple operator
such as $\widetilde{{u}}_{\theta}$ exactly satisfies
the boundary conditions.
Using~$\widetilde{{u}}_{\theta}$, the residual loss becomes,
instead of \eqref{eq:chap67_MC1}:
\begin{equation}\label{eq:chap67_MC1_with_BC}
	\mathcal{J}_r(\theta) =
	\int_{\mathbb{P}}
	\int_{\Omega}
	\big\| \mathcal{D}(\widetilde{{u}}_{\theta}, x; {\mu}) - f(x; {\mu}) \big\|_2^2
	\, dx d{\mu}.
\end{equation}
Examples of such functions $\mathcal{B}$ are provided in
\cref{sec:application_numerical_results}.
{\rb We emphasize that, with this strategy, the PINN output~${u}_{\theta}$
no longer represents an approximate solution to the BVP \eqref{eq:generic_BVP}.
Instead, approximate solutions to the BVP, and thus priors,
will be given by $\widetilde{{u}}_{\theta}$.}

With this approach, we have presented one method for
offline construction of our prior for a family of equilibria.
Note that it is possible to further enhance this prior with
data from previous simulations,
thanks to the loss function $\mathcal{J}_\text{data}$.
Even though training PINNs may be harder than training traditional purely data-driven
neural networks, they are much more efficient as priors.
Indeed, the error estimates of \cref{sec:estimate_prior}
show that the error depends on the $(q+1)^\text{th}$
derivative of the ratio between the prior and the solution.
Therefore, to obtain a small error, it is important for the prior
to provide a good approximation not only of the steady solution,
but also of its derivatives.
Since the PINN loss function~\eqref{eq:chap67_MC1_with_BC}
inherently contains derivatives of ${u}_{\theta}$,
the resulting trained PINN will be more efficient in this respect.
Note that a purely data-driven network could also be interesting
if the data contains information on the derivatives.

\subsection{Algorithm}
\label{sec:algorithm}

Now that we have discussed the strategy we use to obtain our prior,
we give some details on the offline and online algorithms
that we developed to construct the modified DG bases in practice.
We start by describing the offline step,
where the families of priors are computed.
Then, we move on to an online algorithm,
explaining how to construct the DG bases using the prior,
and how to apply them to the actual DG time iterations.

\begin{algorithm}[H]
	\caption{Offline part: neural network training}\label{alg:offline}
	\begin{algorithmic}[1]
		\Require
		space domain $\Omega$,
		parameter set $\mathbb{P}$,
		initial neural network ${u}_{\theta_0}(x;{\mu})$,
		learning rate $\eta$,
		number $N$ of collocation points,
		number $n_\text{epochs}$ of training epochs
		\Ensure trained neural network ${u}_{\theta}(x;{\mu})$

		\State initialize the weights: $\theta = \theta_0$
		\For{$n \le n_\text{epochs}$}
		\State sample $N$ values of $x$ in $\Omega$
		and ${\mu}$ in $\mathbb{P}$
		\State compute the loss function $\mathcal{J}(\theta)$
		\State update $\theta$ using
		the gradient of $\mathcal{J}(\theta)$:
		\(
		\theta = \theta - \eta \nabla_{\theta} J(\theta)
		\)
		\EndFor
	\end{algorithmic}
\end{algorithm}

In practice, we do not use a classical gradient descent
to update the weights, but rather the Adam algorithm.
Moreover, sampling is done through a uniform law on
the space and parameter domains.
It would also be possible use non-uniform sampling
like in \cite{WuZhuTanKarLu2023} for instance,
but we elected to use uniform sampling for the sake of simplicity.
Note that \cref{alg:offline} does not contain solution data
in its inputs.
Indeed, almost all numerical experiments from
\cref{sec:application_numerical_results}
do not require data on the solution.
This avoids the cost of data production,
which would otherwise require sampling the exact solution if it is known,
or using a numerical scheme otherwise.

\begin{algorithm}[H]
	\caption{Online part: using the neural network in the DG scheme}
	\label{alg:online}
	\begin{algorithmic}[1]
		\Require prior $\widetilde u_{\theta}$,
		degree $n_q$ of the Gauss-Lobatto quadrature rule,
		initial data $u_0$,
		space mesh $\Omega_h$,
		parameters~${\mu}$,
		number~$n_t$ of time steps
		\Ensure numerical solution $u_k(t,x)$
		on each cell $\Omega_k$
		\State use the mesh $\Omega_h$ to obtain all quadrature points
		$x_{k,p}$ in each cell $\Omega_k$
		\State evaluate the prior at each point $x_{k,p}$:
		we obtain $\widetilde u_{k,p} \coloneqq \widetilde u_{\theta}(x_{k,p};{\mu})$
		\State reconstruct $u_k(0,x)$ using $\widetilde u_{k,p}$
		\For{$n \le n_{t}$}
		\State construct the mass matrix $\mathcal{M}$,
		the nonlinear flux $\mathcal{V}$,
		the interface flux $\mathcal{I}$
		and the source term $\mathcal{S}$
		using $\widetilde u_{k,p}$ and the quadrature rule
		\State update the solution $u_k$
		at the next time step, using $u_k$ at the previous time step
		as well as the terms computed in the previous step
		\EndFor
	\end{algorithmic}
\end{algorithm}

In this second step,
the additional computational overhead associated to our method,
compared to the classical~DG scheme, comes from two distinct sources.
The first one is a preprocessing phase,
where we evaluate the prior on the quadrature points
(step 2 of \cref{alg:online}).
Even though such networks have been made
to be quickly evaluated on GPUs,
this evaluation step remains fast on CPUs.
The second source of computational cost is associated
to the quadrature rule.
Indeed, in some cases, we will require a quadrature rule
with a higher degree than the traditional DG scheme.
The classical approach is to use $n_q = q+2$ quadrature points
for bases made of $q+1$ polynomial functions,
since the quadrature is exact for polynomials of degree $q$.
However, in our case, our basis is non-polynomial.
Hence, to have a good approximation of the integral of the prior,
we may need to increase the degree of the quadrature.
In most cases, this increase is slight; for a few test cases,
especially to approximate functions with large derivatives,
we will need to use fine quadrature rules.

\section{Applications and numerical results}
\label{sec:application_numerical_results}

This section is dedicated to a validation of the approach
on several parametric hyperbolic systems of balance laws:
the linear advection equation in \cref{sec:transport},
the 1D shallow water equations in \cref{sec:shallow_water},
the Euler-Poisson system in \cref{sec:euler_poisson},
and the 2D shallow water equations in \cref{sec:shallow_water_2d}.
In the first two cases, there exist some exact well-balanced schemes
in the literature.
However, for the Euler-Poisson system and the
2D shallow water equations,
exact (or even approximate) WB schemes are either not available
or very complicated to implement.
Code replicating the experiments of \cref{sec:transport}
is freely available~\cite{Mic2023} on \textsf{GitHub}\footnote{\href{https://github.com/Victor-MichelDansac/DG-PINNs.git}{https://github.com/Victor-MichelDansac/DG-PINNs.git}}.

In this section, we denote by $K$
the number of cells.
We test both bases $V_h^{*}$ and $V_h^{+}$ at first,
showing that both display similar results.
To cut down on the number of tables,
we then only present the results for the additive basis~$V_h^{+}$.
	{\rb Since almost all experiments concern either the preservation
		of steady solutions or the study of perturbed steady solutions,
		the exact steady solution is prescribed as
		inhomogeneous Dirichlet boundary conditions in the DG scheme,
		unless otherwise mentioned.}

Moreover, the time step $\Delta t$ is given by
\begin{equation}
	\label{eq:time_step}
	\Delta t =
	C_\text{CFL} \,
	C_\text{RK} \,
	\frac{\min_{k \in \{1, \dots, K\}} \Delta x_k}{\lambda},
\end{equation}
where $\lambda$ is the maximal wave speed of the system,
$C_\text{CFL}$ is a CFL (Courant-Friedrichs-Lewy) coefficient,
and $C_\text{RK}$
is the stability coefficient associated to the time discretization.
All experiments are run using a strong stability-preserving
Runge-Kutta (SSPRK) time discretization of the correct order.
The time discretizations, with their associated
stability coefficients $C_\text{RK}$,
are collected in \cref{tab:time_discretizations}.
To determine $C_\text{CFL}$,
we run a study of the stability condition for the first experiment;
this study is not repeated for the other experiments,
since the new bases do not influence the stability condition.

\begin{table}[!ht]
	\centering
	\begin{tabular}{ccccc}
		\toprule
		degree $q$ of the DG basis          & $0$ & $1$ & $2$    & $3$ \\
		\cmidrule(lr){1-5}
		time discretization                 &
		explicit Euler                      &
		SSPRK2 \cite{SpiRuu2002}            &
		SSPRK3(5) \cite{SpiRuu2002}         &
		SSPRK4(10) \cite{Ket2008}                                      \\
		stability coefficient $C_\text{RK}$ & $1$ & $1$ & $2.65$ & $3$ \\
		\bottomrule
	\end{tabular}
	\caption{%
		Stability coefficients $C_\text{RK}$ of the
		high-order time discretizations used in the numerical experiments,
		with respect to the degree $q$ of the DG basis.%
	}
	\label{tab:time_discretizations}
\end{table}

\subsection{Linear advection}
\label{sec:transport}

We first consider the case of a linear advection equation
with a source term, on the space domain $\Omega = (0, 1)$.
The equation is given as follows:
\begin{equation}
	\label{eq:advection}
	\begin{cases}
		\partial_t u + \partial_{x}u = s(u; {\mu}),
		 & \text{ for } x \in \Omega,     \\
		u(t=0,x)=u_\text{ini}(x; {\mu}), \\
		u(t,x=0)=u_0,                    \\
	\end{cases}
\end{equation}
Here, the parameter vector ${\mu}$ is made of three elements:
\begin{equation*}
	{\mu} =
	\begin{pmatrix}
		\alpha \\
		\beta  \\
		u_0
	\end{pmatrix} \in \mathbb{P} \subset \mathbb{R}^3,
	\quad \alpha \in \mathbb{R}_+,
	\quad \beta \in \mathbb{R}_+,
	\quad u_0 \in \mathbb{R}_+^*.
\end{equation*}
The source term depends on ${\mu}$ as follows:
\begin{equation*}
	s(u; {\mu}) = \alpha u + \beta u^2,
\end{equation*}
and straightforward computations show
that the associated steady solutions take the form
\begin{equation}
	\label{eq:steady_solution_transport}
	u_\text{eq}(x; {\mu}) =
	\frac{\alpha u_0} {(\alpha + \beta u_0 ) e^{- \alpha x} - \beta u_0}.
\end{equation}
To compute the time step,
we take $\lambda = 1$ in \eqref{eq:time_step},
since the advection velocity in \eqref{eq:advection} is equal to $1$.
The first paragraph of this section shows how to choose $C_\text{CFL}$
to complete the determination of the time step.
Unless otherwise stated, we prescribe Dirichlet boundary conditions
consisting in the steady solution.

To obtain a suitable prior $\widetilde u_\theta$,
we train a PINN with parameters $\theta$.
To avoid cumbersome penalization of boundary conditions,
we define $\smash{\widetilde {u}_\theta}$
using a boundary operator $\mathcal{B}$, as follows:
\begin{equation*}
	\widetilde {u}_\theta(x; {\mu}) =
	\mathcal{B}({u}_\theta, x; {\mu}) =
	u_0 + x u_\theta(x; {\mu}),
\end{equation*}
so that the boundary condition
$\widetilde {u}_\theta(0; {\mu}) = u_0$
is automatically satisfied by $\widetilde {u}_\theta$.
	{\rb We highlight once again that the~PINN output ${u}_\theta$ does not represent the prior;
		rather, the prior is the function $\widetilde {u}_\theta$, which contains the boundary conditions.}
The parameter space $\mathbb{P}$ is chosen such that
the steady solution is well-defined, and we take
\begin{equation}
	\label{eq:parameter_space_advection}
	\mathbb{P} = [0.5, 1] \times [0.5, 1] \times [0.1, 0.2].
\end{equation}
Thanks to the boundary operator $\mathcal{B}$,
the loss function only concerns the ODE residual, and we set
\begin{equation*}
	\mathcal{J}(\theta)
	=
	\left\|
	\partial_x \widetilde{u}_\theta -
	\alpha \widetilde{u}_\theta -
	\beta \widetilde{u}_\theta^2
	\right\|_2^2.
\end{equation*}
{\rb%
This means that the result $u_\theta$ of the PINN will be optimized such that
$\partial_x {\widetilde{u}_\theta}$ is as close as possible to
$\alpha {\widetilde{u}_\theta} + \beta {\widetilde{u}_\theta}^2$,
which is nothing but the equation describing the steady solutions of \eqref{eq:advection}.%
}
We use a neural network with 5 fully connected hidden layers,
and around {\ra $2300$} trainable parameters.
Training takes about 4 minutes on a dual NVIDIA K80 GPU,
until the loss function is equal to about $10^{-6}$.
For this experiment, we increased the order of the quadrature
compared to the baseline for the case with one basis function.
Indeed, we take $n_q = \max(q + 2, 3)$,
to ensure a sufficient precision when integrating the prior.

In this section, we compare four strategies:
the basis $V_h$ \eqref{eq:modal_basis},
the basis $\smash{V_h^*}$ with multiplicative prior \eqref{eq:modal_basis_multiplicative},
the basis $\smash{V_h^+}$ with additive prior \eqref{eq:modal_basis_additive},
and the basis $\smash{V_h^{\text{ex},+}}$ which uses
the exact steady solution \eqref{eq:steady_solution_transport} as a prior.
First, we study the stability condition in \cref{sec:transport_stability}.
Then, we tackle the approximation of a steady solution
without perturbation in \cref{sec:transport_steady}
and with perturbation in \cref{sec:transport_steady_perturbed}.
Afterwards, the approximation of an unsteady solution
is computed in \cref{sec:transport_unsteady}.
{\ra Finally, computation time comparisons are provided in \cref{sec:transport_computation_time}.}

\subsubsection{Study of the stability condition}
\label{sec:transport_stability}

The very first experiment we run aims at making sure that the new bases
do not alter the stability condition of the DG scheme.
To that end, we slowly increase $C_\text{CFL}$
until the time step $\Delta t$ is too large for the scheme to be stable.
For this experiment, the initial condition is made of the steady solution
\begin{equation}
	\label{eq:initial_condition_advection_unperturbed}
	u_\text{ini}(x; \mu)
	= u_\text{eq}(x; \mu),
\end{equation}
and the final time is $T = 0.5$.
\cref{tab:advection_stability} contains the optimal values of $C_\text{CFL}$
(larger values leading to instabilities)
obtained with the four bases and for $q \in \{0, 1, 2, 3\}$.
We observe that the new bases do not change the stability condition,
except for $\smash{V_h^{\text{ex},+}}$ with $q=1$, which is slightly more stable.
This study will not be repeated for other experiments,
since it would yield similar results.
In practice, we take $C_\text{CFL} = 0.1$ to ensure stability.

\begin{table}[!ht]
	\centering
	\begin{tabular}{cccccc}
		\toprule
		$q$ &  & basis $V_h$ & basis $\smash{V_h^*}$ & basis $\smash{V_h^+}$ & basis $\smash{V_h^{\text{ex},+}}$ \\
		\cmidrule(lr){1-6}
		0   &  & 1.250       & 1.250         & 1.250         & 1.250                     \\
		1   &  & 0.399       & 0.399         & 0.399         & 0.416                     \\
		2   &  & 0.209       & 0.209         & 0.209         & 0.209                     \\
		3   &  & 0.185       & 0.185         & 0.185         & 0.185                     \\
		\bottomrule
	\end{tabular}
	\caption{%
		Maximal values of $C_\text{CFL}$ obtained for the four bases
		and for a number of basis elements $q \in \{0, 1, 2, 3\}$.
	}
	\label{tab:advection_stability}
\end{table}

\subsubsection{Steady solution}
\label{sec:transport_steady}

We now study the approximation of a steady solution,
at first without perturbation.
The goal of this section is to check whether the prior
indeed makes it possible to decrease the error compared to the usual modal basis.
For this experiment, the initial condition remains
\eqref{eq:initial_condition_advection_unperturbed},
and the final time is $T = 0.1$.

As a first step, the values of the parameters ${\mu}$ are set
to the midpoints of the intervals making up the
parameter space~\eqref{eq:parameter_space_advection}.
The $L^2$ errors between the exact and approximate solutions are
collected in \cref{tab:advection_steady}.

In this case, we expect both $\smash{V_h^*}$ and $\smash{V_h^+}$ to show similar behavior.
Moreover, we expect the basis $\smash{V_h^{\text{ex},+}}$ to provide
an exactly well-balanced scheme, up to machine precision.
To that end, only for $\smash{V_h^{\text{ex},+}}$, we take $\smash{n_q = \max(q + 2, 5)}$,
to ensure that the quadrature of the exact prior is also exact, up to machine precision.

\begin{table}[ht!]
	\centering
	\makebox[\textwidth][c]{%
		\subfloat[Errors with a basis made of one element: $q=0$.]{%
			\begin{tabular}{cccccccccccccc}
				\toprule
				    &  & \multicolumn{2}{c}{basis $V_h$} &       & \multicolumn{3}{c}{basis $\smash{V_h^*}$} &                      & \multicolumn{3}{c}{basis $\smash{V_h^+}$} &         & basis $\smash{V_h^{\text{ex},+}}$                                                                     \\
				\cmidrule(lr){3-4}\cmidrule(lr){6-8}\cmidrule(lr){10-12}\cmidrule(lr){14-14}
				$K$ &  & error                           & order &                                   & error                & order                             & gain    &                           & error                & order & gain    &  & error                 \\
				\cmidrule(lr){1-1}\cmidrule(lr){3-4}\cmidrule(lr){6-8}\cmidrule(lr){10-12}\cmidrule(lr){14-14}
				10  &  & $1.75 \cdot 10^{-2}$            & ---   &                                   & $1.45 \cdot 10^{-5}$ & ---                               & 1200.02 &                           & $1.45 \cdot 10^{-5}$ & ---   & 1200.02 &  & $1.66 \cdot 10^{-14}$ \\
				20  &  & $8.75 \cdot 10^{-3}$            & 1.00  &                                   & $7.61 \cdot 10^{-6}$ & 0.93                              & 1149.11 &                           & $7.61 \cdot 10^{-6}$ & 0.93  & 1149.11 &  & $2.78 \cdot 10^{-17}$ \\
				40  &  & $4.38 \cdot 10^{-3}$            & 1.00  &                                   & $3.92 \cdot 10^{-6}$ & 0.96                              & 1118.29 &                           & $3.92 \cdot 10^{-6}$ & 0.96  & 1118.29 &  & $5.89 \cdot 10^{-17}$ \\
				80  &  & $2.19 \cdot 10^{-3}$            & 1.00  &                                   & $2.00 \cdot 10^{-6}$ & 0.97                              & 1098.77 &                           & $2.00 \cdot 10^{-6}$ & 0.97  & 1098.77 &  & $2.78 \cdot 10^{-17}$ \\
				160 &  & $1.10 \cdot 10^{-3}$            & 1.00  &                                   & $1.01 \cdot 10^{-6}$ & 0.98                              & 1085.96 &                           & $1.01 \cdot 10^{-6}$ & 0.98  & 1085.96 &  & $2.19 \cdot 10^{-17}$ \\
				\bottomrule
			\end{tabular}
		}
		\label{tab:advection_steady_nG_1}
	}

	\vspace{12pt}


	\makebox[\textwidth][c]{%
		\subfloat[Errors with a basis made of two elements: $q=1$.]{%
			\begin{tabular}{cccccccccccccc}
				\toprule
				    &  & \multicolumn{2}{c}{basis $V_h$} &       & \multicolumn{3}{c}{basis $\smash{V_h^*}$} &                      & \multicolumn{3}{c}{basis $\smash{V_h^+}$} &        & basis $\smash{V_h^{\text{ex},+}}$                                                                    \\
				\cmidrule(lr){3-4}\cmidrule(lr){6-8}\cmidrule(lr){10-12}\cmidrule(lr){14-14}
				$K$ &  & error                           & order &                                   & error                & order                             & gain   &                           & error                & order & gain   &  & error                 \\
				\cmidrule(lr){1-1}\cmidrule(lr){3-4}\cmidrule(lr){6-8}\cmidrule(lr){10-12}\cmidrule(lr){14-14}
				10  &  & $4.93 \cdot 10^{-4}$            & ---   &                                   & $2.18 \cdot 10^{-6}$ & ---                               & 226.00 &                           & $1.03 \cdot 10^{-6}$ & ---   & 479.78 &  & $2.04 \cdot 10^{-14}$ \\
				20  &  & $1.24 \cdot 10^{-4}$            & 2.00  &                                   & $3.20 \cdot 10^{-7}$ & 2.77                              & 386.66 &                           & $2.74 \cdot 10^{-7}$ & 1.91  & 450.59 &  & $6.48 \cdot 10^{-16}$ \\
				40  &  & $3.09 \cdot 10^{-5}$            & 2.00  &                                   & $8.07 \cdot 10^{-8}$ & 1.99                              & 382.88 &                           & $7.00 \cdot 10^{-8}$ & 1.97  & 441.39 &  & $9.46 \cdot 10^{-16}$ \\
				80  &  & $7.72 \cdot 10^{-6}$            & 2.00  &                                   & $2.05 \cdot 10^{-8}$ & 1.98                              & 376.52 &                           & $1.76 \cdot 10^{-8}$ & 1.99  & 438.98 &  & $1.46 \cdot 10^{-15}$ \\
				160 &  & $1.93 \cdot 10^{-6}$            & 2.00  &                                   & $5.16 \cdot 10^{-9}$ & 1.99                              & 374.49 &                           & $4.40 \cdot 10^{-9}$ & 2.00  & 438.34 &  & $2.13 \cdot 10^{-15}$ \\
				\bottomrule
			\end{tabular}
		}
		\label{tab:advection_steady_nG_2}
	}


	\vspace{12pt}

	\makebox[\textwidth][c]{%
		\subfloat[Errors with a basis made of three elements: $q=2$.]{%
			\begin{tabular}{cccccccccccccc}
				\toprule
				    &  & \multicolumn{2}{c}{basis $V_h$}                                &       & \multicolumn{3}{c}{basis $\smash{V_h^*}$} &                                                                & \multicolumn{3}{c}{basis $\smash{V_h^+}$} &       & basis $\smash{V_h^{\text{ex},+}}$                                                                                                                                                      \\
				\cmidrule(lr){3-4}\cmidrule(lr){6-8}\cmidrule(lr){10-12}\cmidrule(lr){14-14}
				$K$ &  & error                                                          & order &                                   & error                                                          & order                             & gain  &                           & error                                                          & order & gain  &  & error                                                          \\
				\cmidrule(lr){1-1}\cmidrule(lr){3-4}\cmidrule(lr){6-8}\cmidrule(lr){10-12}\cmidrule(lr){14-14}
				10  &  & \tablenum[table-format=1.2e1,exponent-product=\cdot]{7.89e-06} & ---   &                                   & \tablenum[table-format=1.2e2,exponent-product=\cdot]{1.00e-07} & ---                               & 78.58 &                           & \tablenum[table-format=1.2e2,exponent-product=\cdot]{1.05e-07} & ---   & 74.90 &  & \tablenum[table-format=1.2e2,exponent-product=\cdot]{9.92e-13} \\
				20  &  & \tablenum[table-format=1.2e1,exponent-product=\cdot]{9.94e-07} & 2.99  &                                   & \tablenum[table-format=1.2e2,exponent-product=\cdot]{1.33e-08} & 2.91                              & 74.60 &                           & \tablenum[table-format=1.2e2,exponent-product=\cdot]{1.41e-08} & 2.90  & 70.65 &  & \tablenum[table-format=1.2e2,exponent-product=\cdot]{7.84e-15} \\
				40  &  & \tablenum[table-format=1.2e1,exponent-product=\cdot]{1.24e-07} & 3.00  &                                   & \tablenum[table-format=1.2e2,exponent-product=\cdot]{1.72e-09} & 2.95                              & 72.13 &                           & \tablenum[table-format=1.2e2,exponent-product=\cdot]{1.79e-09} & 2.97  & 69.19 &  & \tablenum[table-format=1.2e2,exponent-product=\cdot]{2.84e-15} \\
				80  &  & \tablenum[table-format=1.2e1,exponent-product=\cdot]{1.55e-08} & 3.00  &                                   & \tablenum[table-format=1.2e2,exponent-product=\cdot]{2.17e-10} & 2.99                              & 71.43 &                           & \tablenum[table-format=1.2e2,exponent-product=\cdot]{2.25e-10} & 2.99  & 68.81 &  & \tablenum[table-format=1.2e2,exponent-product=\cdot]{7.81e-15} \\
				160 &  & \tablenum[table-format=1.2e1,exponent-product=\cdot]{1.94e-09} & 3.00  &                                   & \tablenum[table-format=1.2e2,exponent-product=\cdot]{2.72e-11} & 3.00                              & 71.25 &                           & \tablenum[table-format=1.2e2,exponent-product=\cdot]{2.82e-11} & 3.00  & 68.72 &  & \tablenum[table-format=1.2e2,exponent-product=\cdot]{1.15e-14} \\
				\bottomrule
			\end{tabular}
		}
		\label{tab:advection_steady_nG_3}
	}

	\vspace{12pt}

	\makebox[\textwidth][c]{%
		\subfloat[Errors with a basis made of four elements: $q=3$.]{%
			\begin{tabular}{cccccccccccccc}
				\toprule
				    &  & \multicolumn{2}{c}{basis $V_h$}                                &       & \multicolumn{3}{c}{basis $\smash{V_h^*}$} &                                                                & \multicolumn{3}{c}{basis $\smash{V_h^+}$} &       & basis $\smash{V_h^{\text{ex},+}}$                                                                                                                                                      \\
				\cmidrule(lr){3-4}\cmidrule(lr){6-8}\cmidrule(lr){10-12}\cmidrule(lr){14-14}
				$K$ &  & error                                                          & order &                                   & error                                                          & order                             & gain  &                           & error                                                          & order & gain  &  & error                                                          \\
				\cmidrule(lr){1-1}\cmidrule(lr){3-4}\cmidrule(lr){6-8}\cmidrule(lr){10-12}\cmidrule(lr){14-14}
				10  &  & \tablenum[table-format=1.2e2,exponent-product=\cdot]{1.20e-07} & ---   &                                   & \tablenum[table-format=1.2e2,exponent-product=\cdot]{8.31e-09} & ---                               & 14.40 &                           & \tablenum[table-format=1.2e2,exponent-product=\cdot]{1.12e-08} & ---   & 10.67 &  & \tablenum[table-format=1.2e2,exponent-product=\cdot]{4.45e-11} \\
				20  &  & \tablenum[table-format=1.2e2,exponent-product=\cdot]{7.39e-09} & 4.02  &                                   & \tablenum[table-format=1.2e2,exponent-product=\cdot]{5.51e-10} & 3.91                              & 13.40 &                           & \tablenum[table-format=1.2e2,exponent-product=\cdot]{7.28e-10} & 3.95  & 10.15 &  & \tablenum[table-format=1.2e2,exponent-product=\cdot]{7.72e-13} \\
				40  &  & \tablenum[table-format=1.2e2,exponent-product=\cdot]{4.59e-10} & 4.01  &                                   & \tablenum[table-format=1.2e2,exponent-product=\cdot]{3.48e-11} & 3.99                              & 13.19 &                           & \tablenum[table-format=1.2e2,exponent-product=\cdot]{4.56e-11} & 4.00  & 10.06 &  & \tablenum[table-format=1.2e2,exponent-product=\cdot]{1.70e-14} \\
				80  &  & \tablenum[table-format=1.2e2,exponent-product=\cdot]{2.92e-11} & 3.98  &                                   & \tablenum[table-format=1.2e2,exponent-product=\cdot]{2.20e-12} & 3.99                              & 13.27 &                           & \tablenum[table-format=1.2e2,exponent-product=\cdot]{2.86e-12} & 3.99  & 10.18 &  & \tablenum[table-format=1.2e2,exponent-product=\cdot]{6.93e-15} \\
				160 &  & \tablenum[table-format=1.2e2,exponent-product=\cdot]{1.85e-12} & 3.98  &                                   & \tablenum[table-format=1.2e2,exponent-product=\cdot]{1.29e-13} & 4.10                              & 14.38 &                           & \tablenum[table-format=1.2e2,exponent-product=\cdot]{1.72e-13} & 4.06  & 10.76 &  & \tablenum[table-format=1.2e2,exponent-product=\cdot]{2.59e-14} \\
				\bottomrule
			\end{tabular}
		}
		\label{tab:advection_steady_nG_4}
	}
	\caption{%
	Advection equation: errors, orders of accuracy,
	and gain obtained when approximating a steady solution
	for bases without prior (basis $V_h$),
	with a PINN prior (bases $\smash{V_h^*}$ and $\smash{V_h^+}$),
	and with an exact prior (basis $\smash{V_h^{\text{ex},+}}$).
	}
	\label{tab:advection_steady}
\end{table}

We observe that the bases with and without prior allow
a convergence of the correct order,
i.e. of the same order as the number of basis elements.
Moreover, we observe a consistent gain for all mesh resolutions,
for a given size of the modal basis,
which is lower the larger the size of the basis.
Bases $\smash{V_h^*}$ and $\smash{V_h^+}$ seem to have comparable performance,
with $\smash{V_h^*}$ being somewhat better for large values of $q$,
and $\smash{V_h^+}$ taking the lead for small values of $q$.
Finally, we observe that the basis $\smash{V_h^{\text{ex},+}}$
is indeed able to provide a solution that is exact up to machine precision,
thus validating the exact well-balanced property of the scheme using this basis.

As a second step, to refine this study,
we now consider $10^3$ parameters,
randomly sampled from the parameter space \eqref{eq:parameter_space_advection}.
For $q \in \{0, 1, 2, 3\}$ and $K = 10$ discretization cells,
we compute the minimum, average and maximum gains
obtained with both bases $\smash{V_h^*}$ and $\smash{V_h^+}$.
These values are reported in \cref{tab:advection_steady_stats}.
We observe, on average, a significant gain in all cases,
with larger gains obtained for smaller values of $q$.
Furthermore, the minimum gain is always greater than one.
Like in the previous experiment, we observe that,
even though both bases display similar behavior and very good results,
$\smash{V_h^+}$ behaves better than $\smash{V_h^*}$ for small values of $q$, and vise versa.
Consequently, and to limit the number of tables in the remainder of this section,
we perform all subsequent experiments with the basis $\smash{V_h^+}$.

\begin{table}[!ht]
	\centering
	\begin{tabular}{ccccccccc}
		\toprule
		      &  & \multicolumn{3}{c}{gains in basis $\smash{V_h^*}$} &                                     & \multicolumn{3}{c}{gains in basis $\smash{V_h^+}$}                                                                                                                      \\
		\cmidrule(lr){3-5}\cmidrule(lr){7-9}
		{$q$} &  & {minimum}                                  & {average}                           & {maximum}                                  &  & {minimum}                          & {average}                           & {maximum}                            \\
		\cmidrule(lr){1-1}\cmidrule(lr){3-5}\cmidrule(lr){7-9}
		0     &  & \tablenum[table-format=2.2]{63.46}         & \tablenum[table-format=3.2]{735.08} & \tablenum[table-format=4.2]{4571.89}       &  & \tablenum[table-format=2.2]{63.46} & \tablenum[table-format=3.2]{735.08} & \tablenum[table-format=4.2]{4571.89} \\
		1     &  & \tablenum[table-format=2.2]{32.22}         & \tablenum[table-format=3.2]{149.38} & \tablenum[table-format=4.2]{450.74}        &  & \tablenum[table-format=2.2]{26.01} & \tablenum[table-format=3.2]{190.08} & \tablenum[table-format=4.2]{830.20}  \\
		2     &  & \tablenum[table-format=2.2]{6.20}          & \tablenum[table-format=3.2]{54.16}  & \tablenum[table-format=4.2]{118.45}        &  & \tablenum[table-format=2.2]{5.92}  & \tablenum[table-format=3.2]{45.47}  & \tablenum[table-format=4.2]{313.07}  \\
		3     &  & \tablenum[table-format=2.2]{1.55}          & \tablenum[table-format=3.2]{19.54}  & \tablenum[table-format=4.2]{108.10}        &  & \tablenum[table-format=2.2]{1.56}  & \tablenum[table-format=3.2]{13.69}  & \tablenum[table-format=4.2]{184.17}  \\
		\bottomrule
	\end{tabular}
	\caption{%
		Advection equation: statistics of the gains obtained
		for the approximation of a steady solution
		in bases $\smash{V_h^*}$ and $\smash{V_h^+}$ with respect to basis $V_h$.
	}
	\label{tab:advection_steady_stats}
\end{table}

\subsubsection{Perturbed steady solution}
\label{sec:transport_steady_perturbed}

\begin{figure}[!ht]
	\centering
	\makebox[\textwidth][c]{\includegraphics{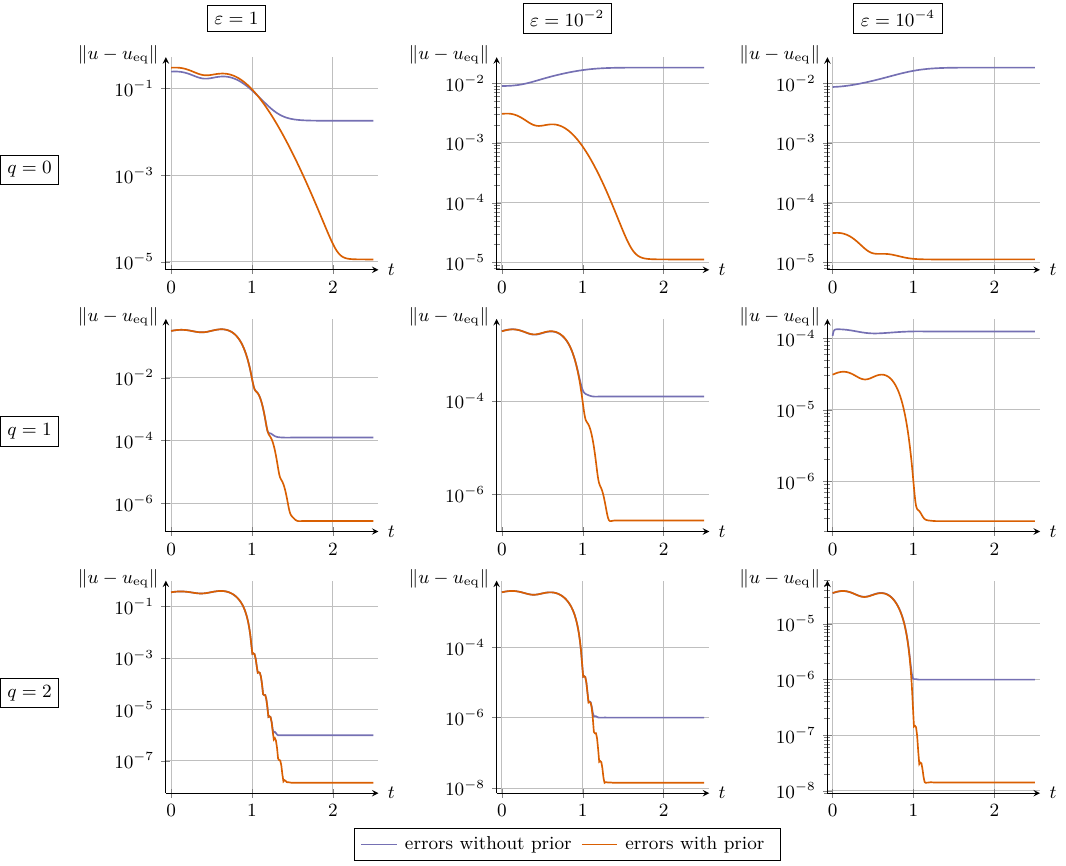}}
	\caption{%
	Advection equation:
	{\rb distance between the DG solution $u$
	and the underlying steady solution $u_\text{eq}$},
	with respect to time,
	for the approximation of a perturbed steady solution
	for bases with and without prior.
	}
	\label{fig:advection_steady_perturbed}
\end{figure}

{\ra We now test the scheme on a perturbed steady solution.}
For this experiment, the initial condition is similar
to~\eqref{eq:initial_condition_advection_unperturbed}, but with a perturbation.
Indeed, we take
\begin{equation*}
	u_\text{ini}(x; \mu)
	= \big( 1 + \varepsilon \sin (2 \pi x) \big)\,
	u_\text{eq}(x; \mu),
\end{equation*}
where $\varepsilon$ is taken nonzero or zero,
to control the strength of the perturbation.
The final time is $T = 2$,
and we study the impact of the perturbation by
taking $\varepsilon \in \{10^{-4},10^{-2}, 1\}$,
and~$K = 10$ discretization cells.
	{\rb The results are collected in \cref{fig:advection_steady_perturbed},
		where we display the errors between the DG approximation of $u$
		and the underlying steady solution $u_\text{eq}$.}
We observe two different {\rb situations}:
first, while the perturbation is being dissipated,
the errors with the two bases are similar.
Then, we note that the introduction of the prior has made it possible
for the approximate solution to converge towards a final solution
that is closer to the exact, unperturbed steady solution.

\subsubsection{Unsteady solution}
\label{sec:transport_unsteady}

Next, we seek to confirm that our proposed basis does not
deteriorate the approximation of unsteady solutions.
To that end, we consider an unsteady solution of
the homogeneous problem,
i.e. a solution to \eqref{eq:advection}
with $s(u; {\mu}) = 0$.
We take the following initial condition:
\begin{equation*}
	u_0(x) = 0.1 \left(1 + \exp \left( - 100 (x - 0.5)^2 \right) \right),
\end{equation*}
so that $u(t, x) = u_0(x - t)$.
The final time is set to $T=1$, and periodic boundary conditions are prescribed.

We compute the approximate solution with the two bases,
for several values of $q$.
The results are collected in \cref{tab:advection_unsteady}.
We note that the basis with prior
does not affect the approximate solution for $q \geq 1$,
while the results are slightly worse with the prior for $q = 0$.
To improve the results here, one could introduce a space-time basis
in a space-time discontinuous Galerkin method;
this will be the object of future work.

\begin{table}[ht!]
	\centering
	\makebox[\textwidth][c]{%
		\subfloat[Errors with a basis made of one element: $q=0$.]{%
			\begin{tabular}{cccccc}
				\toprule
				    & \multicolumn{2}{c}{basis $V_h$} & \multicolumn{3}{c}{basis $\smash{V_h^+}$}                                       \\
				\cmidrule(lr){2-3}\cmidrule(lr){4-6}
				$K$ & error                           & order                             & error                & order & gain \\
				\cmidrule(lr){1-1}\cmidrule(lr){2-3}\cmidrule(lr){4-6}
				10  & $4.04 \cdot 10^{-2}$            & ---                               & $5.04 \cdot 10^{-2}$ & ---   & 0.80 \\
				20  & $3.46 \cdot 10^{-2}$            & 0.22                              & $4.28 \cdot 10^{-2}$ & 0.24  & 0.81 \\
				40  & $2.84 \cdot 10^{-2}$            & 0.28                              & $3.50 \cdot 10^{-2}$ & 0.29  & 0.81 \\
				80  & $2.15 \cdot 10^{-2}$            & 0.40                              & $2.64 \cdot 10^{-2}$ & 0.40  & 0.81 \\
				160 & $1.47 \cdot 10^{-2}$            & 0.55                              & $1.81 \cdot 10^{-2}$ & 0.55  & 0.81 \\
				\bottomrule
			\end{tabular}
		}
		\label{tab:advection_unsteady_nG_1}
		\subfloat[Errors with a basis made of two elements: $q=1$.]{%
			\begin{tabular}{cccccc}
				\toprule
				    & \multicolumn{2}{c}{basis $V_h$} & \multicolumn{3}{c}{basis $\smash{V_h^+}$}                                       \\
				\cmidrule(lr){2-3}\cmidrule(lr){4-6}
				$K$ & error                           & order                             & error                & order & gain \\
				\cmidrule(lr){1-1}\cmidrule(lr){2-3}\cmidrule(lr){4-6}
				10  & $1.92 \cdot 10^{-2}$            & ---                               & $1.93 \cdot 10^{-2}$ & ---   & 1.00 \\
				20  & $6.26 \cdot 10^{-3}$            & 1.62                              & $6.27 \cdot 10^{-3}$ & 1.62  & 1.00 \\
				40  & $1.19 \cdot 10^{-3}$            & 2.39                              & $1.20 \cdot 10^{-3}$ & 2.39  & 1.00 \\
				80  & $1.99 \cdot 10^{-4}$            & 2.59                              & $1.99 \cdot 10^{-4}$ & 2.59  & 1.00 \\
				160 & $4.19 \cdot 10^{-5}$            & 2.24                              & $4.20 \cdot 10^{-5}$ & 2.24  & 1.00 \\
				\bottomrule
			\end{tabular}
		}
		\label{tab:advection_unsteady_nG_2}
	}


	\vspace{12pt}

	\makebox[\textwidth][c]{%
		\subfloat[Errors with a basis made of three elements: $q=2$.]{%
			\begin{tabular}{cccccc}
				\toprule
				    & \multicolumn{2}{c}{basis $V_h$} & \multicolumn{3}{c}{basis $\smash{V_h^+}$}                                       \\
				\cmidrule(lr){2-3}\cmidrule(lr){4-6}
				$K$ & error                           & order                             & error                & order & gain \\
				\cmidrule(lr){1-1}\cmidrule(lr){2-3}\cmidrule(lr){4-6}
				10  & $5.15 \cdot 10^{-3}$            & ---                               & $5.15 \cdot 10^{-3}$ & ---   & 1.00 \\
				20  & $4.56 \cdot 10^{-4}$            & 3.50                              & $4.56 \cdot 10^{-4}$ & 3.50  & 1.00 \\
				40  & $4.55 \cdot 10^{-5}$            & 3.32                              & $4.55 \cdot 10^{-5}$ & 3.32  & 1.00 \\
				80  & $5.42 \cdot 10^{-6}$            & 3.07                              & $5.42 \cdot 10^{-6}$ & 3.07  & 1.00 \\
				160 & $6.75 \cdot 10^{-7}$            & 3.01                              & $6.75 \cdot 10^{-7}$ & 3.01  & 1.00 \\
				\bottomrule
			\end{tabular}
		}
		\label{tab:advection_unsteady_nG_3}
		\subfloat[Errors with a basis made of four elements: $q=3$.]{%
			\begin{tabular}{cccccc}
				\toprule
				    & \multicolumn{2}{c}{basis $V_h$} & \multicolumn{3}{c}{basis $\smash{V_h^+}$}                                       \\
				\cmidrule(lr){2-3}\cmidrule(lr){4-6}
				$K$ & error                           & order                             & error                & order & gain \\
				\cmidrule(lr){1-1}\cmidrule(lr){2-3}\cmidrule(lr){4-6}
				10  & $4.72 \cdot 10^{-4}$            & ---                               & $4.72 \cdot 10^{-4}$ & ---   & 1.00 \\
				20  & $2.87 \cdot 10^{-5}$            & 4.04                              & $2.87 \cdot 10^{-5}$ & 4.04  & 1.00 \\
				40  & $1.81 \cdot 10^{-6}$            & 3.99                              & $1.81 \cdot 10^{-6}$ & 3.99  & 1.00 \\
				80  & $1.14 \cdot 10^{-7}$            & 3.98                              & $1.14 \cdot 10^{-7}$ & 3.98  & 1.00 \\
				160 & $7.20 \cdot 10^{-9}$            & 3.99                              & $7.20 \cdot 10^{-9}$ & 3.99  & 1.00 \\
				\bottomrule
			\end{tabular}
		}
		\label{tab:advection_unsteady_nG_4}
	}
	\caption{%
		Advection equation: errors, orders of accuracy,
		and gain obtained when approximating an unsteady solution
		for bases with and without prior.
	}
	\label{tab:advection_unsteady}
\end{table}

{\ra
\subsubsection{Computation time}
\label{sec:transport_computation_time}

This last part of the study on the advection equation
is dedicated to computation time comparisons.
They were performed on a server equipped with an NVIDIA Tesla V100 GPU.

\paragraph{PINN performance}
We are first concerned with comparing the performance of the PINN
with respect to the number of parameters.
We propose three different PINNs,
denoted by $\mathcal{M}_1$, $\mathcal{M}_2$, and $\mathcal{M}_3$,
all based on multilayer perceptrons whose architectures are described in
\cref{tab:advection_computation_time_pinn_1}.
In each case, we use $\num{5000}$ collocation points
to approximate the integrals in the loss function.
Each of these PINNs is trained 10 distinct times,
for a maximum of $\num{25000}$ epochs.

We first report, in \cref{tab:advection_computation_time_pinn_1},
the computation time for each epoch.
We observe that increasing the number of parameters
does not have a large impact on the computation time.
Moreover, we also report the best values
(over the 10 times each network has been trained) of the loss function
$\mathcal{J}_\theta$ that has been obtained after the $\num{25000}$ epochs.
We observe that increasing the number of parameters by a factor of 2
divides the best value of~$\mathcal{J}_\theta$ by a factor smaller than 2.

\begin{table}[!ht]
	\centering
	\begin{tabular}{rccc}
		\toprule                                                      &
		model $\mathcal{M}_1$                                         &
		model $\mathcal{M}_2$                                         &
		model $\mathcal{M}_3$                                             \\
		\cmidrule(lr){1-4}
		number of neurons in hidden layers                            &
		$16 \times 32 \times 16$                                      &
		$16 \times 32 \times 32 \times 16$                            &
		$16 \times 32 \times 64 \times 32 \times 16$                      \\
		number of parameters                                          &
		1243                                                          &
		2299                                                          &
		5435                                                              \\
		computation time for one epoch (in seconds)                   &
		$7.89 \cdot 10^{-3}$                                          &
		$8.60 \cdot 10^{-3}$                                          &
		$9.32 \cdot 10^{-3}$                                              \\
		best value of $\mathcal{J}_\theta$ after $\num{25000}$ epochs &
		$1.78 \cdot 10^{-7}$                                          &
		$8.39 \cdot 10^{-8}$                                          &
		$5.51 \cdot 10^{-8}$                                              \\
		\bottomrule
	\end{tabular}
	\caption{\ra %
		For each of the three PINNs used to approximate the steady advection equation,
		we report the number of neurons in the hidden layers,
		the total number of parameters,
		the computation time for one epoch,
		and the best value of $\mathcal{J}_\theta$ after $\num{25000}$ epochs
		(obtained by training the network 10 times and
		selecting the best value of~$\mathcal{J}_\theta$).
	}
	\label{tab:advection_computation_time_pinn_1}
\end{table}

To further investigate the performance of the PINN,
we refine this study in \cref{tab:advection_computation_time_pinn_2}
by reporting the number of epochs and the computation time
required to reach different values of $\mathcal{J}_\theta$.
To that end, we averaged the results obtained over the $\num{10}$ times
the networks have been trained.
We observe that, as the models become larger,
they are able to reach lower values of $\mathcal{J}_\theta$ more quickly.
In addition, the first model, $\mathcal{M}_1$, has been unable to reach
a value of $\mathcal{J}_\theta$ lower than $10^{-7}$ in $\num{25000}$ epochs.
Since the results of $\mathcal{M}_2$ and $\mathcal{M}_3$ are mostly similar
in terms of the best value reached for the loss function,
we kept on using $\mathcal{M}_2$ for the rest of this study.

\begin{table}[!ht]
	\centering
	\begin{tabular}{ccccccc}
		\toprule
		                                          &
		\multicolumn{2}{c}{model $\mathcal{M}_1$} &
		\multicolumn{2}{c}{model $\mathcal{M}_2$} &
		\multicolumn{2}{c}{model $\mathcal{M}_3$}                        \\
		\cmidrule(lr){2-3} \cmidrule(lr){4-5} \cmidrule(lr){6-7}
		value of $\mathcal{J}_\theta$             &
		epochs                                    & computation time &
		epochs                                    & computation time &
		epochs                                    & computation time     \\
		\cmidrule(lr){1-1} \cmidrule(lr){2-3} \cmidrule(lr){4-5} \cmidrule(lr){6-7}
		$10^{-4}$                                 &
		967                                       & 8.61             &
		607                                       & 6.29             &
		408                                       & 4.83                 \\
		$10^{-5}$                                 &
		2560                                      & 21.46            &
		1502                                      & 14.41            &
		1131                                      & 11.99                \\
		$10^{-6}$                                 &
		11390                                     & 91.09            &
		9513                                      & 83.51            &
		8681                                      & 82.35                \\
		$10^{-7}$                                 &
		---                                       & ---              &
		22365                                     & 193.68           &
		18611                                     & 174.23               \\
		\bottomrule
	\end{tabular}
	\caption{\ra %
		Comparison between the three PINN models
		described in \cref{tab:advection_computation_time_pinn_1}.
		For each model, we report the number of epochs
		and the computation time required to reach
		different values of the loss function $\mathcal{J}_\theta$
		(averaged over 10 distinct trainings).
	}
	\label{tab:advection_computation_time_pinn_2}
\end{table}

\paragraph{Enhanced DG performance}
After checking the performance of the PINN,
we now compare the computation times of the DG scheme
with and without the prior (made of model $\mathcal{M}_2$).
To that end, we run the DG scheme~$100$ times on the
test case from \cref{sec:transport_steady}, without perturbation,
and keep the 3 lowest computation times to avoid outliers,
before averaging them.
In each case, we report two different computation times:
the time it took to mesh the domain and assemble the DG matrices
(called ``assembly''),
and the time it took to perform a full time loop of the DG scheme
(called ``scheme'').
We expect the scheme time to be the same for both bases,
since the only difference between them is the computation of the prior,
which only intervenes during the assembly phase.

\begin{table}[!ht]
	\centering
	\makebox[\textwidth][c]{%
		\subfloat[Computation time for a basis made of one element: $q=0$.]{%
			\begin{tabular}{cccccccc}
				\toprule
				    & \multicolumn{3}{c}{basis $V_h$} & \multicolumn{4}{c}{basis $\smash{V_h^+}$}                                                                                                     \\
				\cmidrule(lr){2-4}\cmidrule(lr){5-8}
				$K$ & assembly                        & scheme                            & total                & assembly             & scheme               & total                & ratio \\
				\cmidrule(lr){1-1}\cmidrule(lr){2-4}\cmidrule(lr){5-8}
				10  & $7.08 \cdot 10^{-4}$            & $1.98 \cdot 10^{-2}$              & $2.06 \cdot 10^{-2}$ & $5.98 \cdot 10^{-3}$ & $1.98 \cdot 10^{-2}$ & $2.58 \cdot 10^{-2}$ & 1.26  \\
				20  & $6.78 \cdot 10^{-4}$            & $3.77 \cdot 10^{-2}$              & $3.84 \cdot 10^{-2}$ & $5.81 \cdot 10^{-3}$ & $3.78 \cdot 10^{-2}$ & $4.36 \cdot 10^{-2}$ & 1.14  \\
				40  & $7.18 \cdot 10^{-4}$            & $7.70 \cdot 10^{-2}$              & $7.77 \cdot 10^{-2}$ & $6.04 \cdot 10^{-3}$ & $7.70 \cdot 10^{-2}$ & $8.30 \cdot 10^{-2}$ & 1.07  \\
				80  & $7.41 \cdot 10^{-4}$            & $1.57 \cdot 10^{-1}$              & $1.57 \cdot 10^{-1}$ & $6.24 \cdot 10^{-3}$ & $1.57 \cdot 10^{-1}$ & $1.63 \cdot 10^{-1}$ & 1.04  \\
				160 & $7.61 \cdot 10^{-4}$            & $3.16 \cdot 10^{-1}$              & $3.17 \cdot 10^{-1}$ & $6.40 \cdot 10^{-3}$ & $3.13 \cdot 10^{-1}$ & $3.20 \cdot 10^{-1}$ & 1.01  \\
				\bottomrule
			\end{tabular}
		}
		\label{tab:advection_computation_time_q_0}
	}


	\vspace{12pt}

	\centering
	\makebox[\textwidth][c]{%
		\subfloat[Computation time for a basis made of two elements: $q=1$.]{%
			\begin{tabular}{cccccccc}
				\toprule
				    & \multicolumn{3}{c}{basis $V_h$} & \multicolumn{4}{c}{basis $\smash{V_h^+}$}                                                                                                     \\
				\cmidrule(lr){2-4}\cmidrule(lr){5-8}
				$K$ & assembly                        & scheme                            & total                & assembly             & scheme               & total                & ratio \\
				\cmidrule(lr){1-1}\cmidrule(lr){2-4}\cmidrule(lr){5-8}
				10  & $7.45 \cdot 10^{-4}$            & $4.07 \cdot 10^{-2}$              & $4.14 \cdot 10^{-2}$ & $5.81 \cdot 10^{-3}$ & $4.19 \cdot 10^{-2}$ & $4.78 \cdot 10^{-2}$ & 1.15  \\
				20  & $7.52 \cdot 10^{-4}$            & $8.19 \cdot 10^{-2}$              & $8.27 \cdot 10^{-2}$ & $5.94 \cdot 10^{-3}$ & $8.15 \cdot 10^{-2}$ & $8.74 \cdot 10^{-2}$ & 1.06  \\
				40  & $8.09 \cdot 10^{-4}$            & $1.72 \cdot 10^{-1}$              & $1.73 \cdot 10^{-1}$ & $6.09 \cdot 10^{-3}$ & $1.68 \cdot 10^{-1}$ & $1.74 \cdot 10^{-1}$ & 1.01  \\
				80  & $8.43 \cdot 10^{-4}$            & $3.52 \cdot 10^{-1}$              & $3.53 \cdot 10^{-1}$ & $6.12 \cdot 10^{-3}$ & $3.50 \cdot 10^{-1}$ & $3.56 \cdot 10^{-1}$ & 1.01  \\
				160 & $1.00 \cdot 10^{-3}$            & $7.45 \cdot 10^{-1}$              & $7.46 \cdot 10^{-1}$ & $6.50 \cdot 10^{-3}$ & $7.44 \cdot 10^{-1}$ & $7.50 \cdot 10^{-1}$ & 1.00  \\
				\bottomrule
			\end{tabular}
		}
		\label{tﬂab:advection_computation_time_q_1}
	}

	\vspace{12pt}

	\centering
	\makebox[\textwidth][c]{%
		\subfloat[Computation time for a basis made of three elements: $q=2$.]{%
			\begin{tabular}{cccccccc}
				\toprule
				    & \multicolumn{3}{c}{basis $V_h$} & \multicolumn{4}{c}{basis $\smash{V_h^+}$}                                                                                                     \\
				\cmidrule(lr){2-4}\cmidrule(lr){5-8}
				$K$ & assembly                        & scheme                            & total                & assembly             & scheme               & total                & ratio \\
				\cmidrule(lr){1-1}\cmidrule(lr){2-4}\cmidrule(lr){5-8}
				10  & $8.86 \cdot 10^{-4}$            & $4.04 \cdot 10^{-2}$              & $4.13 \cdot 10^{-2}$ & $6.17 \cdot 10^{-3}$ & $4.10 \cdot 10^{-2}$ & $4.72 \cdot 10^{-2}$ & 1.14  \\
				20  & $9.28 \cdot 10^{-4}$            & $8.21 \cdot 10^{-2}$              & $8.31 \cdot 10^{-2}$ & $6.14 \cdot 10^{-3}$ & $8.41 \cdot 10^{-2}$ & $9.03 \cdot 10^{-2}$ & 1.09  \\
				40  & $9.84 \cdot 10^{-4}$            & $1.75 \cdot 10^{-1}$              & $1.76 \cdot 10^{-1}$ & $6.32 \cdot 10^{-3}$ & $1.75 \cdot 10^{-1}$ & $1.81 \cdot 10^{-1}$ & 1.03  \\
				80  & $1.04 \cdot 10^{-3}$            & $3.60 \cdot 10^{-1}$              & $3.61 \cdot 10^{-1}$ & $6.39 \cdot 10^{-3}$ & $3.63 \cdot 10^{-1}$ & $3.69 \cdot 10^{-1}$ & 1.02  \\
				160 & $1.15 \cdot 10^{-3}$            & $7.95 \cdot 10^{-1}$              & $7.96 \cdot 10^{-1}$ & $6.94 \cdot 10^{-3}$ & $7.98 \cdot 10^{-1}$ & $8.05 \cdot 10^{-1}$ & 1.01  \\
				\bottomrule
			\end{tabular}
		}
		\label{tab:advection_computation_time_q_2}
	}

	\vspace{12pt}

	\centering
	\makebox[\textwidth][c]{%
		\subfloat[Computation time for a basis made of four elements: $q=3$.]{%
			\begin{tabular}{cccccccc}
				\toprule
				    & \multicolumn{3}{c}{basis $V_h$} & \multicolumn{4}{c}{basis $\smash{V_h^+}$}                                                                                                     \\
				\cmidrule(lr){2-4}\cmidrule(lr){5-8}
				$K$ & assembly                        & scheme                            & total                & assembly             & scheme               & total                & ratio \\
				\cmidrule(lr){1-1}\cmidrule(lr){2-4}\cmidrule(lr){5-8}
				10  & $1.02 \cdot 10^{-3}$            & $7.25 \cdot 10^{-2}$              & $7.35 \cdot 10^{-2}$ & $6.31 \cdot 10^{-3}$ & $7.25 \cdot 10^{-2}$ & $7.88 \cdot 10^{-2}$ & 1.07  \\
				20  & $1.03 \cdot 10^{-3}$            & $1.47 \cdot 10^{-1}$              & $1.48 \cdot 10^{-1}$ & $6.35 \cdot 10^{-3}$ & $1.47 \cdot 10^{-1}$ & $1.53 \cdot 10^{-1}$ & 1.04  \\
				40  & $1.11 \cdot 10^{-3}$            & $3.07 \cdot 10^{-1}$              & $3.08 \cdot 10^{-1}$ & $6.49 \cdot 10^{-3}$ & $3.00 \cdot 10^{-1}$ & $3.06 \cdot 10^{-1}$ & 1.00  \\
				80  & $1.24 \cdot 10^{-3}$            & $6.47 \cdot 10^{-1}$              & $6.48 \cdot 10^{-1}$ & $6.67 \cdot 10^{-3}$ & $6.50 \cdot 10^{-1}$ & $6.56 \cdot 10^{-1}$ & 1.01  \\
				160 & $1.34 \cdot 10^{-3}$            & $1.49 \cdot 10^{0}$               & $1.49 \cdot 10^{0}$  & $7.12 \cdot 10^{-3}$ & $1.50 \cdot 10^{0}$  & $1.50 \cdot 10^{0}$  & 1.01  \\
				\bottomrule
			\end{tabular}
		}
		\label{tab:advection_computation_time_q_3}
	}
	\caption{\ra %
		For $q \in \{0, 1, 2, 3\}$, we report the computation time
		needed for the matrix assembly and the full time loop of the DG scheme,
		for bases $V_h$ and $\smash{V_h^+}$, with and without prior.
		The ratio between the total computation time taken by each base is also reported.
	}
	\label{tab:advection_computation_time}
\end{table}

We compare the computation times with bases $V_h$ and $\smash{V_h^+}$ for different values of $q$,
and collect the results in \cref{tab:advection_computation_time}.
In each case, we observe that the time needed to run the scheme
is equivalent for both bases.
However, the assembly time can be up to 8 times larger with the prior.
This is expected, since the prior requires the evaluation of the PINN
at each quadrature point of each element of the mesh.
This potentially discouraging result is mitigated by the fact that
the assembly time is negligible compared to the scheme time,
especially when the number $K$ of mesh points,
or the number $q$ of basis functions, is large.
To highlight this, we have reported in the last column of each table
the ratio between the total computation time taken by each basis.
We observe that this ratio is very close to $1$,
except for the smallest values of $K$ and $q$.
These remarks validate the use of the prior to enhance DG bases,
since it is both efficient (as evidenced by the current section)
and accurate (as evidenced by
\cref{sec:transport_steady,sec:transport_steady_perturbed,sec:transport_unsteady}).
}

\subsection{Shallow water equations}
\label{sec:shallow_water}

After studying a scalar linear advection equation in
\cref{sec:transport}, we now turn to
a nonlinear system of conservation laws.
Namely, we tackle the shallow water equations
\begin{equation}
	\label{eq:shallow_water}
	\begin{dcases}
		\partial_t h+ \partial_{x} Q = 0, \\
		\partial_t Q + \partial_{x} \left( \frac{Q^2}{h} + \frac{1}{2}gh^2 \right) = -gh\partial_x Z(x; \alpha, \beta),
	\end{dcases}
\end{equation}
where $h > 0$ is the water height,
$Q$ the water discharge,
$g = 9.81$ the gravity constant,
and where the parameterized topography function is
\begin{equation}
	\label{eq:topography_function}
	Z(x; \alpha, \beta) = \beta \, \omega\left( \alpha \left(x - \frac 1 2\right)\right).
\end{equation}
In \eqref{eq:topography_function},
the function $\omega \in \{\omega_g, \omega_c\}$
is either a Gaussian bump function
\begin{equation}
	\label{eq:Gaussian_topography}
	\omega_g(x) = \frac 1 4 e^{-50 x^2}
\end{equation}
or a compactly supported bump function,
with parameter $\eth = 0.15$:
\begin{equation}
	\label{eq:compact_topography}
	\omega_c(x) =
	\begin{dcases}
		\exp \left(
		1 - \frac 1 {1 - \left(\dfrac{x}{\eth}\right)^2}
		\right) & \text{ if } |x| < \eth, \\
		0       & \text{ otherwise.}
	\end{dcases}
\end{equation}
Unless otherwise mentioned, the final physical time is $T = 0.05$,
and the space domain is $\Omega = (0, 1)$.
For each experiment, Dirichlet boundary conditions corresponding
to the steady solution are prescribed.

The steady solutions are given by cancelling
the time derivatives in \eqref{eq:shallow_water},
and we get the following characterization:
\begin{equation}
	\label{eq:shallow_water_steady}
	Q_\text{eq} = \text{constant} \eqqcolon Q_0
	\text{\qquad and \qquad}
	\left(1-\frac{Q_0^2}{g h_\text{eq}(x; \mu)^3}\right)
	\partial_x h_\text{eq}(x; \mu) + \partial_x Z(x; \alpha, \beta) = 0.
\end{equation}
To solve the nonlinear ODE on $h$,
we impose $h = h_0$ at some point in space.
Without loss of generality,
we restrict the study to the case $Q_0 > 0$.
This leads us to a family of steady solutions with four parameters,
and thus a parameter vector ${\mu}$ made of four elements:
\begin{equation*}
	{\mu} =
	\begin{pmatrix}
		\alpha \\
		\beta  \\
		h_0    \\
		Q_0
	\end{pmatrix}
	\in \mathbb{P} \subset \left( \mathbb{R}_+^* \right)^4.
\end{equation*}
To compute $\Delta t$ in \eqref{eq:time_step},
we take $\lambda = \frac{Q_0}{h_0} + \sqrt{g h_0}$.

Depending on the values of these parameters,
the Froude number
\begin{equation*}
	\text{Fr} = \sqrt{\frac{Q^2}{g h^3}}
\end{equation*}
controls the so-called
flow regime for the steady solution.
They can be in three distinct regimes:
subcritical ($\text{Fr} < 1$ everywhere),
supercritical ($\text{Fr} > 1$ everywhere) or
transcritical ($\text{Fr} = 1$ somewhere in the domain).
Each regime has its own parameter space
for $h_0$ and $Q_0$, described later,
but in all cases we take, unless otherwise stated,
\begin{equation}
	\label{eq:parameter_space_shallow_water_a_b}
	0.5 \leq \alpha \leq 1.5
	\quad ; \qquad
	0.5 \leq \beta \leq 1.5.
\end{equation}
To approximate the steady water height
within this parameter space,
we use a fully-connected PINN
with about~$\num{4000}$ trainable parameters.
Its result ${h}_\theta$
is modified through a boundary function $\mathcal{B}$
that will be defined for each regime.
The loss function is once again
made only of the steady ODE,
and we minimize
\begin{equation*}
	\mathcal{J}(\theta) =
	\left\|
	\left(
	1 - \frac{Q_0^2}{g \widetilde{h}_\theta(x; {\mu})^3}
	\right)
	\partial_x {\widetilde{h}_\theta}(x; {\mu})
	+
	\partial_x Z(x; \alpha, \beta)
	\right\|.
\end{equation*}
{\rb This means that the result ${h}_\theta$ of the PINN is defined
such that $\widetilde{h}_\theta$ is close to $h_\text{eq}$,
since $\widetilde{h}_\theta$ will (approximately)
satisfy the steady ODE \eqref{eq:shallow_water_steady}.}
Training takes about 5 minutes on a dual NVIDIA K80 GPU,
and lasts until the loss function is about $10^{-4}$,
depending on the regime.

\subsubsection{Subcritical flow}
\label{sec:shallow_water_subcritical}

{\rb
We start with a subcritical flow,
where we impose $h = h_0$ at the boundaries.
The parameter space for $h_0$ and $Q_0$ is:
\begin{equation}
	\label{eq:parameter_space_shallow_water_sub}
	2 \leq h_0 \leq 3
	\quad ; \qquad
	3 \leq Q_0 \leq 4.
\end{equation}
To test the preservation of the steady solution,
we set the initial water height to $h_\text{eq}$.

To strongly enforce the boundary conditions,
the prior $\smash{\widetilde{h}_\theta}$ is obtained as follows
from the result $h_\theta$ of the PINN:
\begin{equation}
	\label{eq:BC_PINN_SW_sub}
	\widetilde{h}_\theta(x; {\mu})
	=
	\mathcal{B}(h_\theta, x; {\mu})
	=
	h_0 + Z(x; \alpha, \beta) \, h_\theta(x; {\mu}).
\end{equation}
Since $Z$ is very close to $0$ at the boundaries
(or even equal to $0$ in the compactly supported case),
the expression~\eqref{eq:BC_PINN_SW_sub} ensures that
$\smash{\widetilde{h}_\theta(x; {\mu}) \simeq h_0}$ at the boundaries.
}

A goal of this section is to better understand
the differences between the two topography functions:
the Gaussian bump \eqref{eq:Gaussian_topography}
and the compactly supported bump
\eqref{eq:compact_topography}.
It is well-known that compactly supported functions
exhibit large derivatives close to the support,
see for instance \cite{SpiHuyDeB2015}.
As a consequence,
to get a good approximation of these derivatives
when computing integrals involving the PINN,
we take $n_q = q + 6$ when $\omega = \omega_c$.
Note that this choice
is also motivated by the results in \cite{SpiHuyDeB2015},
where the authors had to take larger polynomial degrees
to observe the correct orders of convergence.
The Gaussian topography also suffers from the same
drawback, but to a lesser extent,
and we take $n_q = q + 3$ when $\omega = \omega_g$
when integrating the result of the PINN.
For the compactly supported topography, the results are reported in \cref{tab:SW_steady_subcritical};
for the Gaussian topography, the results are reported in \cref{tab:SW_steady_subcritical_Gaussian}.

\begin{table}[!ht]
	\centering
	\subfloat[Error with a basis made of one element: $q=0$.]{%
		\begin{tabular}{ccccccccccc}
			\toprule
			    & \multicolumn{2}{c}{$h$, basis $V_h$} & \multicolumn{2}{c}{$Q$, basis $V_h$} & \multicolumn{3}{c}{$h$, basis $V_h^+$} & \multicolumn{3}{c}{$Q$, basis $V_h^+$}                                                                                                                                           \\
			\cmidrule(lr){2-3}\cmidrule(lr){4-5}\cmidrule(lr){6-8}\cmidrule(lr){9-11}
			$K$ & error                                & order                                & error                                  & order                                  & error                & order & gain                                & error                & order & gain                                \\
			\cmidrule(lr){1-1}\cmidrule(lr){2-3}\cmidrule(lr){4-5}\cmidrule(lr){6-8}\cmidrule(lr){9-11}
			20  & $5.76 \cdot 10^{-2}$                 & ---                                  & $1.89 \cdot 10^{-1}$                   & ---                                    & $6.20 \cdot 10^{-4}$ & ---   & \tablenum[table-format=3.2]{92.82}  & $2.90 \cdot 10^{-3}$ & ---   & \tablenum[table-format=3.2]{65.38}  \\
			40  & $3.06 \cdot 10^{-2}$                 & 0.91                                 & $1.50 \cdot 10^{-1}$                   & 0.34                                   & $5.65 \cdot 10^{-5}$ & 3.46  & \tablenum[table-format=3.2]{541.59} & $3.94 \cdot 10^{-4}$ & 2.88  & \tablenum[table-format=3.2]{380.39} \\
			80  & $1.82 \cdot 10^{-2}$                 & 0.75                                 & $8.30 \cdot 10^{-2}$                   & 0.85                                   & $3.46 \cdot 10^{-5}$ & 0.71  & \tablenum[table-format=3.2]{525.20} & $1.70 \cdot 10^{-4}$ & 1.21  & \tablenum[table-format=3.2]{488.25} \\
			160 & $9.94 \cdot 10^{-3}$                 & 0.87                                 & $4.53 \cdot 10^{-2}$                   & 0.87                                   & $1.94 \cdot 10^{-5}$ & 0.84  & \tablenum[table-format=3.2]{511.96} & $9.28 \cdot 10^{-5}$ & 0.87  & \tablenum[table-format=3.2]{488.31} \\
			320 & $5.26 \cdot 10^{-3}$                 & 0.92                                 & $2.37 \cdot 10^{-2}$                   & 0.93                                   & $1.04 \cdot 10^{-5}$ & 0.91  & \tablenum[table-format=3.2]{507.63} & $4.89 \cdot 10^{-5}$ & 0.92  & \tablenum[table-format=3.2]{484.02} \\
			\bottomrule
		\end{tabular}
	}

	\hspace{\fill}

	\vspace{12pt}

	\subfloat[Error with a basis made of two elements: $q=1$.]{%
		\begin{tabular}{ccccccccccc}
			\toprule
			    & \multicolumn{2}{c}{$h$, basis $V_h$} & \multicolumn{2}{c}{$Q$, basis $V_h$} & \multicolumn{3}{c}{$h$, basis $V_h^+$} & \multicolumn{3}{c}{$Q$, basis $V_h^+$}                                                                                 \\
			\cmidrule(lr){2-3}\cmidrule(lr){4-5}\cmidrule(lr){6-8}\cmidrule(lr){9-11}
			$K$ & error                                & order                                & error                                  & order                                  & error                & order & gain   & error                & order & gain   \\
			\cmidrule(lr){1-1}\cmidrule(lr){2-3}\cmidrule(lr){4-5}\cmidrule(lr){6-8}\cmidrule(lr){9-11}
			20  & $2.13 \cdot 10^{-2}$                 & ---                                  & $6.69 \cdot 10^{-2}$                   & ---                                    & $1.05 \cdot 10^{-4}$ & ---   & 202.69 & $3.96 \cdot 10^{-4}$ & ---   & 168.97 \\
			40  & $3.90 \cdot 10^{-3}$                 & 2.45                                 & $1.37 \cdot 10^{-2}$                   & 2.28                                   & $1.93 \cdot 10^{-5}$ & 2.44  & 202.12 & $8.14 \cdot 10^{-5}$ & 2.28  & 168.62 \\
			80  & $8.35 \cdot 10^{-4}$                 & 2.22                                 & $2.91 \cdot 10^{-3}$                   & 2.24                                   & $1.59 \cdot 10^{-6}$ & 3.60  & 525.18 & $7.27 \cdot 10^{-6}$ & 3.48  & 399.73 \\
			160 & $2.04 \cdot 10^{-4}$                 & 2.03                                 & $6.72 \cdot 10^{-4}$                   & 2.11                                   & $3.67 \cdot 10^{-7}$ & 2.11  & 556.12 & $1.55 \cdot 10^{-6}$ & 2.23  & 432.74 \\
			320 & $5.13 \cdot 10^{-5}$                 & 1.99                                 & $1.65 \cdot 10^{-4}$                   & 2.02                                   & $9.06 \cdot 10^{-8}$ & 2.02  & 566.17 & $3.62 \cdot 10^{-7}$ & 2.10  & 455.78 \\
			\bottomrule
		\end{tabular}
	}

	\hspace{\fill}

	\vspace{12pt}

	\subfloat[Error with a basis made of three elements: $q=2$.]{%
		\begin{tabular}{ccccccccccc}
			\toprule
			    & \multicolumn{2}{c}{$h$, basis $V_h$} & \multicolumn{2}{c}{$Q$, basis $V_h$} & \multicolumn{3}{c}{$h$, basis $V_h^+$} & \multicolumn{3}{c}{$Q$, basis $V_h^+$}                                                                                                                                           \\
			\cmidrule(lr){2-3}\cmidrule(lr){4-5}\cmidrule(lr){6-8}\cmidrule(lr){9-11}
			$K$ & error                                & order                                & error                                  & order                                  & error                & order & gain                                & error                & order & gain                                \\
			\cmidrule(lr){1-1}\cmidrule(lr){2-3}\cmidrule(lr){4-5}\cmidrule(lr){6-8}\cmidrule(lr){9-11}
			20  & $6.08 \cdot 10^{-3}$                 & ---                                  & $1.89 \cdot 10^{-2}$                   & ---                                    & $1.44 \cdot 10^{-4}$ & ---   & \tablenum[table-format=3.2]{42.26}  & $6.14 \cdot 10^{-4}$ & ---   & \tablenum[table-format=3.2]{30.75}  \\
			40  & $7.98 \cdot 10^{-4}$                 & 2.93                                 & $2.57 \cdot 10^{-3}$                   & 2.88                                   & $2.52 \cdot 10^{-6}$ & 5.83  & \tablenum[table-format=3.2]{316.56} & $7.71 \cdot 10^{-6}$ & 6.32  & \tablenum[table-format=3.2]{333.56} \\
			80  & $1.05 \cdot 10^{-4}$                 & 2.93                                 & $3.93 \cdot 10^{-4}$                   & 2.71                                   & $2.24 \cdot 10^{-7}$ & 3.49  & \tablenum[table-format=3.2]{467.99} & $8.54 \cdot 10^{-7}$ & 3.17  & \tablenum[table-format=3.2]{460.48} \\
			160 & $1.71 \cdot 10^{-5}$                 & 2.61                                 & $7.02 \cdot 10^{-5}$                   & 2.49                                   & $4.09 \cdot 10^{-8}$ & 2.45  & \tablenum[table-format=3.2]{418.00} & $1.76 \cdot 10^{-7}$ & 2.28  & \tablenum[table-format=3.2]{399.05} \\
			320 & $2.22 \cdot 10^{-6}$                 & 2.94                                 & $1.01 \cdot 10^{-5}$                   & 2.80                                   & $6.02 \cdot 10^{-9}$ & 2.77  & \tablenum[table-format=3.2]{369.32} & $2.91 \cdot 10^{-8}$ & 2.59  & \tablenum[table-format=3.2]{345.73} \\
			\bottomrule
		\end{tabular}
	}
	\caption{%
		Shallow water system,
		compactly supported topography
		\eqref{eq:compact_topography}:
		errors, orders of accuracy,
		and gain obtained when approximating
		a subcritical steady solution
		for bases with and without prior.
	}
	\label{tab:SW_steady_subcritical}
\end{table}

\begin{table}[!ht]
	\centering
	\subfloat[Error with a basis made of one element: $q=0$.]{%
		\begin{tabular}{ccccccccccc}
			\toprule
			    & \multicolumn{2}{c}{$h$, basis $V_h$} & \multicolumn{2}{c}{$Q$, basis $V_h$} & \multicolumn{3}{c}{$h$, basis $V_h^+$} & \multicolumn{3}{c}{$Q$, basis $V_h^+$}                                                                                 \\
			\cmidrule(lr){2-3}\cmidrule(lr){4-5}\cmidrule(lr){6-8}\cmidrule(lr){9-11}
			$K$ & error                                & order                                & error                                  & order                                  & error                & order & gain   & error                & order & gain   \\
			\cmidrule(lr){1-1}\cmidrule(lr){2-3}\cmidrule(lr){4-5}\cmidrule(lr){6-8}\cmidrule(lr){9-11}
			20  & $4.07 \cdot 10^{-2}$                 & ---                                  & $1.65 \cdot 10^{-1}$                   & ---                                    & $9.27 \cdot 10^{-5}$ & ---   & 439.14 & $3.87 \cdot 10^{-4}$ & ---   & 425.63 \\
			40  & $2.30 \cdot 10^{-2}$                 & 0.83                                 & $1.04 \cdot 10^{-1}$                   & 0.67                                   & $5.85 \cdot 10^{-5}$ & 0.67  & 393.05 & $2.65 \cdot 10^{-4}$ & 0.55  & 391.09 \\
			80  & $1.26 \cdot 10^{-2}$                 & 0.86                                 & $5.86 \cdot 10^{-2}$                   & 0.82                                   & $3.28 \cdot 10^{-5}$ & 0.83  & 384.93 & $1.60 \cdot 10^{-4}$ & 0.73  & 366.15 \\
			160 & $6.74 \cdot 10^{-3}$                 & 0.91                                 & $3.13 \cdot 10^{-2}$                   & 0.90                                   & $1.74 \cdot 10^{-5}$ & 0.91  & 386.04 & $8.72 \cdot 10^{-5}$ & 0.88  & 359.19 \\
			320 & $3.50 \cdot 10^{-3}$                 & 0.95                                 & $1.62 \cdot 10^{-2}$                   & 0.95                                   & $9.27 \cdot 10^{-6}$ & 0.91  & 377.48 & $4.56 \cdot 10^{-5}$ & 0.94  & 356.17 \\
			\bottomrule
		\end{tabular}
	}

	\hspace{\fill}

	\vspace{12pt}

	\subfloat[Error with a basis made of two elements: $q=1$.]{%
		\begin{tabular}{ccccccccccc}
			\toprule
			    & \multicolumn{2}{c}{$h$, basis $V_h$} & \multicolumn{2}{c}{$Q$, basis $V_h$} & \multicolumn{3}{c}{$h$, basis $V_h^+$} & \multicolumn{3}{c}{$Q$, basis $V_h^+$}                                                                                 \\
			\cmidrule(lr){2-3}\cmidrule(lr){4-5}\cmidrule(lr){6-8}\cmidrule(lr){9-11}
			$K$ & error                                & order                                & error                                  & order                                  & error                & order & gain   & error                & order & gain   \\
			\cmidrule(lr){1-1}\cmidrule(lr){2-3}\cmidrule(lr){4-5}\cmidrule(lr){6-8}\cmidrule(lr){9-11}
			20  & $3.21 \cdot 10^{-3}$                 & ---                                  & $9.80 \cdot 10^{-3}$                   & ---                                    & $2.37 \cdot 10^{-5}$ & ---   & 135.38 & $8.94 \cdot 10^{-5}$ & ---   & 109.61 \\
			40  & $7.96 \cdot 10^{-4}$                 & 2.01                                 & $2.35 \cdot 10^{-3}$                   & 2.06                                   & $5.53 \cdot 10^{-6}$ & 2.10  & 143.75 & $1.89 \cdot 10^{-5}$ & 2.24  & 124.54 \\
			80  & $1.99 \cdot 10^{-4}$                 & 2.00                                 & $5.82 \cdot 10^{-4}$                   & 2.01                                   & $1.36 \cdot 10^{-6}$ & 2.02  & 145.47 & $4.53 \cdot 10^{-6}$ & 2.06  & 128.58 \\
			160 & $4.96 \cdot 10^{-5}$                 & 2.00                                 & $1.45 \cdot 10^{-4}$                   & 2.00                                   & $3.39 \cdot 10^{-7}$ & 2.01  & 146.20 & $1.12 \cdot 10^{-6}$ & 2.02  & 129.69 \\
			320 & $1.24 \cdot 10^{-5}$                 & 2.00                                 & $3.63 \cdot 10^{-5}$                   & 2.00                                   & $8.46 \cdot 10^{-8}$ & 2.00  & 146.56 & $2.79 \cdot 10^{-7}$ & 2.00  & 129.97 \\
			\bottomrule
		\end{tabular}
	}

	\hspace{\fill}

	\vspace{12pt}

	\subfloat[Error with a basis made of three elements: $q=2$.]{%
		\begin{tabular}{ccccccccccc}
			\toprule
			    & \multicolumn{2}{c}{$h$, basis $V_h$} & \multicolumn{2}{c}{$Q$, basis $V_h$} & \multicolumn{3}{c}{$h$, basis $V_h^+$} & \multicolumn{3}{c}{$Q$, basis $V_h^+$}                                                                               \\
			\cmidrule(lr){2-3}\cmidrule(lr){4-5}\cmidrule(lr){6-8}\cmidrule(lr){9-11}
			$K$ & error                                & order                                & error                                  & order                                  & error                & order & gain  & error                & order & gain  \\
			\cmidrule(lr){1-1}\cmidrule(lr){2-3}\cmidrule(lr){4-5}\cmidrule(lr){6-8}\cmidrule(lr){9-11}
			20  & $3.06 \cdot 10^{-4}$                 & ---                                  & $1.23 \cdot 10^{-3}$                   & ---                                    & $3.90 \cdot 10^{-6}$ & ---   & 78.49 & $1.39 \cdot 10^{-5}$ & ---   & 88.29 \\
			40  & $4.20 \cdot 10^{-5}$                 & 2.86                                 & $1.83 \cdot 10^{-4}$                   & 2.75                                   & $5.87 \cdot 10^{-7}$ & 2.73  & 71.56 & $2.46 \cdot 10^{-6}$ & 2.50  & 74.27 \\
			80  & $5.44 \cdot 10^{-6}$                 & 2.95                                 & $2.43 \cdot 10^{-5}$                   & 2.91                                   & $8.10 \cdot 10^{-8}$ & 2.86  & 67.24 & $3.66 \cdot 10^{-7}$ & 2.75  & 66.26 \\
			160 & $6.88 \cdot 10^{-7}$                 & 2.98                                 & $3.09 \cdot 10^{-6}$                   & 2.98                                   & $1.05 \cdot 10^{-8}$ & 2.95  & 65.54 & $4.86 \cdot 10^{-8}$ & 2.91  & 63.52 \\
			320 & $8.62 \cdot 10^{-8}$                 & 3.00                                 & $3.88 \cdot 10^{-7}$                   & 2.99                                   & $1.33 \cdot 10^{-9}$ & 2.99  & 65.05 & $6.18 \cdot 10^{-9}$ & 2.98  & 62.76 \\
			\bottomrule
		\end{tabular}
	}
	\caption{%
		Shallow water system,
		Gaussian topography
		\eqref{eq:Gaussian_topography}:
		errors, orders of accuracy,
		and gain obtained when approximating
		a subcritical steady solution
		for bases with and without prior.
	}
	\label{tab:SW_steady_subcritical_Gaussian}
\end{table}

As a conclusion of this first test case,
we observe that using a Gaussian topography compared to a compactly supported topography
leads to a more stable order of accuracy, but with lower gains,
except for small values of~$K$ where the compactly supported topography
is not well-approximated.
The most important point is that the Gaussian topography requires a lower order quadrature to converge.
These results are in line with~\cite{SpiHuyDeB2015}.
As a consequence, we use the Gaussian topography in the remainder of this section.

Like in the previous section, we now consider $10^3$ parameters in $\mathbb{P}$,
and we compute the minimum, average and maximum gains for
$q \in \{0, 1, 2\}$.
To that end, we take $K = 20$ discretization cells.
The results are reported in \cref{tab:SW_steady_subcritical_stats},
where we observe that the average gains are substantial,
whatever the value of $q$,
and that the minimum gain is always greater than $1$.

\begin{table}[!ht]
	\centering
	\begin{tabular}{ccS[table-format=2.2]S[table-format=2.2]cS[table-format=3.2]S[table-format=3.2]cS[table-format=4.2]S[table-format=4.2]}
		\toprule
		      &  & \multicolumn{2}{c}{minimum gain} &       & \multicolumn{2}{c}{average gain} &        & \multicolumn{2}{c}{maximum gain}                        \\
		\cmidrule(lr){3-4} \cmidrule(lr){6-7} \cmidrule(lr){9-10}
		{$q$} &  & {$h$}                            & {$Q$} &                                  & {$h$}  & {$Q$}                            &  & {$h$}   & {$Q$}   \\
		\cmidrule(lr){1-1} \cmidrule(lr){3-4} \cmidrule(lr){6-7} \cmidrule(lr){9-10}
		0     &  & 21.28                            & 17.40 &                                  & 309.84 & 269.59                           &  & 1562.20 & 1628.39 \\
		1     &  & 7.47                             & 5.47  &                                  & 161.16 & 129.90                           &  & 845.97  & 729.03  \\
		2     &  & 4.37                             & 5.02  &                                  & 96.54  & 102.36                           &  & 707.41  & 704.55  \\
		\bottomrule
	\end{tabular}
	\caption{%
		Shallow water system,
		Gaussian topography \eqref{eq:Gaussian_topography}:
		statistics of the gains obtained
		for the approximation of a subcritical steady solution
		in basis $V_h^+$ with respect to basis $V_h$.
	}
	\label{tab:SW_steady_subcritical_stats}
\end{table}

\subsubsection{Supercritical flow}

We now turn to a supercritical flow.
In this case, the remaining parameters
$h_0$ and $Q_0$ are taken such that:
\begin{equation}
	\label{eq:parameter_space_shallow_water_sup}
	0.5 \leq h_0 \leq 0.75
	\quad ; \qquad
	4 \leq Q_0 \leq 5.
\end{equation}
{\rb
The boundary conditions are enforced
using the same expression \eqref{eq:BC_PINN_SW_sub}
as in the subcritical case,
and $h = h_0$ is imposed at the boundaries.
}
We check the approximate preservation
of the steady solution
by taking the initial water height
equal to the steady solution.

The results are displayed in \cref{tab:SW_steady_supercritical},
and we note that the gains are in line with the subcritical case,
from \cref{tab:SW_steady_subcritical_Gaussian}.

\begin{table}[!ht]
	\centering
	\subfloat[Error with a basis made of one element: $q=0$.]{%
		\begin{tabular}{ccccccccccc}
			\toprule
			    & \multicolumn{2}{c}{$h$, basis $V_h$} & \multicolumn{2}{c}{$Q$, basis $V_h$} & \multicolumn{3}{c}{$h$, basis $V_h^+$} & \multicolumn{3}{c}{$Q$, basis $V_h^+$}                                                                                 \\
			\cmidrule(lr){2-3}\cmidrule(lr){4-5}\cmidrule(lr){6-8}\cmidrule(lr){9-11}
			$K$ & error                                & order                                & error                                  & order                                  & error                & order & gain   & error                & order & gain   \\
			\cmidrule(lr){1-1}\cmidrule(lr){2-3}\cmidrule(lr){4-5}\cmidrule(lr){6-8}\cmidrule(lr){9-11}
			20  & $1.25 \cdot 10^{-2}$                 & ---                                  & $4.49 \cdot 10^{-2}$                   & ---                                    & $2.20 \cdot 10^{-5}$ & ---   & 566.64 & $7.25 \cdot 10^{-5}$ & ---   & 619.66 \\
			40  & $8.37 \cdot 10^{-3}$                 & 0.58                                 & $3.21 \cdot 10^{-2}$                   & 0.48                                   & $1.54 \cdot 10^{-5}$ & 0.51  & 542.09 & $5.17 \cdot 10^{-5}$ & 0.49  & 621.10 \\
			80  & $5.03 \cdot 10^{-3}$                 & 0.74                                 & $2.02 \cdot 10^{-2}$                   & 0.67                                   & $9.46 \cdot 10^{-6}$ & 0.71  & 531.46 & $3.20 \cdot 10^{-5}$ & 0.69  & 630.14 \\
			160 & $2.80 \cdot 10^{-3}$                 & 0.84                                 & $1.16 \cdot 10^{-2}$                   & 0.80                                   & $5.18 \cdot 10^{-6}$ & 0.87  & 540.85 & $1.81 \cdot 10^{-5}$ & 0.82  & 638.36 \\
			320 & $1.49 \cdot 10^{-3}$                 & 0.91                                 & $6.23 \cdot 10^{-3}$                   & 0.89                                   & $2.81 \cdot 10^{-6}$ & 0.88  & 529.25 & $1.04 \cdot 10^{-5}$ & 0.80  & 599.95 \\
			\bottomrule
		\end{tabular}
	}

	\hspace{\fill}

	\vspace{12pt}

	\subfloat[Error with a basis made of two elements: $q=1$.]{%
		\begin{tabular}{ccccccccccc}
			\toprule
			    & \multicolumn{2}{c}{$h$, basis $V_h$} & \multicolumn{2}{c}{$Q$, basis $V_h$} & \multicolumn{3}{c}{$h$, basis $V_h^+$} & \multicolumn{3}{c}{$Q$, basis $V_h^+$}                                                                                 \\
			\cmidrule(lr){2-3}\cmidrule(lr){4-5}\cmidrule(lr){6-8}\cmidrule(lr){9-11}
			$K$ & error                                & order                                & error                                  & order                                  & error                & order & gain   & error                & order & gain   \\
			\cmidrule(lr){1-1}\cmidrule(lr){2-3}\cmidrule(lr){4-5}\cmidrule(lr){6-8}\cmidrule(lr){9-11}
			20  & $5.32 \cdot 10^{-4}$                 & ---                                  & $1.96 \cdot 10^{-3}$                   & ---                                    & $4.96 \cdot 10^{-6}$ & ---   & 107.22 & $1.57 \cdot 10^{-5}$ & ---   & 124.92 \\
			40  & $1.16 \cdot 10^{-4}$                 & 2.19                                 & $4.50 \cdot 10^{-4}$                   & 2.12                                   & $8.11 \cdot 10^{-7}$ & 2.61  & 143.39 & $3.21 \cdot 10^{-6}$ & 2.29  & 140.36 \\
			80  & $2.80 \cdot 10^{-5}$                 & 2.05                                 & $1.11 \cdot 10^{-4}$                   & 2.03                                   & $1.67 \cdot 10^{-7}$ & 2.28  & 167.57 & $7.30 \cdot 10^{-7}$ & 2.14  & 151.60 \\
			160 & $6.93 \cdot 10^{-6}$                 & 2.02                                 & $2.76 \cdot 10^{-5}$                   & 2.01                                   & $3.97 \cdot 10^{-8}$ & 2.07  & 174.56 & $1.78 \cdot 10^{-7}$ & 2.04  & 154.79 \\
			320 & $1.73 \cdot 10^{-6}$                 & 2.00                                 & $6.88 \cdot 10^{-6}$                   & 2.00                                   & $9.82 \cdot 10^{-9}$ & 2.02  & 175.99 & $4.42 \cdot 10^{-8}$ & 2.01  & 155.59 \\
			\bottomrule
		\end{tabular}
	}

	\hspace{\fill}

	\vspace{12pt}

	\subfloat[Error with a basis made of three elements: $q=2$.]{%
		\begin{tabular}{ccccccccccc}
			\toprule
			    & \multicolumn{2}{c}{$h$, basis $V_h$} & \multicolumn{2}{c}{$Q$, basis $V_h$} & \multicolumn{3}{c}{$h$, basis $V_h^+$} & \multicolumn{3}{c}{$Q$, basis $V_h^+$}                                                                                                                                                       \\
			\cmidrule(lr){2-3}\cmidrule(lr){4-5}\cmidrule(lr){6-8}\cmidrule(lr){9-11}
			$K$ & error                                & order                                & error                                  & order                                  & error                                                          & order & gain                                & error                & order & gain  \\
			\cmidrule(lr){1-1}\cmidrule(lr){2-3}\cmidrule(lr){4-5}\cmidrule(lr){6-8}\cmidrule(lr){9-11}
			20  & $8.33 \cdot 10^{-5}$                 & ---                                  & $2.50 \cdot 10^{-4}$                   & ---                                    & \tablenum[table-format=1.2e2,exponent-product=\cdot]{6.98e-7}  & ---   & \tablenum[table-format=3.2]{119.47} & $2.59 \cdot 10^{-6}$ & ---   & 96.70 \\
			40  & $1.33 \cdot 10^{-5}$                 & 2.64                                 & $3.85 \cdot 10^{-5}$                   & 2.70                                   & \tablenum[table-format=1.2e2,exponent-product=\cdot]{1.45e-7}  & 2.27  & \tablenum[table-format=3.2]{92.12}  & $4.45 \cdot 10^{-7}$ & 2.54  & 86.58 \\
			80  & $1.83 \cdot 10^{-6}$                 & 2.87                                 & $5.21 \cdot 10^{-6}$                   & 2.89                                   & \tablenum[table-format=1.2e2,exponent-product=\cdot]{2.47e-8}  & 2.55  & \tablenum[table-format=3.2]{73.83}  & $7.19 \cdot 10^{-8}$ & 2.63  & 72.48 \\
			160 & $2.34 \cdot 10^{-7}$                 & 2.96                                 & $6.67 \cdot 10^{-7}$                   & 2.97                                   & \tablenum[table-format=1.2e2,exponent-product=\cdot]{3.50e-9}  & 2.82  & \tablenum[table-format=3.2]{67.04}  & $1.00 \cdot 10^{-8}$ & 2.84  & 66.59 \\
			320 & $2.95 \cdot 10^{-8}$                 & 2.99                                 & $8.39 \cdot 10^{-8}$                   & 2.99                                   & \tablenum[table-format=1.2e2,exponent-product=\cdot]{4.54e-10} & 2.95  & \tablenum[table-format=3.2]{65.00}  & $1.30 \cdot 10^{-9}$ & 2.95  & 64.74 \\
			\bottomrule
		\end{tabular}
	}
	\caption{%
		Shallow water system,
		Gaussian topography
		\eqref{eq:Gaussian_topography}:
		errors, orders of accuracy,
		and gain obtained when approximating
		a supercritical steady solution
		for bases with and without prior.
	}
	\label{tab:SW_steady_supercritical}
\end{table}

Furthermore, in \cref{tab:SW_steady_supercritical_stats},
we display some statistics on the gains obtained by using the prior,
in the same configuration as for the subcritical regime.
We draw similar conclusions to the subcritical case.

\begin{table}[!ht]
	\centering
	\begin{tabular}{ccS[table-format=2.2]S[table-format=2.2]cS[table-format=3.2]S[table-format=3.2]cS[table-format=4.2]S[table-format=4.2]}
		\toprule
		      &  & \multicolumn{2}{c}{minimum gain} &       & \multicolumn{2}{c}{average gain} &        & \multicolumn{2}{c}{maximum gain}                        \\
		\cmidrule(lr){3-4} \cmidrule(lr){6-7} \cmidrule(lr){9-10}
		{$q$} &  & {$h$}                            & {$Q$} &                                  & {$h$}  & {$Q$}                            &  & {$h$}   & {$Q$}   \\
		\cmidrule(lr){1-1} \cmidrule(lr){3-4} \cmidrule(lr){6-7} \cmidrule(lr){9-10}
		0     &  & 19.83                            & 23.50 &                                  & 309.13 & 314.36                           &  & 1789.56 & 1923.34 \\
		1     &  & 5.36                             & 5.54  &                                  & 111.41 & 120.11                           &  & 354.89  & 376.47  \\
		2     &  & 7.29                             & 7.18  &                                  & 123.58 & 104.49                           &  & 468.92  & 381.27  \\
		\bottomrule
	\end{tabular}
	\caption{%
		Shallow water system,
		Gaussian topography \eqref{eq:Gaussian_topography}:
		statistics of the gains obtained
		for the approximation of a supercritical steady solution
		in basis $V_h^+$ with respect to basis $V_h$.
	}
	\label{tab:SW_steady_supercritical_stats}
\end{table}

\subsubsection{Transcritical flow}

The last steady experiment we study is the preservation of a
transcritical steady solution.
Such steady solutions are
significantly harder to capture.
Indeed, when $\text{Fr} = 1$,
the steady ODE \eqref{eq:shallow_water_steady}
yields $\partial_x Z = 0$, and therefore
the derivative of the steady water height
is not defined using only \eqref{eq:shallow_water_steady}.
This is a well-known issue when approximating
transcritical steady solutions,
see for instance \cite{CasCha2016,GomCasParRus2021}.
{\rb
In this case, the Froude number is known to be equal to $1$
at the top of the topography bump,
i.e. at $x = 1/2$,
where $\partial_x Z = 0$.
At this location, we know that
the water height becomes equal to
\smash{$h_c({\mu}) = Q_0^{2/3} g^{-1/3}$}.
Since fixing the discharge $Q_0$ imposes a
fixed value of the water height at $x = 1/2$,
this means that $h_0$ is no longer a degree of freedom in this transcritical case.
For this regime,
we choose $2 \leq Q_0 \leq 3$
and we take $0.75 \leq \alpha \leq 1.25$.

Then, to ensure a correct treatment of the boundary conditions, we take
\begin{equation*}
	\label{eq:BC_PINN_SW_trans}
	\widetilde{h}_\theta(x; {\mu})
	=
	\left(
        h_R({\mu}) + \left(
        1 - \tanh\left(15 \left(x - \frac 1 2\right)\right)
        \right) \frac {h_L({\mu}) - h_R({\mu})} 2
    \right) + Z(x; \alpha, \beta) \ h_\theta(x; {\mu}).
\end{equation*}
Since $Z$ is (very) close to $0$ at the boundaries,
this expression ensures that
$\smash{\widetilde{h}_\theta(x; {\mu})} \simeq h_R({\mu})$ at the right boundary, and that
$\smash{\widetilde{h}_\theta(x; {\mu})} \simeq h_L({\mu})$ at the left boundary.
Hence, $h_L({\mu})$ and $h_R({\mu})$ represent the left and right boundary conditions,
which need to be computed according to the flow characteristics.
}
Since we consider a smooth steady solution,
relations \eqref{eq:shallow_water_steady} lead to
\begin{equation*}
	E(h, x; {\mu}) \coloneqq
	\frac {Q_0^2}{2 h^2} + g \big(h + Z(x; \alpha, \beta)\big) =
	\text{constant}.
\end{equation*}
{\rb
We assume, without loss of generality, that $h_L({\mu}) > h_R({\mu})$.
This corresponds to a subcritical flow on the left and a supercritical flow on the right,
on either side of the topography bump.
Denoting by $E_c(\mu) = E(h_c({\mu}), 1/2, x)$ the value of $E$ at the critical point,
$h_L({\mu})$ and $h_R({\mu})$ are then solution to the following nonlinear equations:
\begin{equation*}
	\frac {Q_0^2}{2 h_L({\mu})^2} + g (h_L({\mu}) + Z(0; \alpha, \beta))
    =
    E_c(\mu)
    \text{\qquad and \qquad}
	\frac {Q_0^2}{2 h_R({\mu})^2} + g (h_R({\mu}) + Z(1; \alpha, \beta))
    =
    E_c(\mu).
\end{equation*}
Each of these equations has two solutions, and we choose the one that corresponds
to the subcritical flow ($\text{Fr} < 1$) on the left
and the supercritical flow ($\text{Fr} > 1$) on the right.
}

\cref{tab:SW_steady_transcritical} contains
the errors, the orders of convergence and the gains.
We observe that the gains are lower than in the other two cases,
but that was to be expected since the transcritical solution
comes from a singular ODE,
and it is harder for the PINN to approximate its solutions.

\begin{table}[!ht]
	\centering
	\subfloat[Error with a basis made of one element: $q=0$.]{%
		\begin{tabular}{ccccccccccc}
			\toprule
			    & \multicolumn{2}{c}{$h$, basis $V_h$} & \multicolumn{2}{c}{$Q$, basis $V_h$} & \multicolumn{3}{c}{$h$, basis $V_h^+$} & \multicolumn{3}{c}{$Q$, basis $V_h^+$}                                                                                                              \\
			\cmidrule(lr){2-3}\cmidrule(lr){4-5}\cmidrule(lr){6-8}\cmidrule(lr){9-11}
			$K$ & error                                & order                                & error                                  & order                                  & error                & order & gain                                & error                & order & gain   \\
			\cmidrule(lr){1-1}\cmidrule(lr){2-3}\cmidrule(lr){4-5}\cmidrule(lr){6-8}\cmidrule(lr){9-11}
			40  & $4.81 \cdot 10^{-2}$                 & ---                                  & $4.29 \cdot 10^{-2}$                   & ---                                    & $1.79 \cdot 10^{-4}$ & ---   & \tablenum[table-format=3.2]{268.84} & $2.10 \cdot 10^{-4}$ & ---   & 204.54 \\
			80  & $2.58 \cdot 10^{-2}$                 & 0.90                                 & $2.55 \cdot 10^{-2}$                   & 0.75                                   & $1.37 \cdot 10^{-4}$ & 0.39  & \tablenum[table-format=3.2]{189.00} & $1.53 \cdot 10^{-4}$ & 0.45  & 165.82 \\
			160 & $1.34 \cdot 10^{-2}$                 & 0.94                                 & $1.40 \cdot 10^{-2}$                   & 0.86                                   & $9.50 \cdot 10^{-5}$ & 0.52  & \tablenum[table-format=3.2]{141.21} & $1.00 \cdot 10^{-4}$ & 0.61  & 139.28 \\
			320 & $6.84 \cdot 10^{-3}$                 & 0.97                                 & $7.35 \cdot 10^{-3}$                   & 0.93                                   & $6.00 \cdot 10^{-5}$ & 0.66  & \tablenum[table-format=3.2]{114.00} & $6.03 \cdot 10^{-5}$ & 0.74  & 121.99 \\
			640 & $3.46 \cdot 10^{-3}$                 & 0.98                                 & $3.77 \cdot 10^{-3}$                   & 0.96                                   & $3.56 \cdot 10^{-5}$ & 0.75  & \tablenum[table-format=3.2]{97.20}  & $3.40 \cdot 10^{-5}$ & 0.82  & 110.81 \\
			\bottomrule
		\end{tabular}
	}

	\hspace{\fill}

	\vspace{12pt}

	\subfloat[Error with a basis made of two elements: $q=1$.]{%
		\begin{tabular}{ccccccccccc}
			\toprule
			    & \multicolumn{2}{c}{$h$, basis $V_h$} & \multicolumn{2}{c}{$Q$, basis $V_h$} & \multicolumn{3}{c}{$h$, basis $V_h^+$} & \multicolumn{3}{c}{$Q$, basis $V_h^+$}                                                                               \\
			\cmidrule(lr){2-3}\cmidrule(lr){4-5}\cmidrule(lr){6-8}\cmidrule(lr){9-11}
			$K$ & error                                & order                                & error                                  & order                                  & error                & order & gain  & error                & order & gain  \\
			\cmidrule(lr){1-1}\cmidrule(lr){2-3}\cmidrule(lr){4-5}\cmidrule(lr){6-8}\cmidrule(lr){9-11}
			40  & $6.69 \cdot 10^{-4}$                 & ---                                  & $6.15 \cdot 10^{-4}$                   & ---                                    & $1.85 \cdot 10^{-5}$ & ---   & 36.18 & $1.16 \cdot 10^{-5}$ & ---   & 52.90 \\
			80  & $1.67 \cdot 10^{-4}$                 & 2.00                                 & $1.53 \cdot 10^{-4}$                   & 2.01                                   & $3.11 \cdot 10^{-6}$ & 2.57  & 53.69 & $2.11 \cdot 10^{-6}$ & 2.46  & 72.26 \\
			160 & $4.17 \cdot 10^{-5}$                 & 2.00                                 & $3.81 \cdot 10^{-5}$                   & 2.00                                   & $6.77 \cdot 10^{-7}$ & 2.20  & 61.66 & $4.36 \cdot 10^{-7}$ & 2.28  & 87.42 \\
			320 & $1.04 \cdot 10^{-5}$                 & 2.00                                 & $9.53 \cdot 10^{-6}$                   & 2.00                                   & $1.65 \cdot 10^{-7}$ & 2.04  & 63.28 & $1.05 \cdot 10^{-7}$ & 2.06  & 91.04 \\
			640 & $2.61 \cdot 10^{-6}$                 & 2.00                                 & $2.38 \cdot 10^{-6}$                   & 2.00                                   & $4.10 \cdot 10^{-8}$ & 2.01  & 63.67 & $2.59 \cdot 10^{-8}$ & 2.01  & 91.90 \\
			\bottomrule
		\end{tabular}
	}

	\hspace{\fill}

	\vspace{12pt}

	\subfloat[Error with a basis made of three elements: $q=2$.]{%
		\begin{tabular}{ccccccccccc}
			\toprule
			    & \multicolumn{2}{c}{$h$, basis $V_h$} & \multicolumn{2}{c}{$Q$, basis $V_h$} & \multicolumn{3}{c}{$h$, basis $V_h^+$} & \multicolumn{3}{c}{$Q$, basis $V_h^+$}                                                                                                            \\
			\cmidrule(lr){2-3}\cmidrule(lr){4-5}\cmidrule(lr){6-8}\cmidrule(lr){9-11}
			$K$ & error                                & order                                & error                                  & order                                  & error                & order & gain                               & error                & order & gain  \\
			\cmidrule(lr){1-1}\cmidrule(lr){2-3}\cmidrule(lr){4-5}\cmidrule(lr){6-8}\cmidrule(lr){9-11}
			40  & $9.76 \cdot 10^{-5}$                 & ---                                  & $7.21 \cdot 10^{-5}$                   & ---                                    & $2.95 \cdot 10^{-6}$ & ---   & \tablenum[table-format=2.2]{33.10} & $1.81 \cdot 10^{-6}$ & ---   & 39.88 \\
			80  & $1.91 \cdot 10^{-5}$                 & 2.35                                 & $1.19 \cdot 10^{-5}$                   & 2.60                                   & $6.89 \cdot 10^{-7}$ & 2.10  & \tablenum[table-format=2.2]{27.74} & $3.86 \cdot 10^{-7}$ & 2.23  & 30.84 \\
			160 & $3.25 \cdot 10^{-6}$                 & 2.56                                 & $1.80 \cdot 10^{-6}$                   & 2.72                                   & $1.66 \cdot 10^{-7}$ & 2.05  & \tablenum[table-format=2.2]{19.54} & $8.68 \cdot 10^{-8}$ & 2.15  & 20.76 \\
			320 & $5.01 \cdot 10^{-7}$                 & 2.70                                 & $2.51 \cdot 10^{-7}$                   & 2.84                                   & $3.82 \cdot 10^{-8}$ & 2.12  & \tablenum[table-format=2.2]{13.10} & $1.78 \cdot 10^{-8}$ & 2.29  & 14.16 \\
			640 & $7.42 \cdot 10^{-8}$                 & 2.76                                 & $3.34 \cdot 10^{-8}$                   & 2.91                                   & $7.74 \cdot 10^{-9}$ & 2.30  & \tablenum[table-format=2.2]{9.58}  & $3.13 \cdot 10^{-9}$ & 2.50  & 10.66 \\
			\bottomrule
		\end{tabular}
	}
	\caption{%
		Shallow water system,
		Gaussian topography
		\eqref{eq:Gaussian_topography}:
		errors, orders of accuracy,
		and gain obtained when approximating
		a transcritical steady solution
		for bases with and without prior.
	}
	\label{tab:SW_steady_transcritical}
\end{table}

Finally, we report in \cref{tab:SW_steady_transcritical_stats}
the minimum, average and maximum gains
obtained by using the basis $V_h^+$ instead of the basis $V_h$.
We draw the same conclusions as in the other two regimes,
even though the gains are, on average, lower.
This was to be expected, since the transcritical regime
is harder to capture than the subcritical and supercritical ones,
and therefore that the prior is of lower quality.
Nevertheless, the gains remain substantial for all values of $q$.

\begin{table}[!ht]
	\centering
	\begin{tabular}{ccS[table-format=2.2]S[table-format=2.2]cS[table-format=3.2]S[table-format=3.2]cS[table-format=3.2]S[table-format=3.2]}
		\toprule
		      &  & \multicolumn{2}{c}{minimum gain} &       & \multicolumn{2}{c}{average gain} &        & \multicolumn{2}{c}{maximum gain}                      \\
		\cmidrule(lr){3-4} \cmidrule(lr){6-7} \cmidrule(lr){9-10}
		{$q$} &  & {$h$}                            & {$Q$} &                                  & {$h$}  & {$Q$}                            &  & {$h$}  & {$Q$}  \\
		\cmidrule(lr){1-1} \cmidrule(lr){3-4} \cmidrule(lr){6-7} \cmidrule(lr){9-10}
		0     &  & 35.82                            & 26.19 &                                  & 254.53 & 177.02                           &  & 928.03 & 668.73 \\
		1     &  & 5.51                             & 4.73  &                                  & 30.83  & 38.69                            &  & 134.83 & 142.11 \\
		2     &  & 4.55                             & 6.16  &                                  & 16.49  & 24.29                            &  & 96.95  & 109.94 \\
		\bottomrule
	\end{tabular}
	\caption{%
		Shallow water system,
		Gaussian topography \eqref{eq:Gaussian_topography}:
		statistics of the gains obtained
		for the approximation of a transcritical steady solution
		in basis $V_h^+$ with respect to basis $V_h$.
	}
	\label{tab:SW_steady_transcritical_stats}
\end{table}

\subsubsection{Perturbation of a steady flow}

This last experiment related to the shallow water equations
concerns a perturbed steady flow.
We only perform this study on the subcritical flow,
but the other regimes behave the same.
We take $\smash{\varepsilon \in \{ 5 \cdot 10^{-k} \}_{k \in \{1, 2, 3\}}}$
and $20$ space cells,
and set the initial water height to
$h(0, x; \mu) = \big( 1 + \varepsilon \sin (2 \pi x) \big) \, h_\text{eq}(x; \mu)$.
{\rb The PINN is the same as in \cref{sec:shallow_water_subcritical}.}

\begin{figure}[!ht]
	\centering
	\makebox[\textwidth][c]{\includegraphics{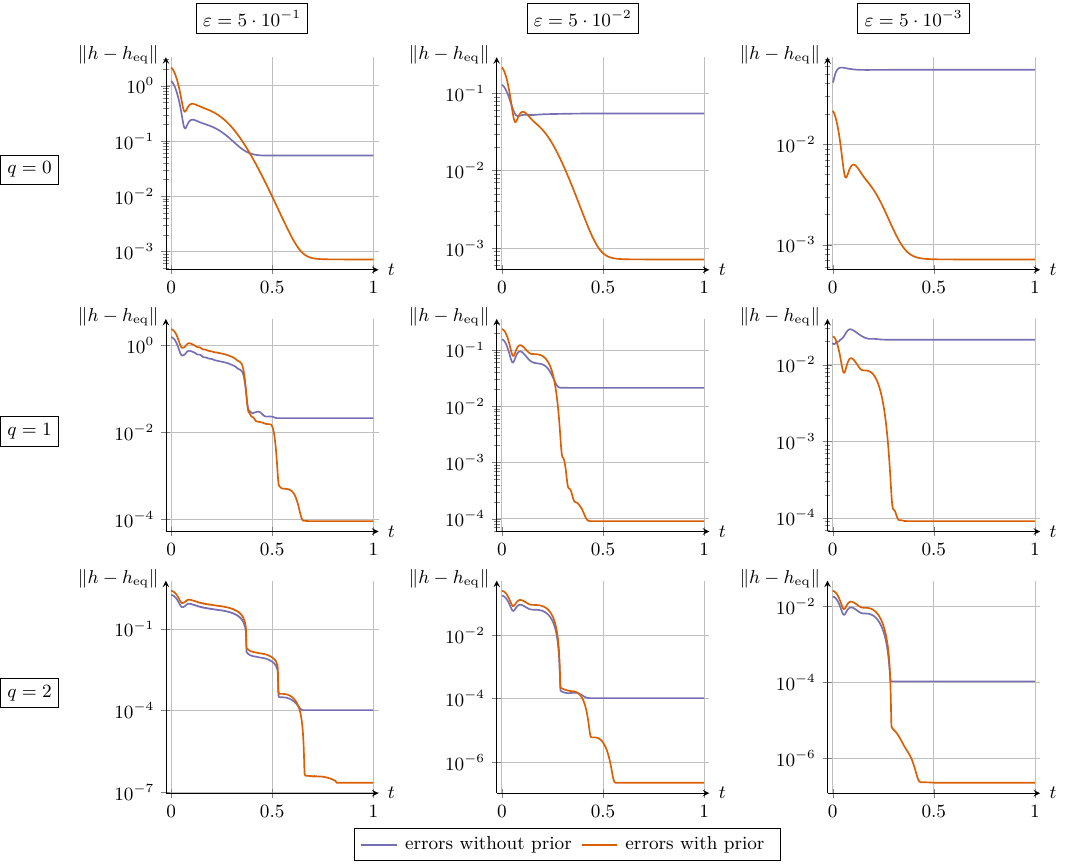}}
	\caption{%
		Shallow water equations, compactly supported topography:
        {\rb distance between the DG solution $h$
        and the underlying steady solution $h_\text{eq}$},
		with respect to time,
		for the approximation of a perturbed subcritical steady solution
		for bases with and without prior.
	}
	\label{fig:SW_subcritical_h_perturbed}
\end{figure}

{\rb The errors between the DG approximation of $h$
and the underlying steady solution $h_\text{eq}$}
with respect to time are displayed in
\cref{fig:SW_subcritical_h_perturbed},
until the final physical time $T = 1$.
Like in \cref{sec:transport}, with the prior,
the error decreases to a much lower level than without the prior.
This good behavior was expected since the prior
makes it possible for the enhanced DG scheme
to achieve higher accuracy on steady solutions.

\subsection{Euler-Poisson equations in spherical geometry}
\label{sec:euler_poisson}

We now consider the Euler-Poisson equations in spherical geometry.
This system is used in astrophysics, for instance,
where it serves to model stars held together by gravitation,
see e.g. \cite{Cha1967,CriBal2018,Kae2022}.
They are given by
\begin{equation}
	\label{eq:euler_poisson}
	\begin{dcases}
		\partial_t \rho + \partial_r Q
		=
		- \frac 2 r Q,                                          \\
		\partial_t Q + \partial_r \left(\frac {Q^2} \rho + p \right)
		=
		- \frac 2 r \frac {Q^2} \rho - \rho \partial_r \phi,    \\
		\partial_t E + \partial_r \left(\frac {Q} \rho (E + p)\right)
		=
		- \frac 2 r \frac {Q} \rho (E + p) - Q \partial_r \phi, \\
		{\rb \frac 1 {r^2} \partial_r (r^2 \partial_r \phi) = 4 \pi G \rho},
	\end{dcases}
\end{equation}
where $G$ is a gravity constant, fixed to $G = 1$ in our applications,
and where we take $p$ as a function of $\rho$, $Q$ and~$E$
through a pressure law to be specified.
In \eqref{eq:euler_poisson}, $\rho$ is the density,
$Q$ is the momentum, $E$ is the energy,
and $\phi$ is the gravitational potential.
Unless otherwise mentioned, the boundary conditions are of Dirichlet type,
with the value of the steady solution prescribed of the boundaries.

The space domain is $r \in (0, 1)$.
The apparent singularity at $r = 0$
is resolved by imposing suitable boundary conditions,
namely $\rho(0) = 1$ and $\partial_r \rho(0)$ given according to the pressure law.
Indeed, the assumption that there is no gravity at $r = 0$
leads to $\partial_r p(0) = 0$,
which makes it possible to determine $\partial_r \rho(0)$.
For more information on the boundary conditions
and on the DG discretization of \eqref{eq:euler_poisson},
the reader is referred to \cite{ZhaXinEnd2022}.

The steady solutions at rest are given by
\begin{equation*}
	\begin{cases}
		Q = 0,                                   \\
		\partial_r p + \rho \partial_r \phi = 0, \\
		{\rb \partial_r(r^2 \partial_r \phi) = 4 \pi r^2 G \rho}.
	\end{cases}
\end{equation*}
For the steady solutions,
we shall distinguish two cases for the pressure law:
a polytropic pressure law,
and a temperature-dependent pressure law.

\subsubsection{Polytropic pressure law}

In this case, we introduce two parameters $\kappa$ and $\gamma$,
so the parameter vector ${\mu}$ is composed of two elements:
\begin{equation*}
	{\mu} =
	\begin{pmatrix}
		\kappa \\
		\gamma \\
	\end{pmatrix}
	\in \mathbb{P} \subset \mathbb{R}^2,
	\quad \kappa \in \mathbb{R}_+,
	\quad \gamma \in (1, +\infty).
\end{equation*}
Equipped with this parameter vector, we define the polytropic pressure law
\begin{equation*}
	p(\rho; {\mu}) = \kappa \rho^\gamma,
\end{equation*}
and the steady solutions are then given as solutions to
the following nonlinear second-order ordinary differential equation:
\begin{equation*}
	\frac{d}{dr} \left(
	r^2 \kappa \gamma \rho^{\gamma - 2}
	\frac{d \rho}{d r}
	\right)
	= 4 \pi r^2 G \rho.
\end{equation*}
In general, this ODE does not have analytic solutions.
However, it turns out that,
for specific values of $\gamma$,
there exists an analytic solution to this ODE.
For instance, with $\gamma = 2$, we obtain
\begin{equation*}
	\rho(r) = \frac{\sin(\alpha r)}{\alpha r},
	\text{\qquad with \quad}
	\alpha = \sqrt{\frac{2 \pi G} \kappa}.
\end{equation*}
Regarding the boundary conditions,
the condition $\partial_r p(0) = 0$ leads to $\partial_r \rho(0) = 0$
for this pressure law.
In this case, we take $\lambda = 1 + \sqrt{\gamma}$
in \eqref{eq:time_step} to compute the time step $\Delta t$.

To obtain a prior $\rho_\theta$, as usual,
we train a PINN with about $1400$ trainable parameters on $7$ fully connected layers.
The boundary conditions are taken into account by setting
\begin{equation*}
	\widetilde{\rho}_\theta(r; {\mu}) =
	1 + r^2 \rho_\theta(r; {\mu}),
\end{equation*}
where $\rho_\theta$ is the result of the PINN.
{\rb
This expression ensures that $\widetilde{\rho}_\theta(0; {\mu}) = 1$
and $\partial_r \widetilde{\rho}_\theta(0; {\mu}) = 0$.
}
The PINN is trained on the parameter space
\begin{equation}
	\label{eq:parameter_space_euler_polytropic}
	\mathbb{P} = [2, 5] \times [1.5, 3.5],
\end{equation}
with only the physics-based loss function
corresponding to the steady solution:
\begin{equation*}
	\mathcal{J}(\theta) = \left\|
	\frac{d}{dr} \left(
	r^2 \kappa \gamma \widetilde{\rho}_\theta^{\gamma - 2}
	\frac{d \widetilde{\rho}_\theta}{d r}
	\right)
	- 4 \pi r^2 G \widetilde{\rho}_\theta
	\right\|.
\end{equation*}
In addition, the prior for $q$ is set to $Q_\theta = 1$
since we wish to approximate a constant momentum.
Finally, the prior for $E$ is set to
$E_\theta = p(\widetilde{\rho}_\theta; {\mu}) / (\gamma - 1)$.
Training takes about 5 minutes on a dual NVIDIA K80 GPU,
until the loss function is equal to about $5 \cdot 10^{-5}$.
In the DG discretization,
the degree of the quadrature formula is the usual $n_q = q + 2$:
there is no need to further increase the order
of the quadrature rule in this case.

We first collect, in \cref{tab:euler_steady_polytropic},
the results of the approximation in both bases (with and without prior),
for $\kappa = 2$ and $\gamma = 2.5$, and until the final time $T = 0.01$.
As usual, the observed gain is larger for smaller number of basis elements.
We observe a slight superconvergence on the momentum $Q$
when using the prior with $q = 0$.
For these values of $\kappa$ and $\gamma$,
gains on the density are not very large for $q = 2$,
but this is compensated by larger gains on the energy.

\begin{table}[ht!]
	\centering
	\subfloat[Errors with a basis made of one element: $q=0$.]{%
		\makebox[0.9\textwidth][c]{
			\small
			\begin{tabular}{c@{\hspace{10pt}}c@{\hspace{1pt}}c@{\hspace{10pt}}c@{\hspace{1pt}}c@{\hspace{10pt}}c@{\hspace{1pt}}c@{\hspace{10pt}}c@{\hspace{1pt}}c@{\hspace{1pt}}c@{\hspace{10pt}}c@{\hspace{1pt}}c@{\hspace{1pt}}c@{\hspace{10pt}}c@{\hspace{1pt}}c@{\hspace{1pt}}c}
				\toprule
				    & \multicolumn{2}{c}{$\rho$, basis $V_h$}                           & \multicolumn{2}{c}{$Q$, basis $V_h$} & \multicolumn{2}{c}{$E$, basis $V_h$}                           & \multicolumn{3}{c}{$\rho$, basis $V_h^+$} & \multicolumn{3}{c}{$Q$, basis $V_h^+$}                         & \multicolumn{3}{c}{$E$, basis $V_h^+$}                                                                                                                                                                                                                                                                                   \\
				\cmidrule(lr){2-3}\cmidrule(lr){4-5}\cmidrule(lr){6-7}\cmidrule(lr){8-10} \cmidrule(lr){11-13} \cmidrule(lr){14-16}
				$K$ & error                                                          & order                                & error                                                          & order                                  & error                                                          & order                                  & error                                                          & order & gain   & error                                                          & order & gain                               & error                                                          & order & gain   \\
				\cmidrule(lr){2-3}\cmidrule(lr){4-5}\cmidrule(lr){6-7}\cmidrule(lr){8-10} \cmidrule(lr){11-13} \cmidrule(lr){14-16}
				10  & \tablenum[table-format=1.2e3,exponent-product=\cdot]{3.37e-02} & ---                                  & \tablenum[table-format=1.2e3,exponent-product=\cdot]{2.60e-03} & ---                                    & \tablenum[table-format=1.2e3,exponent-product=\cdot]{7.55e-02} & ---                                    & \tablenum[table-format=1.2e3,exponent-product=\cdot]{1.08e-04} & ---   & 312.50 & \tablenum[table-format=1.2e3,exponent-product=\cdot]{8.94e-04} & ---   & \tablenum[table-format=2.2]{2.91}  & \tablenum[table-format=1.2e3,exponent-product=\cdot]{3.43e-04} & ---   & 219.99 \\
				20  & \tablenum[table-format=1.2e3,exponent-product=\cdot]{1.69e-02} & 1.00                                 & \tablenum[table-format=1.2e3,exponent-product=\cdot]{1.51e-03} & 0.79                                   & \tablenum[table-format=1.2e3,exponent-product=\cdot]{3.78e-02} & 1.00                                   & \tablenum[table-format=1.2e3,exponent-product=\cdot]{6.49e-05} & 0.73  & 259.80 & \tablenum[table-format=1.2e3,exponent-product=\cdot]{3.34e-04} & 1.42  & \tablenum[table-format=2.2]{4.51}  & \tablenum[table-format=1.2e3,exponent-product=\cdot]{2.10e-04} & 0.71  & 180.21 \\
				40  & \tablenum[table-format=1.2e3,exponent-product=\cdot]{8.44e-03} & 1.00                                 & \tablenum[table-format=1.2e3,exponent-product=\cdot]{8.27e-04} & 0.87                                   & \tablenum[table-format=1.2e3,exponent-product=\cdot]{1.89e-02} & 1.00                                   & \tablenum[table-format=1.2e3,exponent-product=\cdot]{3.41e-05} & 0.93  & 247.75 & \tablenum[table-format=1.2e3,exponent-product=\cdot]{1.06e-04} & 1.65  & \tablenum[table-format=2.2]{7.78}  & \tablenum[table-format=1.2e3,exponent-product=\cdot]{1.11e-04} & 0.93  & 171.17 \\
				80  & \tablenum[table-format=1.2e3,exponent-product=\cdot]{4.22e-03} & 1.00                                 & \tablenum[table-format=1.2e3,exponent-product=\cdot]{4.60e-04} & 0.85                                   & \tablenum[table-format=1.2e3,exponent-product=\cdot]{9.46e-03} & 1.00                                   & \tablenum[table-format=1.2e3,exponent-product=\cdot]{1.72e-05} & 0.99  & 246.26 & \tablenum[table-format=1.2e3,exponent-product=\cdot]{3.24e-05} & 1.71  & \tablenum[table-format=2.2]{14.20} & \tablenum[table-format=1.2e3,exponent-product=\cdot]{5.57e-05} & 0.99  & 169.99 \\
				160 & \tablenum[table-format=1.2e3,exponent-product=\cdot]{2.11e-03} & 1.00                                 & \tablenum[table-format=1.2e3,exponent-product=\cdot]{2.59e-04} & 0.83                                   & \tablenum[table-format=1.2e3,exponent-product=\cdot]{4.73e-03} & 1.00                                   & \tablenum[table-format=1.2e3,exponent-product=\cdot]{8.17e-06} & 1.07  & 258.59 & \tablenum[table-format=1.2e3,exponent-product=\cdot]{9.15e-06} & 1.82  & \tablenum[table-format=2.2]{28.25} & \tablenum[table-format=1.2e3,exponent-product=\cdot]{2.64e-05} & 1.07  & 178.96 \\
				\bottomrule
			\end{tabular}
		}
	}
	\label{tab:euler_steady_polytropic_nG_1}

	\hspace{\fill}

	\vspace{12pt}

	\subfloat[Errors with a basis made of two elements: $q=1$.]{%
		\makebox[0.9\textwidth][c]{
			\small
			\begin{tabular}{c@{\hspace{10pt}}c@{\hspace{1pt}}c@{\hspace{10pt}}c@{\hspace{1pt}}c@{\hspace{10pt}}c@{\hspace{1pt}}c@{\hspace{10pt}}c@{\hspace{1pt}}c@{\hspace{1pt}}c@{\hspace{10pt}}c@{\hspace{1pt}}c@{\hspace{1pt}}c@{\hspace{10pt}}c@{\hspace{1pt}}c@{\hspace{1pt}}c}
				\toprule
				    & \multicolumn{2}{c}{$\rho$, basis $V_h$}                           & \multicolumn{2}{c}{$Q$, basis $V_h$} & \multicolumn{2}{c}{$E$, basis $V_h$}                           & \multicolumn{3}{c}{$\rho$, basis $V_h^+$} & \multicolumn{3}{c}{$Q$, basis $V_h^+$}                         & \multicolumn{3}{c}{$E$, basis $V_h^+$}                                                                                                                                                                                                                                                                                  \\
				\cmidrule(lr){2-3}\cmidrule(lr){4-5}\cmidrule(lr){6-7}\cmidrule(lr){8-10} \cmidrule(lr){11-13} \cmidrule(lr){14-16}
				$K$ & error                                                          & order                                & error                                                          & order                                  & error                                                          & order                                  & error                                                          & order & gain  & error                                                          & order & gain                                & error                                                          & order & gain  \\
				\cmidrule(lr){1-1}\cmidrule(lr){2-3}\cmidrule(lr){4-5}\cmidrule(lr){6-7}\cmidrule(lr){8-10} \cmidrule(lr){11-13} \cmidrule(lr){14-16}
				10  & \tablenum[table-format=1.2e3,exponent-product=\cdot]{1.03e-03} & ---                                  & \tablenum[table-format=1.2e3,exponent-product=\cdot]{1.30e-03} & ---                                    & \tablenum[table-format=1.2e3,exponent-product=\cdot]{2.23e-03} & ---                                    & \tablenum[table-format=1.2e3,exponent-product=\cdot]{1.06e-05} & ---   & 96.73 & \tablenum[table-format=1.2e3,exponent-product=\cdot]{1.11e-05} & ---   & \tablenum[table-format=3.2]{117.58} & \tablenum[table-format=1.2e3,exponent-product=\cdot]{3.29e-05} & ---   & 67.74 \\
				20  & \tablenum[table-format=1.2e3,exponent-product=\cdot]{2.57e-04} & 2.00                                 & \tablenum[table-format=1.2e3,exponent-product=\cdot]{4.17e-04} & 1.64                                   & \tablenum[table-format=1.2e3,exponent-product=\cdot]{5.61e-04} & 1.99                                   & \tablenum[table-format=1.2e3,exponent-product=\cdot]{2.74e-06} & 1.95  & 93.71 & \tablenum[table-format=1.2e3,exponent-product=\cdot]{4.33e-06} & 1.36  & \tablenum[table-format=3.2]{96.42}  & \tablenum[table-format=1.2e3,exponent-product=\cdot]{8.44e-06} & 1.96  & 66.43 \\
				40  & \tablenum[table-format=1.2e3,exponent-product=\cdot]{6.43e-05} & 2.00                                 & \tablenum[table-format=1.2e3,exponent-product=\cdot]{1.11e-04} & 1.91                                   & \tablenum[table-format=1.2e3,exponent-product=\cdot]{1.41e-04} & 2.00                                   & \tablenum[table-format=1.2e3,exponent-product=\cdot]{6.96e-07} & 1.98  & 92.37 & \tablenum[table-format=1.2e3,exponent-product=\cdot]{1.25e-06} & 1.79  & \tablenum[table-format=3.2]{88.87}  & \tablenum[table-format=1.2e3,exponent-product=\cdot]{2.14e-06} & 1.98  & 65.80 \\
				80  & \tablenum[table-format=1.2e3,exponent-product=\cdot]{1.61e-05} & 2.00                                 & \tablenum[table-format=1.2e3,exponent-product=\cdot]{2.80e-05} & 1.99                                   & \tablenum[table-format=1.2e3,exponent-product=\cdot]{3.51e-05} & 2.00                                   & \tablenum[table-format=1.2e3,exponent-product=\cdot]{1.74e-07} & 2.00  & 92.17 & \tablenum[table-format=1.2e3,exponent-product=\cdot]{3.22e-07} & 1.96  & \tablenum[table-format=3.2]{86.81}  & \tablenum[table-format=1.2e3,exponent-product=\cdot]{5.34e-07} & 2.00  & 65.70 \\
				160 & \tablenum[table-format=1.2e3,exponent-product=\cdot]{4.01e-06} & 2.00                                 & \tablenum[table-format=1.2e3,exponent-product=\cdot]{6.99e-06} & 2.00                                   & \tablenum[table-format=1.2e3,exponent-product=\cdot]{8.74e-06} & 2.01                                   & \tablenum[table-format=1.2e3,exponent-product=\cdot]{4.33e-08} & 2.01  & 92.65 & \tablenum[table-format=1.2e3,exponent-product=\cdot]{8.10e-08} & 1.99  & \tablenum[table-format=3.2]{86.38}  & \tablenum[table-format=1.2e3,exponent-product=\cdot]{1.32e-07} & 2.01  & 65.99 \\
				\bottomrule
			\end{tabular}
		}
	}
	\label{tab:euler_steady_polytropic_nG_2}

	\hspace{\fill}

	\vspace{12pt}

	\subfloat[Errors with a basis made of three elements: $q=2$.]{%
		\makebox[0.9\textwidth][c]{
			\small
			\begin{tabular}{c@{\hspace{10pt}}c@{\hspace{1pt}}c@{\hspace{10pt}}c@{\hspace{1pt}}c@{\hspace{10pt}}c@{\hspace{1pt}}c@{\hspace{10pt}}c@{\hspace{1pt}}c@{\hspace{1pt}}c@{\hspace{10pt}}c@{\hspace{1pt}}c@{\hspace{1pt}}c@{\hspace{10pt}}c@{\hspace{1pt}}c@{\hspace{1pt}}c}
				\toprule
				    & \multicolumn{2}{c}{$\rho$, basis $V_h$}                           & \multicolumn{2}{c}{$Q$, basis $V_h$} & \multicolumn{2}{c}{$E$, basis $V_h$}                           & \multicolumn{3}{c}{$\rho$, basis $V_h^+$} & \multicolumn{3}{c}{$Q$, basis $V_h^+$}                         & \multicolumn{3}{c}{$E$, basis $V_h^+$}                                                                                                                                                                                                                                                  \\
				\cmidrule(lr){2-3}\cmidrule(lr){4-5}\cmidrule(lr){6-7}\cmidrule(lr){8-10} \cmidrule(lr){11-13} \cmidrule(lr){14-16}
				$K$ & error                                                          & order                                & error                                                          & order                                  & error                                                          & order                                  & error                                                          & order & gain & error                                                          & order & gain & error                                                          & order & gain  \\
				\cmidrule(lr){1-1}\cmidrule(lr){2-3}\cmidrule(lr){4-5}\cmidrule(lr){6-7}\cmidrule(lr){8-10} \cmidrule(lr){11-13} \cmidrule(lr){14-16}
				10  & \tablenum[table-format=1.2e3,exponent-product=\cdot]{1.38e-06} & ---                                  & \tablenum[table-format=1.2e3,exponent-product=\cdot]{1.74e-06} & ---                                    & \tablenum[table-format=1.2e3,exponent-product=\cdot]{4.38e-05} & ---                                    & \tablenum[table-format=1.2e3,exponent-product=\cdot]{8.88e-07} & ---   & 1.55 & \tablenum[table-format=1.2e3,exponent-product=\cdot]{9.49e-07} & ---   & 1.83 & \tablenum[table-format=1.2e3,exponent-product=\cdot]{2.74e-06} & ---   & 15.97 \\
				20  & \tablenum[table-format=1.2e3,exponent-product=\cdot]{1.85e-07} & 2.90                                 & \tablenum[table-format=1.2e3,exponent-product=\cdot]{4.05e-07} & 2.10                                   & \tablenum[table-format=1.2e3,exponent-product=\cdot]{5.74e-06} & 2.93                                   & \tablenum[table-format=1.2e3,exponent-product=\cdot]{1.25e-07} & 2.83  & 1.48 & \tablenum[table-format=1.2e3,exponent-product=\cdot]{2.58e-07} & 1.88  & 1.57 & \tablenum[table-format=1.2e3,exponent-product=\cdot]{3.58e-07} & 2.94  & 16.01 \\
				40  & \tablenum[table-format=1.2e3,exponent-product=\cdot]{2.89e-08} & 2.68                                 & \tablenum[table-format=1.2e3,exponent-product=\cdot]{5.84e-08} & 2.79                                   & \tablenum[table-format=1.2e3,exponent-product=\cdot]{7.40e-07} & 2.95                                   & \tablenum[table-format=1.2e3,exponent-product=\cdot]{1.79e-08} & 2.81  & 1.62 & \tablenum[table-format=1.2e3,exponent-product=\cdot]{3.31e-08} & 2.96  & 1.76 & \tablenum[table-format=1.2e3,exponent-product=\cdot]{4.80e-08} & 2.90  & 15.43 \\
				80  & \tablenum[table-format=1.2e3,exponent-product=\cdot]{3.52e-09} & 3.04                                 & \tablenum[table-format=1.2e3,exponent-product=\cdot]{6.95e-09} & 3.07                                   & \tablenum[table-format=1.2e3,exponent-product=\cdot]{9.26e-08} & 3.00                                   & \tablenum[table-format=1.2e3,exponent-product=\cdot]{2.24e-09} & 2.99  & 1.57 & \tablenum[table-format=1.2e3,exponent-product=\cdot]{4.25e-09} & 2.96  & 1.63 & \tablenum[table-format=1.2e3,exponent-product=\cdot]{5.96e-09} & 3.01  & 15.55 \\
				160 & \tablenum[table-format=1.2e3,exponent-product=\cdot]{4.45e-10} & 2.98                                 & \tablenum[table-format=1.2e3,exponent-product=\cdot]{8.87e-10} & 2.97                                   & \tablenum[table-format=1.2e3,exponent-product=\cdot]{1.16e-08} & 3.00                                   & \tablenum[table-format=1.2e3,exponent-product=\cdot]{2.83e-10} & 2.98  & 1.57 & \tablenum[table-format=1.2e3,exponent-product=\cdot]{5.50e-10} & 2.95  & 1.61 & \tablenum[table-format=1.2e3,exponent-product=\cdot]{7.47e-10} & 2.99  & 15.50 \\
				\bottomrule
			\end{tabular}
		}
	}
	\label{tab:euler_steady_polytropic_nG_3}
	\caption{%
		Euler-Poisson system, polytropic pressure law:
		errors, orders of accuracy,
		and gain obtained when approximating a steady solution
		for bases with and without prior.
	}
	\label{tab:euler_steady_polytropic}
\end{table}

To extend this study, we compute the statistics over the whole parameter space
\eqref{eq:parameter_space_euler_polytropic} by uniformly sampling~$10^3$ values
and taking $10$ cells in the mesh.
The results are reported in \cref{tab:euler_steady_polytropic_stats}.
Just like before, the average gain is substantial,
while the minimum rarely falls below $1$.
Moreover, note that the gains recorded in \cref{tab:euler_steady_polytropic}
correspond to a rather bad set of parameters compared to the average.

\begin{table}[!ht]
	\centering
	\begin{tabular}{ccS[table-format=2.2]SS[table-format=2.2]cS[table-format=3.2]S[table-format=3.2]S[table-format=3.2]cS[table-format=4.2]S[table-format=4.2]S[table-format=4.2]}
		\toprule
		      &  & \multicolumn{3}{c}{minimum gain} &       & \multicolumn{3}{c}{average gain} &  & \multicolumn{3}{c}{maximum gain}                                                     \\
		\cmidrule(lr){3-5} \cmidrule(lr){7-9} \cmidrule(lr){11-13}
		{$q$} &  & {$\rho$}                         & {$Q$} & {$E$}                            &  & {$\rho$}                         & {$Q$}  & {$E$}  &  & {$\rho$} & {$Q$}   & {$E$}   \\
		\cmidrule(lr){1-1} \cmidrule(lr){3-5} \cmidrule(lr){7-9} \cmidrule(lr){11-13}
		0     &  & 19.14                            & 2.33  & 17.04                            &  & 233.48                           & 3.73   & 197.28 &  & 510.42   & 4.48    & 371.87  \\
		1     &  & 7.61                             & 8.28  & 6.98                             &  & 158.25                           & 188.92 & 130.57 &  & 1095.68  & 1291.90 & 1024.59 \\
		2     &  & 0.14                             & 0.22  & 2.99                             &  & 12.11                            & 16.55  & 23.73  &  & 89.47    & 109.93  & 169.28  \\
		\bottomrule
	\end{tabular}
	\caption{%
		Euler-Poisson system, polytropic pressure law:
        statistics of the gains obtained
        for the approximation of a steady solution
		in basis $V_h^+$ with respect to basis $V_h$.
	}
	\label{tab:euler_steady_polytropic_stats}
\end{table}

\subsubsection{Temperature-dependent pressure law}

In this case, we take a given smooth temperature
function $T(r, {\mu})$ parameterized by ${\mu}$,
where the parameter vector ${\mu}$
is composed of two elements:
\begin{equation*}
	{\mu} =
	\begin{pmatrix}
		\kappa \\
		\alpha \\
	\end{pmatrix}
	\in \mathbb{P} \subset \mathbb{R}^2,
	\quad \kappa \in \mathbb{R}_+,
	\quad \alpha \in \mathbb{R}_+.
\end{equation*}
This allows us to define the parameterized temperature function
$T(r; \alpha) = e^{- \alpha r}$,
and so we get the following temperature-based pressure law:
\begin{equation*}
	p(\rho; {\mu})
	=
	\kappa \rho T.
\end{equation*}
For this pressure law, the steady solutions are given by
the following nonlinear second-order ODE:
\begin{equation*}
	\frac{d}{dr} \left(
	r^2 \kappa \frac T \rho
	\frac{d \rho}{d r}
	\right)
	+ \frac{d}{dr} \left(
	r^2 \kappa
	\frac{d T}{d r}
	\right)
	= 4 \pi r^2 G \rho,
\end{equation*}
and the boundary condition $\partial_r p(0) = 0$ leads to
$\partial_r \rho(0) = \alpha$.
For this pressure law, we also take $\lambda = 1 + \sqrt{\gamma}$
in~\eqref{eq:time_step} to compute $\Delta t$.

The prior is obtained \emph{via} a PINN
with the same characteristics as in the polytropic case,
and whose result is still denoted by $\rho_\theta$.
To impose the boundary conditions, this time, we set
\begin{equation*}
	\widetilde{\rho}_\theta(r; {\mu}) =
	1 + \alpha r + r^2 \rho_\theta(r; {\mu}).
\end{equation*}
{\rb
Thanks to this expression, the conditions $\widetilde{\rho}_\theta(0; {\mu}) = 1$
and $\partial_r \widetilde{\rho}_\theta(0; {\mu}) = \alpha$
are automatically satisfied.
}
The parameter space is
\begin{equation}
	\label{eq:parameter_space_euler_temperature}
	\mathbb{P} = [2, 5] \times [0.5, 1.5],
\end{equation}
and the PINN is trained using only the physics-based loss function
\begin{equation*}
	\mathcal{J}(\theta) = \left\|
	\frac{d}{dr} \left(
	r^2 \kappa \frac T {\widetilde{\rho}_\theta}
	\frac{d \widetilde{\rho}_\theta}{d r}
	\right)
	+ \frac{d}{dr} \left(
	r^2 \kappa
	\frac{d T}{d r}
	\right)
	- 4 \pi r^2 G \widetilde{\rho}_\theta
	\right\|.
\end{equation*}
Training takes about 5 minutes on a dual NVIDIA K80 GPU,
until the loss function is equal to about $5 \cdot 10^{-4}$.
The priors~$Q_\theta$ and $E_\theta$ are then defined
in the same way as in the polytropic case.
In this case, we also take $n_q = q + 2$.

As is becoming usual, we first report, in \cref{tab:euler_steady_temperature},
the results of the approximation in both bases (with and without prior).
The final time is set to $T = 0.01$,
and we take $\kappa = 3.5$ and $\alpha = 0.5$.
As usual, using the prior provides significant gains,
especially for low values of $q$.
Compared to the polytropic case,
gains are consistently better for the large values of $q$.

\begin{table}[ht!]
	\centering
	\subfloat[Errors with a basis made of one element: $q=0$.]{%
		\makebox[0.9\textwidth][c]{
			\small
			\small
			\begin{tabular}{c@{\hspace{10pt}}c@{\hspace{1pt}}c@{\hspace{10pt}}c@{\hspace{1pt}}c@{\hspace{10pt}}c@{\hspace{1pt}}c@{\hspace{10pt}}c@{\hspace{1pt}}c@{\hspace{1pt}}c@{\hspace{10pt}}c@{\hspace{1pt}}c@{\hspace{1pt}}c@{\hspace{10pt}}c@{\hspace{1pt}}c@{\hspace{1pt}}c}
				\toprule
				    & \multicolumn{2}{c}{$\rho$, basis $V_h$} & \multicolumn{2}{c}{$Q$, basis $V_h$} & \multicolumn{2}{c}{$E$, basis $V_h$} & \multicolumn{3}{c}{$\rho$, basis $V_h^+$} & \multicolumn{3}{c}{$Q$, basis $V_h^+$} & \multicolumn{3}{c}{$E$, basis $V_h^+$}                                                                                                                                                     \\
				\cmidrule(lr){2-3}\cmidrule(lr){4-5}\cmidrule(lr){6-7}\cmidrule(lr){8-10} \cmidrule(lr){11-13} \cmidrule(lr){14-16}
				$K$ & error                                & order                                & error                                & order                                  & error                                  & order                                  & error                & order & gain   & error                & order & gain                               & error                & order & gain   \\
				\cmidrule(lr){1-1}\cmidrule(lr){2-3}\cmidrule(lr){4-5}\cmidrule(lr){6-7}\cmidrule(lr){8-10} \cmidrule(lr){11-13} \cmidrule(lr){14-16}
				10  & $3.91 \cdot 10^{-2}$                 & ---                                  & $2.80 \cdot 10^{-3}$                 & ---                                    & $1.56 \cdot 10^{-1}$                   & ---                                    & $2.65 \cdot 10^{-4}$ & ---   & 147.33 & $1.18 \cdot 10^{-3}$ & ---   & \tablenum[table-format=2.2]{2.37}  & $8.41 \cdot 10^{-4}$ & ---   & 186.02 \\
				20  & $1.96 \cdot 10^{-2}$                 & 1.00                                 & $1.66 \cdot 10^{-3}$                 & 0.76                                   & $7.83 \cdot 10^{-2}$                   & 1.00                                   & $1.39 \cdot 10^{-4}$ & 0.93  & 140.95 & $4.18 \cdot 10^{-4}$ & 1.50  & \tablenum[table-format=2.2]{3.97}  & $4.96 \cdot 10^{-4}$ & 0.76  & 157.73 \\
				40  & $9.81 \cdot 10^{-3}$                 & 1.00                                 & $9.02 \cdot 10^{-4}$                 & 0.88                                   & $3.92 \cdot 10^{-2}$                   & 1.00                                   & $7.04 \cdot 10^{-5}$ & 0.98  & 139.37 & $1.23 \cdot 10^{-4}$ & 1.77  & \tablenum[table-format=2.2]{7.35}  & $2.58 \cdot 10^{-4}$ & 0.94  & 151.77 \\
				80  & $4.91 \cdot 10^{-3}$                 & 1.00                                 & $5.30 \cdot 10^{-4}$                 & 0.77                                   & $1.96 \cdot 10^{-2}$                   & 1.00                                   & $3.61 \cdot 10^{-5}$ & 0.96  & 135.81 & $3.85 \cdot 10^{-5}$ & 1.67  & \tablenum[table-format=2.2]{13.75} & $1.41 \cdot 10^{-4}$ & 0.87  & 138.76 \\
				160 & $2.46 \cdot 10^{-3}$                 & 1.00                                 & $2.94 \cdot 10^{-4}$                 & 0.85                                   & $9.80 \cdot 10^{-3}$                   & 1.00                                   & $1.80 \cdot 10^{-5}$ & 1.00  & 136.36 & $1.09 \cdot 10^{-5}$ & 1.82  & \tablenum[table-format=2.2]{26.86} & $6.98 \cdot 10^{-5}$ & 1.02  & 140.52 \\
				\bottomrule
			\end{tabular}
		}
	}
	\label{tab:euler_steady_temperature_nG_1}

	\hspace{\fill}

	\vspace{12pt}

	\subfloat[Errors with a basis made of two elements: $q=1$.]{%
		\makebox[0.9\textwidth][c]{
			\small
			\begin{tabular}{c@{\hspace{10pt}}c@{\hspace{1pt}}c@{\hspace{10pt}}c@{\hspace{1pt}}c@{\hspace{10pt}}c@{\hspace{1pt}}c@{\hspace{10pt}}c@{\hspace{1pt}}c@{\hspace{1pt}}c@{\hspace{10pt}}c@{\hspace{1pt}}c@{\hspace{1pt}}c@{\hspace{10pt}}c@{\hspace{1pt}}c@{\hspace{1pt}}c}
				\toprule
				    & \multicolumn{2}{c}{$\rho$, basis $V_h$} & \multicolumn{2}{c}{$Q$, basis $V_h$} & \multicolumn{2}{c}{$E$, basis $V_h$} & \multicolumn{3}{c}{$\rho$, basis $V_h^+$} & \multicolumn{3}{c}{$Q$, basis $V_h^+$} & \multicolumn{3}{c}{$E$, basis $V_h^+$}                                                                                                                      \\
				\cmidrule(lr){2-3}\cmidrule(lr){4-5}\cmidrule(lr){6-7}\cmidrule(lr){8-10} \cmidrule(lr){11-13} \cmidrule(lr){14-16}
				$K$ & error                                & order                                & error                                & order                                  & error                                  & order                                  & error                & order & gain  & error                & order & gain  & error                & order & gain  \\
				\cmidrule(lr){1-1}\cmidrule(lr){2-3}\cmidrule(lr){4-5}\cmidrule(lr){6-7}\cmidrule(lr){8-10} \cmidrule(lr){11-13} \cmidrule(lr){14-16}
				10  & $1.89 \cdot 10^{-3}$                 & ---                                  & $2.19 \cdot 10^{-3}$                 & ---                                    & $4.61 \cdot 10^{-3}$                   & ---                                    & $2.85 \cdot 10^{-5}$ & ---   & 66.41 & $2.67 \cdot 10^{-5}$ & ---   & 81.99 & $8.80 \cdot 10^{-5}$ & ---   & 52.42 \\
				20  & $4.74 \cdot 10^{-4}$                 & 2.00                                 & $7.95 \cdot 10^{-4}$                 & 1.46                                   & $1.16 \cdot 10^{-3}$                   & 1.99                                   & $7.46 \cdot 10^{-6}$ & 1.93  & 63.53 & $1.16 \cdot 10^{-5}$ & 1.21  & 68.68 & $2.39 \cdot 10^{-5}$ & 1.88  & 48.51 \\
				40  & $1.19 \cdot 10^{-4}$                 & 2.00                                 & $2.31 \cdot 10^{-4}$                 & 1.79                                   & $2.90 \cdot 10^{-4}$                   & 2.00                                   & $1.92 \cdot 10^{-6}$ & 1.96  & 61.82 & $3.68 \cdot 10^{-6}$ & 1.65  & 62.63 & $6.28 \cdot 10^{-6}$ & 1.93  & 46.22 \\
				80  & $2.96 \cdot 10^{-5}$                 & 2.00                                 & $6.04 \cdot 10^{-5}$                 & 1.93                                   & $7.24 \cdot 10^{-5}$                   & 2.00                                   & $4.83 \cdot 10^{-7}$ & 1.99  & 61.27 & $1.00 \cdot 10^{-6}$ & 1.88  & 60.25 & $1.58 \cdot 10^{-6}$ & 2.00  & 45.96 \\
				160 & $7.40 \cdot 10^{-6}$                 & 2.00                                 & $1.51 \cdot 10^{-5}$                 & 2.00                                   & $1.80 \cdot 10^{-5}$                   & 2.01                                   & $1.20 \cdot 10^{-7}$ & 2.00  & 61.43 & $2.54 \cdot 10^{-7}$ & 1.98  & 59.61 & $3.80 \cdot 10^{-7}$ & 2.05  & 47.45 \\
				\bottomrule
			\end{tabular}
		}
	}
	\label{tab:euler_steady_temperature_nG_2}

	\hspace{\fill}

	\vspace{12pt}

	\subfloat[Errors with a basis made of three elements: $q=2$.]{%
		\makebox[0.9\textwidth][c]{
			\small
			\begin{tabular}{c@{\hspace{10pt}}c@{\hspace{1pt}}c@{\hspace{10pt}}c@{\hspace{1pt}}c@{\hspace{10pt}}c@{\hspace{1pt}}c@{\hspace{10pt}}c@{\hspace{1pt}}c@{\hspace{1pt}}c@{\hspace{10pt}}c@{\hspace{1pt}}c@{\hspace{1pt}}c@{\hspace{10pt}}c@{\hspace{1pt}}c@{\hspace{1pt}}c}
				\toprule
				    & \multicolumn{2}{c}{$\rho$, basis $V_h$} & \multicolumn{2}{c}{$Q$, basis $V_h$} & \multicolumn{2}{c}{$E$, basis $V_h$} & \multicolumn{3}{c}{$\rho$, basis $V_h^+$} & \multicolumn{3}{c}{$Q$, basis $V_h^+$} & \multicolumn{3}{c}{$E$, basis $V_h^+$}                                                                                                                                                                                \\
				\cmidrule(lr){2-3}\cmidrule(lr){4-5}\cmidrule(lr){6-7}\cmidrule(lr){8-10} \cmidrule(lr){11-13} \cmidrule(lr){14-16}
				$K$ & error                                & order                                & error                                & order                                  & error                                  & order                                  & error                & order & gain                               & error                & order & gain                               & error                & order & gain  \\
				\cmidrule(lr){1-1}\cmidrule(lr){2-3}\cmidrule(lr){4-5}\cmidrule(lr){6-7}\cmidrule(lr){8-10} \cmidrule(lr){11-13} \cmidrule(lr){14-16}
				10  & $3.83 \cdot 10^{-5}$                 & ---                                  & $4.49 \cdot 10^{-5}$                 & ---                                    & $1.27 \cdot 10^{-4}$                   & ---                                    & $2.77 \cdot 10^{-6}$ & ---   & \tablenum[table-format=2.2]{13.83} & $3.75 \cdot 10^{-6}$ & ---   & \tablenum[table-format=2.2]{11.98} & $8.95 \cdot 10^{-6}$ & ---   & 14.20 \\
				20  & $5.71 \cdot 10^{-6}$                 & 2.75                                 & $8.25 \cdot 10^{-6}$                 & 2.44                                   & $2.67 \cdot 10^{-5}$                   & 2.25                                   & $4.88 \cdot 10^{-7}$ & 2.50  & \tablenum[table-format=2.2]{11.70} & $7.62 \cdot 10^{-7}$ & 2.30  & \tablenum[table-format=2.2]{10.82} & $2.03 \cdot 10^{-6}$ & 2.14  & 13.14 \\
				40  & $7.37 \cdot 10^{-7}$                 & 2.95                                 & $8.72 \cdot 10^{-7}$                 & 3.24                                   & $3.66 \cdot 10^{-6}$                   & 2.87                                   & $7.19 \cdot 10^{-8}$ & 2.76  & \tablenum[table-format=2.2]{10.25} & $9.64 \cdot 10^{-8}$ & 2.98  & \tablenum[table-format=2.2]{9.05}  & $3.07 \cdot 10^{-7}$ & 2.73  & 11.93 \\
				80  & $8.88 \cdot 10^{-8}$                 & 3.05                                 & $1.09 \cdot 10^{-7}$                 & 3.00                                   & $4.48 \cdot 10^{-7}$                   & 3.03                                   & $8.89 \cdot 10^{-9}$ & 3.02  & \tablenum[table-format=2.2]{9.99}  & $1.14 \cdot 10^{-8}$ & 3.08  & \tablenum[table-format=2.2]{9.55}  & $3.85 \cdot 10^{-8}$ & 2.99  & 11.64 \\
				160 & $1.11 \cdot 10^{-8}$                 & 3.00                                 & $1.36 \cdot 10^{-8}$                 & 3.01                                   & $5.61 \cdot 10^{-8}$                   & 3.00                                   & $1.14 \cdot 10^{-9}$ & 2.96  & \tablenum[table-format=2.2]{9.74}  & $1.47 \cdot 10^{-9}$ & 2.96  & \tablenum[table-format=2.2]{9.23}  & $4.96 \cdot 10^{-9}$ & 2.96  & 11.31 \\
				\bottomrule
			\end{tabular}
		}
	}
	\label{tab:euler_steady_temperature_nG_3}
	\caption{%
		Euler-Poisson system, temperature-based pressure law:
		errors, orders of accuracy,
		and gain obtained when approximating a steady solution
		for bases with and without prior.
	}
	\label{tab:euler_steady_temperature}
\end{table}

To understand gains on the whole parameter space
\eqref{eq:parameter_space_euler_temperature},
we uniformly sample $10^3$ values of $\kappa$ and $\alpha$ and take a mesh made of $10$ cells.
We compute the minimum, average and maximum gains.
These values are reported in \cref{tab:euler_steady_temperature_stats}.
For this pressure law, the minimum gain is always larger than $1$,
and we obtain consistently large average gains, even for $q = 2$.

\begin{table}[!ht]
	\centering
	\begin{tabular}{ccS[table-format=2.2]SS[table-format=2.2]cS[table-format=3.2]S[table-format=2.2]S[table-format=3.2]cS[table-format=3.2]S[table-format=3.2]S[table-format=3.2]}
		\toprule
		      &  & \multicolumn{3}{c}{minimum gain} &       & \multicolumn{3}{c}{average gain} &  & \multicolumn{3}{c}{maximum gain}                                                  \\
		\cmidrule(lr){3-5} \cmidrule(lr){7-9} \cmidrule(lr){11-13}
		{$q$} &  & {$\rho$}                         & {$Q$} & {$E$}                            &  & {$\rho$}                         & {$Q$} & {$E$}  &  & {$\rho$} & {$Q$}  & {$E$}  \\
		\cmidrule(lr){1-1} \cmidrule(lr){3-5} \cmidrule(lr){7-9} \cmidrule(lr){11-13}
		0     &  & 13.30                            & 1.05  & 16.24                            &  & 151.96                           & 1.88  & 150.63 &  & 600.13   & 2.91   & 473.83 \\
		1     &  & 6.30                             & 7.53  & 5.40                             &  & 72.63                            & 77.20 & 51.09  &  & 321.20   & 302.58 & 257.19 \\
		2     &  & 3.35                             & 3.45  & 2.20                             &  & 18.96                            & 22.58 & 13.56  &  & 55.47    & 63.45  & 47.83  \\
		\bottomrule
	\end{tabular}
	\caption{%
		Euler-Poisson system, temperature-based pressure law:
        statistics of the gains obtained
        for the approximation of a steady solution
		in basis $V_h^+$ with respect to basis $V_h$.
	}
	\label{tab:euler_steady_temperature_stats}
\end{table}

\subsubsection{Spherical blast wave}

The goal of this last test case is to show that our prior
does not negatively affect the capability of the scheme
to capture discontinuous solutions.
Let us emphasize that numerical viscosity is not
an object of this study,
and therefore that we have not used any regularization procedure.
Consequently, results will show some oscillations.

This experiment is nothing but a Riemann problem in spherical geometry,
inspired by the experiments in \cite{Tor2009}.
As such, the initial condition is piecewise constant
on the space domain $r \in (0, 0.4)$, as follows:
\begin{equation*}
	\rho(0, x) =
	\begin{dcases}
		2 & \text{if } r < 0.2, \\
		1 & \text{otherwise;}
	\end{dcases}
	\qquad
	Q(0, x) = 0;
	\qquad
	p(0, x) =
	\begin{dcases}
		2 & \text{if } r < 0.2, \\
		1 & \text{otherwise.}
	\end{dcases}
\end{equation*}
For this experiment, the pressure law is the standard ideal gas law
\begin{equation*}
	p = (\gamma - 1) \left( E - \frac 1 2 \frac {Q^2} \rho \right),
\end{equation*}
and we take the gas constant $\gamma$ equal to $1.4$.
The experiment is run until the final time $T = 0.1$,
and with Neumann boundary conditions.
We take $25$ discretization cells,
and we use a basis made of $3$ elements.
Moreover, the source term is deactivated:
we set $\phi = 0$, and we merely consider the Euler equations
in spherical geometry, without gravity effects.

\begin{figure}[!ht]
	\centering
	\includegraphics{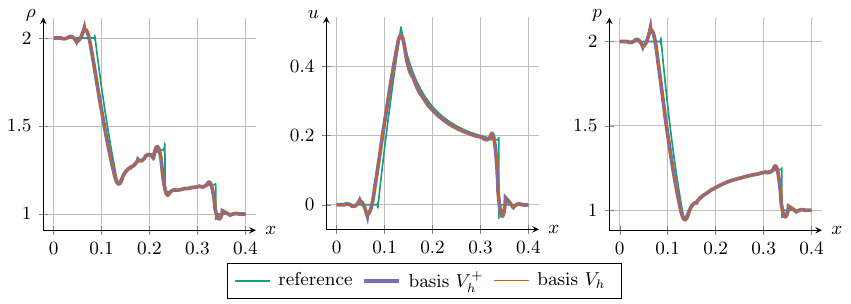}
	\caption{%
        \rb%
		Euler equations, ideal gas pressure law:
		results for the spherical blast wave.
        From left to right, we plot the density,
        the velocity $u = Q / \rho$, and the pressure
        with respect to space.
	}
	\label{fig:euler_spherical_blast_wave}
\end{figure}

The results are depicted in \cref{fig:euler_spherical_blast_wave},
where we compare the two bases
(with and without prior, blue and orange lines respectively)
to a reference solution (green line).
We observe very good agreement with the reference solution,
even though oscillations are present, as expected.
We also note that the graphs for the solutions
with and without prior are superimposed,
which means that the quality of the approximation of this
discontinuous solution has not been degraded
by the introduction of the prior in the basis.

\subsection{Shallow water equations in two space dimensions}
\label{sec:shallow_water_2d}

The last system considered in this series of
experiments is the two-dimensional
shallow water system.
It is given by
\begin{equation}
    \label{eq:shallow_water_2d}
	\begin{dcases}
		\partial_t h + \boldsymbol \nabla \cdot \boldsymbol Q = 0, \\
		\partial_t \boldsymbol Q + \boldsymbol \nabla \cdot \left( \frac{\boldsymbol Q \otimes \boldsymbol Q} h + \frac 1 2 g h^2 \text{Id} \right) = -g h \boldsymbol \nabla Z(\boldsymbol x; {\mu}),
	\end{dcases}
\end{equation}
where $g$ is the gravity constant,
$\text{Id}$ is the $2 \times 2$ identity matrix,
$h$ is the water height,
$\boldsymbol Q$ is the water discharge,
and~$Z$ is the topography.
For this system, $\Delta t$ is computed
by setting $\lambda = 2 + \sqrt{\gamma}$ in \eqref{eq:time_step}.
{\rb The boundary conditions are described later;
we will prescribe Dirichlet boundary conditions
according to the value of the steady solution at the boundaries.}

The space variable $\boldsymbol x = (x_1, x_2)$
belongs to the space domain $\Omega = [-3, 3]^2$,
and we introduce three parameters:
\begin{equation*}
	{\mu} =
	\begin{pmatrix}
		\alpha \\
		\Gamma \\
		r_0    \\
	\end{pmatrix}
	\in \mathbb{P} \subset \mathbb{R}^3,
	\quad \alpha \in \mathbb{R}^*_+,
	\quad \Gamma \in \mathbb{R}^*_+,
	\quad r_0 \in \mathbb{R}^*_+.
\end{equation*}
{\rb%
In practice, the parameter space is
\begin{equation*}
	\mathbb{P} = [ 0.25, 0.75 ] \times [ 0.1, 0.4 ] \times [ 0.5, 1.25 ].
\end{equation*}
}
This enables us to define the topography as
the following Gaussian bump function,
with $r = \|\boldsymbol{x}\|$:
\begin{equation*}
	Z(\boldsymbol x; {\mu}) =
	\Gamma \exp \left( \alpha (r_0^2 - r^2) \right),
\end{equation*}
see for instance \cite{Shu1998} for a similar test case.
{\rb %
Since this topography is radially symmetric,
we can expect that radially symmetric steady solutions also exist.
To derive such a steady solution, let us assume that
\begin{equation}
    \label{eq:shallow_water_2d_steady_Q}
    \boldsymbol{Q}_\text{eq}(\boldsymbol{x}; {\mu})
    =
    - \boldsymbol{x}^\perp \, h_\text{eq}(\boldsymbol{x}; {\mu}) \,
    u_\text{eq}(\boldsymbol{x}; {\mu}),
\end{equation}
with $\boldsymbol{x}^\perp = (-x_2, x_1)$
and where $h_\text{eq}$ and $u_\text{eq}$ are radial functions, to be determined.
One can easily check that such an expression $\boldsymbol{Q}_\text{eq}$ is divergence-free,
which is the first requirement for a steady solution.
The second requirement is given by the second equation of \eqref{eq:shallow_water_2d},
which, in our case, reduces to
\begin{equation*}
    \boldsymbol \nabla \cdot \Big(
        \left({\boldsymbol{x}^\perp \otimes \boldsymbol{x}^\perp}\right)
        h_\text{eq}(\boldsymbol{x}; {\mu}) \,
        u_\text{eq}(\boldsymbol{x}; {\mu})^2
    \Big)
    + g \boldsymbol \nabla \Big(
        h_\text{eq}(\boldsymbol x; {\mu}) + Z(\boldsymbol x; {\mu})
    \Big)
    = 0.
\end{equation*}
Arguing radial symmetry, straightforward computations lead to the following system:
\begin{equation}
    \label{eq:shallow_water_2d_steady_system}
    \begin{aligned}
        - x_1 u_\text{eq}(\boldsymbol{x}; {\mu})^2
        + g \partial_{x_1} \Big(
            h_\text{eq}(\boldsymbol x; {\mu}) + Z(\boldsymbol x; {\mu})
        \Big) &= 0, \\
        - x_2 u_\text{eq}(\boldsymbol{x}; {\mu})^2
        + g \partial_{x_2} \Big(
            h_\text{eq}(\boldsymbol x; {\mu}) + Z(\boldsymbol x; {\mu})
        \Big) &= 0.
    \end{aligned}
\end{equation}
A solution to this system is given by
\begin{equation}
	\label{eq:shallow_water_2d_steady_h_u}
	\begin{dcases}
		h_\text{eq}(\boldsymbol{x}; {\mu})
		=
		2
		- Z(\boldsymbol{x}; {\mu})
		- \frac{\alpha \Gamma}{8 g} Z(\boldsymbol{x}; {\mu})^4, \\
        u_\text{eq}(\boldsymbol{x}; {\mu})
        =
        \alpha Z(\boldsymbol x; {\mu})^2,
	\end{dcases}
\end{equation}
and the full steady solution is thus governed
by \eqref{eq:shallow_water_2d_steady_Q} -- \eqref{eq:shallow_water_2d_steady_h_u}.%
}

To obtain a relevant prior,
we approximate $h_\text{eq}$ and $u_\text{eq}$,
using a different PINN for each of the two functions.
{\rb %
Another possibility would be to strongly impose the divergence-free
constraint, by learning a potential and taking the prior
$\boldsymbol{Q}_\theta(\boldsymbol{x}; {\mu})$
as the curl of this potential.
However, we elected not to do so,
since the current strategy resulted in faster training.
The results of the PINN are then denoted by
$h_\theta$ and $u_\theta$,
and we define the priors
$\smash{\widetilde{h}_\theta}$ and $\smash{\widetilde{u}_\theta}$ as follows,
to include the boundary conditions:
\begin{equation*}
	\widetilde{h}_\theta(\boldsymbol x; {\mu})
	=
	2 - Z(\boldsymbol{x}; {\mu}) \,
	h_\theta(\boldsymbol x; {\mu})^2
	\text{\quad and \quad}
	\widetilde{u}_\theta(\boldsymbol x; {\mu})
	=
	Z(\boldsymbol{x}; {\mu}) \,
	u_\theta(\boldsymbol x; {\mu}).
\end{equation*}
Note that $Z$ almost vanishes at the domain boundaries,
as does $u_\text{eq}$.
The water height $h_\text{eq}$, for its part, is almost constant at the boundaries.
These boundary conditions are thus correctly imposed
on the priors $\smash{\widetilde{h}_\theta}$ and $\smash{\widetilde{u}_\theta}$
by multiplying the PINN outputs $h_\theta$ and $u_\theta$ by $Z$.%
}

The loss function is made in equal parts of the now usual PDE loss function,
and of the minimization with respect to data.
{\rb%
Therefore, as described in \cref{sec:parametric_PINNs},
the loss function is given by
$\mathcal{J}(\theta) = \mathcal{J}_r(\theta) + \mathcal{J}_\text{data}(\theta)$,
where the residual loss function is defined as
\begin{equation*}
    \mathcal{J}_r(\theta) =
    \biggl\|
        - x_1 \widetilde{u}_\theta^2
        + g \partial_{x_1} \Big(
            \widetilde{h}_\theta + Z
        \Big)
    \biggr\|_2^2
    +
    \biggl\|
        - x_2 \widetilde{u}_\theta^2
        + g \partial_{x_2} \Big(
            \widetilde{h}_\theta + Z
        \Big)
    \biggr\|_2^2,
\end{equation*}
and where the data loss function is defined as
\begin{equation*}
    \mathcal{J}_\text{data}(\theta) =
    \biggl\|
        \widetilde{h}_\theta(\boldsymbol{x}; {\mu})
        -
        {h}_\text{eq}(\boldsymbol{x}; {\mu})
    \biggr\|_2^2
    +
    \biggl\|
        \widetilde{u}_\theta(\boldsymbol{x}; {\mu})
        -
        {u}_\text{eq}(\boldsymbol{x}; {\mu})
    \biggr\|_2^2.
\end{equation*}
Thanks to $\mathcal{J}_r(\theta)$, the PINN ensures that
$(\widetilde{h}_\theta, \widetilde{u}_\theta)$
approximately satisfies the steady PDE \eqref{eq:shallow_water_2d_steady_system}.
Data is regenerated at each epoch,
and helps to avoid falling in a local minimum
corresponding to a lake at rest, where
$\widetilde{h}_\theta + Z = \text{constant}$ and $u_\theta = 0$.%
}
Each PINN has about $2500$ parameters,
and training takes about 10 minutes on an NVIDIA V100 GPU,
until the loss function reaches about $4 \times 10^{-7}$.
This prior is integrated with
a quadrature formula of degree $n_q = q + 3$:
we needed to increase the usual quadrature degree by $1$
to obtain the best possible approximation.

\subsubsection{Approximation of a steady solution}

We take the steady solution
\eqref{eq:shallow_water_2d_steady_Q} -- \eqref{eq:shallow_water_2d_steady_h_u}
as the initial condition to test the approximate well-balanced property.
The experiments are run until the final physical time $T = 0.01$.
We prescribe Dirichlet boundary conditions consisting
in the value of the steady solution.

First, we take the parameters as the center of the parameter cube $\mathbb{P}$.
The results are collected in \cref{tab:shallow_water_2d_steady},
and we note that, as expected, the presence of the prior
makes it possible to reach much lower errors,
especially for the water height $h$.

\begin{table}[ht!]
	\centering
	\subfloat[Errors with a basis made of one element: $q=0$.]{%
		\makebox[0.9\textwidth][c]{
			\small
			\begin{tabular}{c@{\hspace{10pt}}c@{\hspace{1pt}}c@{\hspace{10pt}}c@{\hspace{1pt}}c@{\hspace{10pt}}c@{\hspace{1pt}}c@{\hspace{10pt}}c@{\hspace{1pt}}c@{\hspace{1pt}}c@{\hspace{10pt}}c@{\hspace{1pt}}c@{\hspace{1pt}}c@{\hspace{10pt}}c@{\hspace{1pt}}c@{\hspace{1pt}}c}
				\toprule
				    & \multicolumn{2}{c}{$h$, basis $V_h$} & \multicolumn{2}{c}{$Q_1$, basis $V_h$} & \multicolumn{2}{c}{$Q_2$, basis $V_h$} & \multicolumn{3}{c}{$h$, basis $V_h^+$} & \multicolumn{3}{c}{$Q_1$, basis $V_h^+$} & \multicolumn{3}{c}{$Q_2$, basis $V_h^+$}                                                                                                                      \\
				\cmidrule(lr){2-3}\cmidrule(lr){4-5}\cmidrule(lr){6-7}\cmidrule(lr){8-10} \cmidrule(lr){11-13} \cmidrule(lr){14-16}
				$K$ & error                                & order                                  & error                                  & order                                  & error                                    & order                                    & error                & order & gain    & error                & order & gain & error                & order & gain \\
				\cmidrule(lr){1-1}\cmidrule(lr){2-3}\cmidrule(lr){4-5}\cmidrule(lr){6-7}\cmidrule(lr){8-10} \cmidrule(lr){11-13} \cmidrule(lr){14-16}
				20  & $1.94 \cdot 10^{-1}$                 & ---                                    & $4.31 \cdot 10^{-2}$                   & ---                                    & $4.31 \cdot 10^{-2}$                     & ---                                      & $1.31 \cdot 10^{-4}$ & ---   & 1477.51 & $6.24 \cdot 10^{-3}$ & ---   & 6.91 & $6.24 \cdot 10^{-3}$ & ---   & 6.91 \\
				40  & $9.75 \cdot 10^{-2}$                 & 0.99                                   & $2.19 \cdot 10^{-2}$                   & 0.98                                   & $2.19 \cdot 10^{-2}$                     & 0.98                                     & $6.37 \cdot 10^{-5}$ & 1.04  & 1531.52 & $2.84 \cdot 10^{-3}$ & 1.14  & 7.69 & $2.84 \cdot 10^{-3}$ & 1.14  & 7.69 \\
				80  & $4.88 \cdot 10^{-2}$                 & 1.00                                   & $1.09 \cdot 10^{-2}$                   & 1.00                                   & $1.09 \cdot 10^{-2}$                     & 1.00                                     & $3.17 \cdot 10^{-5}$ & 1.01  & 1540.17 & $1.43 \cdot 10^{-3}$ & 0.99  & 7.63 & $1.43 \cdot 10^{-3}$ & 0.99  & 7.63 \\
				160 & $2.44 \cdot 10^{-2}$                 & 1.00                                   & $5.48 \cdot 10^{-3}$                   & 1.00                                   & $5.48 \cdot 10^{-3}$                     & 1.00                                     & $1.59 \cdot 10^{-5}$ & 1.00  & 1539.94 & $7.21 \cdot 10^{-4}$ & 0.99  & 7.60 & $7.21 \cdot 10^{-4}$ & 0.99  & 7.60 \\
				320 & $1.22 \cdot 10^{-2}$                 & 1.00                                   & $2.74 \cdot 10^{-3}$                   & 1.00                                   & $2.74 \cdot 10^{-3}$                     & 1.00                                     & $7.93 \cdot 10^{-6}$ & 1.00  & 1539.59 & $3.61 \cdot 10^{-4}$ & 1.00  & 7.58 & $3.61 \cdot 10^{-4}$ & 1.00  & 7.58 \\
				\bottomrule
			\end{tabular}
		}
	}
	\label{tab:shallow_water_2d_steady_nG_1}

	\hspace{\fill}

	\vspace{12pt}

	\subfloat[Errors with a basis made of two elements: $q=1$.]{%
		\makebox[0.9\textwidth][c]{
			\small
			\begin{tabular}{c@{\hspace{10pt}}c@{\hspace{1pt}}c@{\hspace{10pt}}c@{\hspace{1pt}}c@{\hspace{10pt}}c@{\hspace{1pt}}c@{\hspace{10pt}}c@{\hspace{1pt}}c@{\hspace{1pt}}c@{\hspace{10pt}}c@{\hspace{1pt}}c@{\hspace{1pt}}c@{\hspace{10pt}}c@{\hspace{1pt}}c@{\hspace{1pt}}c}
				\toprule
				    & \multicolumn{2}{c}{$h$, basis $V_h$} & \multicolumn{2}{c}{$Q_1$, basis $V_h$} & \multicolumn{2}{c}{$Q_2$, basis $V_h$} & \multicolumn{3}{c}{$h$, basis $V_h^+$} & \multicolumn{3}{c}{$Q_1$, basis $V_h^+$} & \multicolumn{3}{c}{$Q_2$, basis $V_h^+$}                                                                                                                       \\
				\cmidrule(lr){2-3}\cmidrule(lr){4-5}\cmidrule(lr){6-7}\cmidrule(lr){8-10} \cmidrule(lr){11-13} \cmidrule(lr){14-16}
				$K$ & error                                & order                                  & error                                  & order                                  & error                                    & order                                    & error                & order & gain   & error                & order & gain  & error                & order & gain  \\
				\cmidrule(lr){1-1}\cmidrule(lr){2-3}\cmidrule(lr){4-5}\cmidrule(lr){6-7}\cmidrule(lr){8-10} \cmidrule(lr){11-13} \cmidrule(lr){14-16}
				20  & $2.17 \cdot 10^{-2}$                 & ---                                    & $2.58 \cdot 10^{-2}$                   & ---                                    & $2.58 \cdot 10^{-2}$                     & ---                                      & $8.51 \cdot 10^{-5}$ & ---   & 254.60 & $1.42 \cdot 10^{-3}$ & ---   & 18.21 & $1.42 \cdot 10^{-3}$ & ---   & 18.21 \\
				40  & $5.46 \cdot 10^{-3}$                 & 1.99                                   & $8.88 \cdot 10^{-3}$                   & 1.54                                   & $8.88 \cdot 10^{-3}$                     & 1.54                                     & $3.23 \cdot 10^{-5}$ & 1.40  & 169.11 & $3.70 \cdot 10^{-4}$ & 1.94  & 23.99 & $3.70 \cdot 10^{-4}$ & 1.94  & 23.99 \\
				80  & $1.37 \cdot 10^{-3}$                 & 2.00                                   & $2.50 \cdot 10^{-3}$                   & 1.83                                   & $2.50 \cdot 10^{-3}$                     & 1.83                                     & $9.43 \cdot 10^{-6}$ & 1.78  & 145.10 & $9.35 \cdot 10^{-5}$ & 1.98  & 26.74 & $9.35 \cdot 10^{-5}$ & 1.98  & 26.74 \\
				160 & $3.42 \cdot 10^{-4}$                 & 2.00                                   & $6.46 \cdot 10^{-4}$                   & 1.95                                   & $6.46 \cdot 10^{-4}$                     & 1.95                                     & $2.47 \cdot 10^{-6}$ & 1.94  & 138.89 & $2.35 \cdot 10^{-5}$ & 2.00  & 27.54 & $2.35 \cdot 10^{-5}$ & 2.00  & 27.54 \\
				320 & $8.56 \cdot 10^{-5}$                 & 2.00                                   & $1.62 \cdot 10^{-4}$                   & 2.00                                   & $1.62 \cdot 10^{-4}$                     & 2.00                                     & $6.19 \cdot 10^{-7}$ & 1.99  & 138.25 & $5.87 \cdot 10^{-6}$ & 2.00  & 27.55 & $5.87 \cdot 10^{-6}$ & 2.00  & 27.55 \\
				\bottomrule
			\end{tabular}
		}
	}
	\label{tab:shallow_water_2d_steady_nG_2}

	\hspace{\fill}

	\vspace{12pt}

	\subfloat[Errors with a basis made of three elements: $q=2$.]{%
		\makebox[0.9\textwidth][c]{
			\small
			\begin{tabular}{c@{\hspace{10pt}}c@{\hspace{1pt}}c@{\hspace{10pt}}c@{\hspace{1pt}}c@{\hspace{10pt}}c@{\hspace{1pt}}c@{\hspace{10pt}}c@{\hspace{1pt}}c@{\hspace{1pt}}c@{\hspace{10pt}}c@{\hspace{1pt}}c@{\hspace{1pt}}c@{\hspace{10pt}}c@{\hspace{1pt}}c@{\hspace{1pt}}c}
				\toprule
				    & \multicolumn{2}{c}{$h$, basis $V_h$} & \multicolumn{2}{c}{$Q_1$, basis $V_h$} & \multicolumn{2}{c}{$Q_2$, basis $V_h$} & \multicolumn{3}{c}{$h$, basis $V_h^+$} & \multicolumn{3}{c}{$Q_1$, basis $V_h^+$} & \multicolumn{3}{c}{$Q_2$, basis $V_h^+$}                                                                                                                      \\
				\cmidrule(lr){2-3}\cmidrule(lr){4-5}\cmidrule(lr){6-7}\cmidrule(lr){8-10} \cmidrule(lr){11-13} \cmidrule(lr){14-16}
				$K$ & error                                & order                                  & error                                  & order                                  & error                                    & order                                    & error                & order & gain  & error                & order & gain  & error                & order & gain  \\
				\cmidrule(lr){1-1}\cmidrule(lr){2-3}\cmidrule(lr){4-5}\cmidrule(lr){6-7}\cmidrule(lr){8-10} \cmidrule(lr){11-13} \cmidrule(lr){14-16}
				20  & $1.61 \cdot 10^{-3}$                 & ---                                    & $3.03 \cdot 10^{-3}$                   & ---                                    & $3.03 \cdot 10^{-3}$                     & ---                                      & $1.63 \cdot 10^{-5}$ & ---   & 98.79 & $2.95 \cdot 10^{-4}$ & ---   & 10.27 & $2.95 \cdot 10^{-4}$ & ---   & 10.27 \\
				40  & $2.18 \cdot 10^{-4}$                 & 2.89                                   & $4.83 \cdot 10^{-4}$                   & 2.65                                   & $4.83 \cdot 10^{-4}$                     & 2.65                                     & $2.55 \cdot 10^{-6}$ & 2.68  & 85.60 & $4.03 \cdot 10^{-5}$ & 2.87  & 11.97 & $4.03 \cdot 10^{-5}$ & 2.87  & 11.97 \\
				80  & $2.85 \cdot 10^{-5}$                 & 2.94                                   & $5.77 \cdot 10^{-5}$                   & 3.06                                   & $5.77 \cdot 10^{-5}$                     & 3.06                                     & $3.12 \cdot 10^{-7}$ & 3.03  & 91.29 & $5.11 \cdot 10^{-6}$ & 2.98  & 11.30 & $5.11 \cdot 10^{-6}$ & 2.98  & 11.30 \\
				160 & $3.47 \cdot 10^{-6}$                 & 3.04                                   & $6.86 \cdot 10^{-6}$                   & 3.07                                   & $6.86 \cdot 10^{-6}$                     & 3.07                                     & $3.69 \cdot 10^{-8}$ & 3.08  & 94.23 & $6.33 \cdot 10^{-7}$ & 3.01  & 10.84 & $6.33 \cdot 10^{-7}$ & 3.01  & 10.84 \\
				320 & $4.35 \cdot 10^{-7}$                 & 3.00                                   & $8.56 \cdot 10^{-7}$                   & 3.00                                   & $8.56 \cdot 10^{-7}$                     & 3.00                                     & $4.66 \cdot 10^{-9}$ & 2.98  & 93.43 & $7.85 \cdot 10^{-8}$ & 3.01  & 10.91 & $7.85 \cdot 10^{-8}$ & 3.01  & 10.91 \\
				\bottomrule
			\end{tabular}
		}
	}
	\label{tab:shallow_water_2d_steady_nG_3}
	\caption{%
		Shallow water equations in two space dimensions:
		errors, orders of accuracy,
		and gain obtained when approximating a steady solution
		for bases with and without prior.
	}
	\label{tab:shallow_water_2d_steady}
\end{table}

In addition, we provide some statistics over the whole parameter space $\mathbb{P}$,
computed on a mesh with $25 \times 25$ cells, in
\cref{tab:shallow_water_2d_steady_stats}.
{\ra We note that}, on average,
the gains are substantial.
However, note that the minimum gains may be
smaller than $1$, which denotes a loss of precision
due to the prior.
This happens in around 0.75\% of cases,
so we obtain an improvement in an overwhelming
majority of cases.

\begin{table}[!ht]
	\centering
	\begin{tabular}{ccS[table-format=2.2]SScS[table-format=4.2]S[table-format=2.2]S[table-format=2.2]cS[table-format=4.2]SS}
		\toprule
		      &  & \multicolumn{3}{c}{minimum gain} &         & \multicolumn{3}{c}{average gain} &  & \multicolumn{3}{c}{maximum gain}                                                       \\
		\cmidrule(lr){3-5} \cmidrule(lr){7-9} \cmidrule(lr){11-13}
		{$q$} &  & {$h$}                         & {$Q_1$} & {$Q_2$}                          &  & {$h$}                         & {$Q_1$} & {$Q_2$} &  & {$h$} & {$Q_1$} & {$Q_2$} \\
		\cmidrule(lr){1-1} \cmidrule(lr){3-5} \cmidrule(lr){7-9} \cmidrule(lr){11-13}
		0     &  & 48.21                            & 0.39    & 0.39                             &  & 1131.55                          & 6.18    & 6.18    &  & 1592.52  & 11.24   & 11.24   \\
		1     &  & 0.96                             & 0.16    & 0.16                             &  & 186.82                           & 21.00   & 21.00   &  & 422.68   & 49.05   & 49.05   \\
		2     &  & 0.06                             & 0.02    & 0.02                             &  & 82.45                            & 8.76    & 8.76    &  & 206.29   & 22.43   & 22.43   \\
		\bottomrule
	\end{tabular}
	\caption{%
		Shallow water equations in two space dimensions:
		statistics of the gains obtained
		for the approximation of a steady solution
		in basis $V_h^+$ with respect to basis $V_h$.
	}
	\label{tab:shallow_water_2d_steady_stats}
\end{table}

\subsubsection{Perturbed steady solution}

We now compare the new basis to the classical one
when the initial condition is a perturbed steady solution.
To that end, the initial water height is set to
\begin{equation*}
	h(0, \boldsymbol{x}; {\mu})
	=
	h_\text{eq}(\boldsymbol{x}; {\mu})
	- 0.02 \exp \left( - 2 ((x_1 + 2)^2 + (x_2 + 2)^2) \right),
\end{equation*}
thus creating a bump-shaped perturbation whose center is located at $(-2, -2)$.
For simplicity, we still use the value of the steady solution
as Dirichlet boundary conditions.
Moreover, we set the parameters $\mu$ to be the center of the cube $\mathbb{P}$,
and we take $q = 1$ with $16^2$ discretization cells.

\begin{figure}[!ht]
	\centering
	\includegraphics[width=\textwidth]{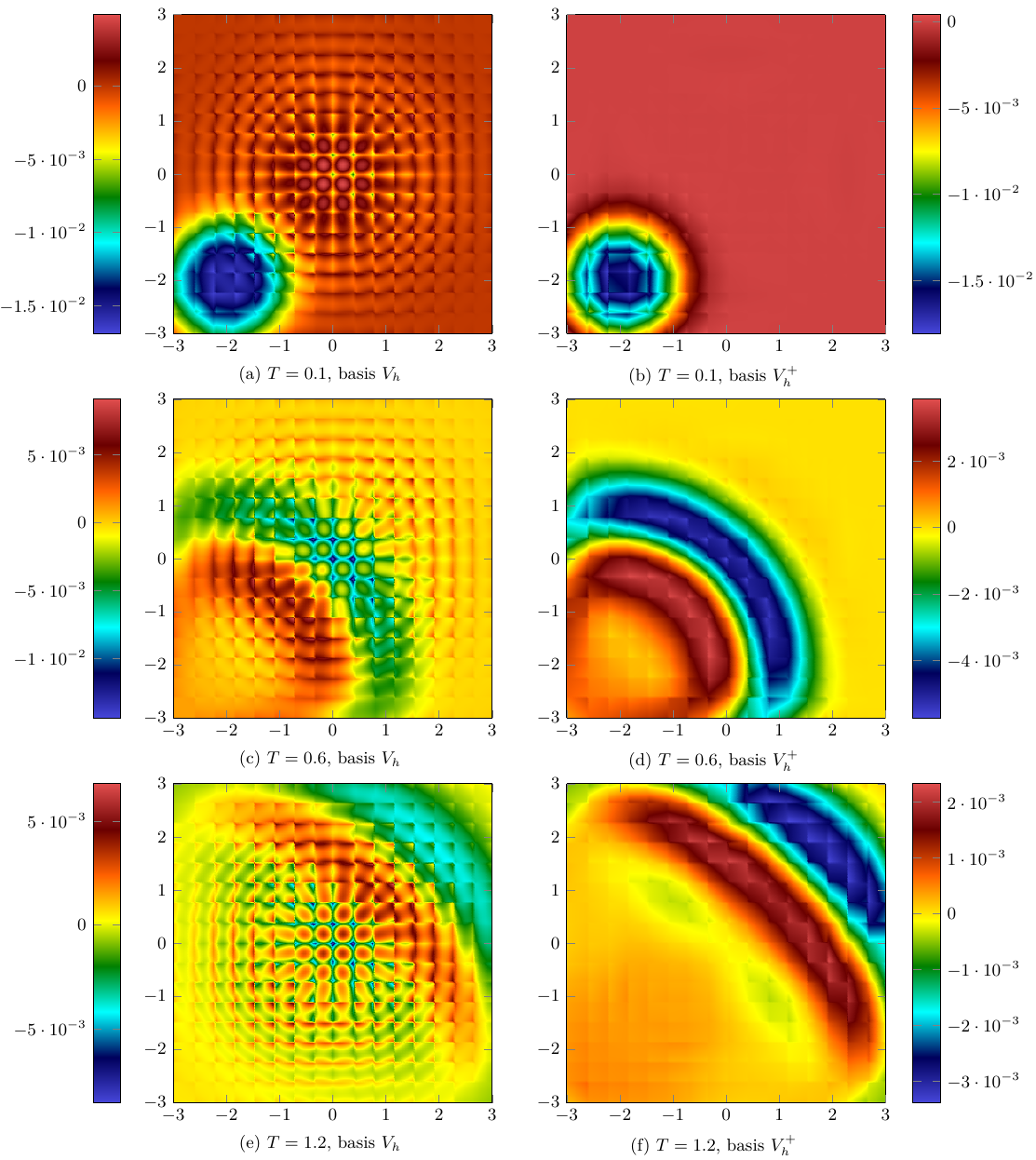}
	\caption{%
		2D shallow water equations:
		representation of the error between the perturbed solution
		and the underlying steady solution.
		Left panels: classical basis $V_h$;
		right panels: prior-enriched basis $V_h^+$.
		From top to bottom, several final times are displayed
		($T = 0.1$, $T = 0.6$, $T = 1.2$).
	}
	\label{fig:2D_SW_perturbed}
\end{figure}

The pointwise difference between $h$ and $h_\text{eq}$
is displayed in \cref{fig:2D_SW_perturbed}.
We observe that the prior-enriched basis~$V_h^+$ (right panels)
is able to capture the perturbation much better
than the classical basis $V_h$ (left panels).
Indeed, the background, underlying steady solution has been
smeared by basis $V_h$,
while is preserved with much greater resolution by basis $V_h^+$.

\section{Conclusion}
\label{sec:conclusion}

In this work, we proposed a Discontinuous Galerkin scheme
whose basis has been enriched by neural networks
to ensure an approximate well-balance property
for a generic PDE and a generic equilibrium.
The offline phase of the algorithm consists in learning
a family of equilibria using parametric PINNs.
Then, during the online phase,
the trained network is used to enrich the DG basis
and to approximate the solution to the PDE.

The results show significant gains in accuracy
compared with the conventional DG method,
particularly for low-dimensional approximation spaces.
To obtain the same accuracy, we can significantly reduce
the number of cells and use larger time steps.
The method has been validated on a wide range of PDEs and equilibria,
showing that it is a general-purpose approach.
Furthermore, it makes it possible to handle complicated equilibria,
on complex geometries,
which are rarely treated by conventional WB schemes,
especially in two space dimensions.
The cost of training the network is low,
as is the cost of inference.
The main additional cost of the method comes from
the quadrature rule, whose order has to be increased to
ensure a good approximation of the integral of the prior.
In most cases, this increase in order is not very important,
and the gain between our approach
and the classical ones remains significant.

There are several possible ways of extending our approach.
From an application point of view,
we wish to deal with more difficult equilibria,
such as equilibria for the magnetohydrodynamics in tokamaks.
From a methodological point of view,
{\rab we could implement an orthogonalization process
		to ensure that the DG mass matrix is diagonal instead of block diagonal,
		which would improve the computation time.
		Moreover, we would like to improve the determination of the prior}
by replacing parametric PINNs with
physics-informed neural
operators~\cite{wang2021learning,goswami2022physics}
in order to widen the family of equilibria that can be considered.
The other approach is to extend the method with time-dependent
priors, in order to increase the accuracy of the scheme
around families of unsteady solutions.
To that end, we wish to move on to space-time DG methods,
see e.g. \cite{PetFarTez2008}.

\bibliographystyle{plain}
\bibliography{biblio}

\end{document}